\documentclass[11pt,twoside]{article} %
\usepackage{amsmath,amssymb}

\usepackage{hyperref}
\hypersetup{
    hypertexnames=false,
    colorlinks,
    citecolor=black,
    filecolor=black,
    linkcolor=black,
    urlcolor=black}
\makeatletter
\newcommand*{\mint}[1]{%
  \mint@l{#1}{}%
}
\newcommand*{\mint@l}[2]{%
  \@ifnextchar\limits{%
    \mint@l{#1}%
  }{%
    \@ifnextchar\nolimits{%
      \mint@l{#1}%
    }{%
      \@ifnextchar\displaylimits{%
        \mint@l{#1}%
      }{%
        \mint@s{#2}{#1}%
      }%
    }%
  }%
}
\newcommand*{\mint@s}[2]{%
  \@ifnextchar_{%
    \mint@sub{#1}{#2}%
  }{%
    \@ifnextchar^{%
      \mint@sup{#1}{#2}%
    }{%
      \mint@{#1}{#2}{}{}%
    }%
  }%
}
\def\mint@sub#1#2_#3{%
  \@ifnextchar^{%
    \mint@sub@sup{#1}{#2}{#3}%
  }{%
    \mint@{#1}{#2}{#3}{}%
  }%
}
\def\mint@sup#1#2^#3{%
  \@ifnextchar_{%
    \mint@sup@sub{#1}{#2}{#3}%
  }{%
    \mint@{#1}{#2}{}{#3}%
  }%
}
\def\mint@sub@sup#1#2#3^#4{%
  \mint@{#1}{#2}{#3}{#4}%
}
\def\mint@sup@sub#1#2#3_#4{%
  \mint@{#1}{#2}{#4}{#3}%
}
\newcommand*{\mint@}[4]{%
  \mathop{}%
  \mkern-\thinmuskip
  \mathchoice{%
    \mint@@{#1}{#2}{#3}{#4}%
        \displaystyle\textstyle\scriptstyle
  }{%
    \mint@@{#1}{#2}{#3}{#4}%
        \textstyle\scriptstyle\scriptstyle
  }{%
    \mint@@{#1}{#2}{#3}{#4}%
        \scriptstyle\scriptscriptstyle\scriptscriptstyle
  }{%
    \mint@@{#1}{#2}{#3}{#4}%
        \scriptscriptstyle\scriptscriptstyle\scriptscriptstyle
  }%
  \mkern-\thinmuskip
  \int#1%
  \ifx\\#3\\\else_{#3}\fi
  \ifx\\#4\\\else^{#4}\fi
}
\newcommand*{\mint@@}[7]{%
  \begingroup
    \sbox0{$#5\int\m@th$}%
    \sbox2{$#5\int_{}\m@th$}%
    \dimen2=\wd0 %
    \let\mint@limits=#1\relax
    \ifx\mint@limits\relax
      \sbox4{$#5\int_{\kern1sp}^{\kern1sp}\m@th$}%
      \ifdim\wd4>\wd2 %
        \let\mint@limits=\nolimits
      \else
        \let\mint@limits=\limits
      \fi
    \fi
    \ifx\mint@limits\displaylimits
      \ifx#5\displaystyle
        \let\mint@limits=\limits
      \fi
    \fi
    \ifx\mint@limits\limits
      \sbox0{$#7#3\m@th$}%
      \sbox2{$#7#4\m@th$}%
      \ifdim\wd0>\dimen2 %
        \dimen2=\wd0 %
      \fi
      \ifdim\wd2>\dimen2 %
        \dimen2=\wd2 %
      \fi
    \fi
    \rlap{%
      $#5%
        \vcenter{%
          \hbox to\dimen2{%
            \hss
            $#6{#2}\m@th$%
            \hss
          }%
        }%
      $%
    }%
  \endgroup
}
\usepackage{mathrsfs}
\usepackage{amssymb}
\usepackage{amsmath}
\usepackage{amsthm}
\usepackage{needspace}
\usepackage{amsfonts}
\usepackage{color}
\usepackage{graphicx}
\usepackage[active]{srcltx}
\usepackage{tikz}
\usepackage{pgflibraryarrows}
\usepackage{pgflibrarysnakes}

\usepackage[cp1252]{inputenc}

\usepackage{mathrsfs}
\usepackage{graphicx}

\usepackage[active]{srcltx}

\allowdisplaybreaks

\usepackage{titletoc}
\titlecontents{section}[0pt]{\addvspace{2pt}\filright}
              {\contentspush{\thecontentslabel\ }}
              {}{\titlerule*[8pt]{.}\contentspage}

\def\bint{{\ifinner\rlap{\bf\kern.35em--}
\int\else\rlap{\bf\kern.45em--}\int\fi}\ignorespaces}

\def\bbint{{\ifinner\rlap{\bf\kern.35em--}
\hspace{0.078cm}\int\else\rlap{\bf\kern.45em--}\int\fi}\ignorespaces}

\def\dfrac{\displaystyle\frac}

\newtheorem{thm}{Theorem}[section]
\newtheorem{lem}[thm]{Lemma}%
\theoremstyle{definition}
\theoremstyle{plain}
\newtheorem{cor}[thm]{Corollary}%
\numberwithin{equation}{section}

\title{Stable solutions to $p$-Laplacian type equations\\
with variable coefficients\thanks{%
Weican Luo was supported by the National Key R\&D Program of China
(No.~2025YFA1018400) and the National Natural Science Foundation of China
(No.~12571209).}
}

\author{%
Weican Luo\\[.4em]
{\small\itshape School of Mathematical Sciences, Beihang University, Changping District,}\\
{\small\itshape Shahe Higher Education Park, South Third Street No.~9, Beijing 102206, PR China}
}
\date{}

\makeatletter
\renewcommand{\@maketitle}{%
  \newpage
  \null
  \vskip .5em
  \begin{center}%
    \let\footnote\thanks
    {\Large\bfseries\boldmath \@title\par}%
    \ifx\@author\@empty
    \else
      \vskip .6em
      {\large\begin{tabular}[t]{c}\@author\end{tabular}\par}%
    \fi
  \end{center}%
  \vskip .5em
}
\makeatother

\providecommand{\Ball}{B}
\providecommand{\dd}{\,dx}

\providecommand{\R}{\mathbb R}
\providecommand{\op}{\mathrm{op}}

\providecommand{\norm}[1]{\left\|#1\right\|}

\newtheorem{lemma}[thm]{Lemma}
\newtheorem{theorem}[thm]{Theorem}

\theoremstyle{definition}
\newtheorem{remark}[thm]{Remark}
\theoremstyle{plain}

\renewenvironment{abstract}
  {\small\quotation\noindent\textbf{\abstractname.}\enspace\ignorespaces}
  {\endquotation}

\begin{document}

\arraycolsep=1pt
\maketitle
\begingroup
\renewcommand{\thefootnote}{}
\footnotetext[0]{E-mail address:
\href{mailto:weicanluo@buaa.edu.cn}{\nolinkurl{weicanluo@buaa.edu.cn}}.}
\endgroup

\begin{abstract}
We study stable regular solutions of the equation
\[
-\operatorname{div}\bigg(\Big(A(x)\nabla u \cdot \nabla u\Big)^{p/2-1}A(x)\nabla u\bigg) = f(u) \qquad \text{in } B_1 \subset \mathbb{R}^n,
\]
where \(p > 1\), \(f \in C^1(\mathbb{R})\) is nonnegative, and \(A\) is a locally Lipschitz, symmetric, uniformly elliptic matrix field.
We prove that stable regular solutions satisfy an interior H\"older estimate for \(n < p + 4p/(p - 1)\) when \(p \ge 2\), and for \(n < 5p\) when \(1 < p < 2\).
This extends the result of Cabr\'e, Miraglio, and Sanch\'on \cite{cms} for the \(p\)-Laplace equation.

Note that the constant in our H\"older estimate depends on the local \(C^1\) norm of \(A\), but not on its \(C^2\) norm.
Our approach is motivated by the interior regularity estimates of Cabr\'e, Miraglio, and Sanch\'on \cite{cms} for the \(p\)-Laplace equation and by Erneta's \cite{erneta} treatment of variable coefficients in semilinear elliptic equations.
\end{abstract}

\smallskip
{\small\raggedright
\par
}
\medskip

\begingroup
\small
\setcounter{tocdepth}{1}
\tableofcontents
\endgroup

\section{Introduction}

Let $B_1\subset\R^n$ be the unit ball and let $p>1$. We consider the equation
\begin{equation}
-\operatorname{div}\!\left(
    \bigl(A(x)\nabla u\cdot\nabla u\bigr)^{\frac p2-1}
    A(x)\nabla u
\right)=f(u)
\qquad\text{in }B_1.
\label{intro:pde}
\end{equation}
Here $A\in W^{1,\infty}_{\mathrm{loc}}(B_1;\R^{n\times n})$ is real symmetric. We use its locally Lipschitz representative and assume uniform ellipticity:
\[
\lambda |\xi|^2
\le A(x)\xi\cdot\xi
\le \Lambda |\xi|^2,
\qquad x\in B_1,\ \xi\in\R^n,
\]
where $0<\lambda\le\Lambda$. At $\nabla u=0$, the flux in \eqref{intro:pde} is defined by its continuous extension, equal to zero. If $A\equiv I$, equation \eqref{intro:pde} reduces to the usual $p$-Laplace equation
\[
-\Delta_p u
=-\operatorname{div}\bigl(|\nabla u|^{p-2}\nabla u\bigr)
=f(u),
\]
while for $p=2$ it becomes the variable-coefficient semilinear equation in divergence form
\[
-\operatorname{div}\bigl(A(x)\nabla u\bigr)=f(u).
\]

For $f\in C^1(\R)$ and $F'=f$, equation \eqref{intro:pde} has the variational functional
\begin{equation}
\mathcal E(v)
=
\int_{B_1}
\left\{
\frac1p\bigl(A(x)\nabla v\cdot\nabla v\bigr)^{p/2}-F(v)
\right\}\,dx.
\label{intro:energy}
\end{equation}
We consider this functional on $v\in W^{1,p}(B_1)$ with $F(v)\in L^1(B_1)$, so that $\mathcal E(v)$ is finite.
The Euler--Lagrange equation of \eqref{intro:energy} is \eqref{intro:pde}. Its second variation leads to a stability inequality whose coefficients may degenerate or become singular on the critical set when $p\ne2$. We use the weighted stability framework of Cabr\'e, Miraglio, and Sanch\'on~\cite{cms}, adapted to the anisotropic principal part.

\subsection{Problem setting, stability, and previous results}
\label{intro:solution-class}

We call $u\in W^{1,p}(B_1)$ an \textbf{energy solution} of \eqref{intro:pde} if $f(u)\in L^1_{\mathrm{loc}}(B_1)$ and
\begin{equation}
\int_{B_1}
\bigl(A(x)\nabla u\cdot\nabla u\bigr)^{\frac p2-1}
A(x)\nabla u\cdot\nabla\varphi\,dx
=
\int_{B_1}f(u)\varphi\,dx
\label{intro:weak-form}
\end{equation}
for every $\varphi\in C_c^1(B_1)$. If, in addition, $f(u)\in L^\infty(B_1)$, we call $u$ a \textbf{regular solution}.
Under the above coefficient assumptions, every regular solution has a
representative in $C^{1,\vartheta}_{\mathrm{loc}}(B_1)$ for some
$\vartheta\in(0,1)$. Indeed, the flux in \eqref{intro:pde} satisfies the
standard $p$-growth and ellipticity conditions and is locally Lipschitz
in $x$, with a Lipschitz bound proportional to $|\xi|^{p-1}$,
where $\xi$ denotes the gradient variable. For $n\ge2$, local application of
\cite[Theorem~1.1]{araujo-zhang} when $p\ge2$, and
\cite[Corollary~1.7]{dong-zhu} when $1<p<2$, gives this regularity.
For the latter application, the Lipschitz coefficient modulus is bounded
by $Cr^\vartheta$, and the bounded right-hand side satisfies
\[
\int_{B_r(x)}|f(u)|\,dy\le C r^n\le C r^{n-1+\vartheta}
\qquad(0<r\le1,\ B_r(x)\Subset B_1),
\]
with $\vartheta>0$ chosen sufficiently small as required by that corollary.
When $n=1$, the same conclusion follows by integrating the flux and
inverting the map $t\mapsto |t|^{p-2}t$.
We use this $C^1$ representative throughout. This regularity is a
qualitative input; its solution-dependent bounds do not enter the final
quantitative constants. On compact subsets of the noncritical region,
the linearized principal part is uniformly elliptic. Appendix~A supplies
the weak second-order regularity and critical-set cutoff justification
used below.
Classical local gradient regularity is developed in~\cite{uhlenbeck,tolksdorf}.

To define stability, set
\[
Z_u:=\{x\in B_1:\nabla u(x)=0\}.
\]
On $B_1\setminus Z_u$, let
\[
\sigma_u(x):=
\bigl(A(x)\nabla u\cdot\nabla u\bigr)^{\frac{p-2}{2}}.
\]
On $Z_u$ we use the extended-value convention
\[
\sigma_u=0 \quad\text{if }p>2,
\qquad
\sigma_u=1 \quad\text{if }p=2,
\qquad
\sigma_u=+\infty \quad\text{if }1<p<2.
\]
For $p\ge2$, define $W^{1,2}_{\sigma_u}(B_1)$ as the completion of those functions in $C^1(B_1)$ for which the following norm is finite:
\[
\norm{\xi}_{W^{1,2}_{\sigma_u}(B_1)}
:=
\norm{\xi}_{L^2(B_1)}
+
\left(
\int_{B_1}
\sigma_u\,A(x)\nabla\xi\cdot\nabla\xi\,dx
\right)^{1/2},
\]
and let $W^{1,2}_{\sigma_u,0}(B_1)$ denote the closure of $C_c^1(B_1)$ in this norm. Following~\cite{cms}, we take the stability test space to be
\begin{equation}
\mathcal T_u:=
\begin{cases}
W^{1,2}_{\sigma_u,0}(B_1), & p\ge2,\\[2mm]
\displaystyle
\left\{
\xi\in W^{1,2}_0(B_1):
\int_{B_1}\sigma_u\,A(x)\nabla\xi\cdot\nabla\xi\,dx<\infty
\right\},
& 1<p<2.
\end{cases}
\label{intro:test-space}
\end{equation}
Uniform ellipticity makes $\sigma_u$ comparable to $|\nabla u|^{p-2}$ off $Z_u$. For $1<p<2$, the weighted density on $Z_u$ is taken to be zero if $\nabla\xi=0$ and $+\infty$ otherwise; thus finite energy forces $\nabla\xi=0$ almost everywhere on $Z_u$. Expressions with negative powers are evaluated off $Z_u$. The gradient-dependent stability tests used below are justified by the critical-set cutoffs and weighted second-order estimates of Appendix~A.

A regular solution $u$ of \eqref{intro:pde} is called a \textbf{stable regular solution} if, for every $\xi\in\mathcal T_u$,
\begin{equation}
\begin{aligned}
&\int_{B_1}
\bigl(A(x)\nabla u\cdot\nabla u\bigr)^{\frac p2-1}
A(x)\nabla\xi\cdot\nabla\xi\,dx
\\
&\quad +(p-2)
\int_{B_1\setminus Z_u}
\bigl(A(x)\nabla u\cdot\nabla u\bigr)^{\frac p2-2}
\bigl(A(x)\nabla u\cdot\nabla\xi\bigr)^2\,dx
\\
&\qquad -\int_{B_1}f'(u)\xi^2\,dx\ge0.
\end{aligned}
\label{intro:stability}
\end{equation}
When $1<p<2$, the first term is understood in the weighted sense described above. Inequality \eqref{intro:stability} expresses nonnegativity of the weighted quadratic form associated with the second variation of \eqref{intro:energy} at $u$. By the Cauchy--Schwarz inequality for the inner product induced by~$A$,\looseness=-1\relax
\[
\bigl(A\nabla u\cdot\nabla\xi\bigr)^2
\le
\bigl(A\nabla u\cdot\nabla u\bigr)
\bigl(A\nabla\xi\cdot\nabla\xi\bigr).
\]
Hence, when $1<p<2$,
\[
\begin{aligned}
&\bigl(A\nabla u\cdot\nabla u\bigr)^{\frac p2-1}
A\nabla\xi\cdot\nabla\xi
\\
&\quad +(p-2)
\bigl(A\nabla u\cdot\nabla u\bigr)^{\frac p2-2}
\bigl(A\nabla u\cdot\nabla\xi\bigr)^2
\\
&\ge
(p-1)
\bigl(A\nabla u\cdot\nabla u\bigr)^{\frac p2-1}
A\nabla\xi\cdot\nabla\xi.
\end{aligned}
\]
Thus the principal part of the stability quadratic form is still nonnegative when $1<p<2$.

For $p=2$ and $A\equiv I$, \eqref{intro:stability} reduces to the classical inequality $\int_{B_1}f'(u)\xi^2\,dx\le\int_{B_1}|\nabla\xi|^2\,dx$.

{\emergencystretch=1em
The study of stable semilinear equations includes the nonlinear eigenvalue problems of Crandall--Rabinowitz~\cite{crandall-rabinowitz} and the low-dimensional regularity results of Nedev~\cite{nedev} and Cabr\'e~\cite{cabre-low-dim} for Dirichlet problems. Cabr\'e, Figalli, Ros-Oton, and Serra~\cite{cfrs} proved interior quantitative H\"older estimates for stable classical solutions with nonnegative nonlinearities in dimensions $n\le9$, with constants independent of the nonlinearity. The dimensional threshold is sharp, since singular stable energy solutions exist for $n\ge10$; see~\cite[Section~5]{brezis-vazquez}.\par}

For the $p$-Laplace equation, Cabr\'e--Sanch\'on~\cite{cabre-sanchon} studied semistable and extremal solutions under growth assumptions on the nonlinearity. Cabr\'e--Capella--Sanch\'on~\cite{cabre-capella-sanchon} obtained nonlinearity-independent pointwise estimates for radially decreasing semistable solutions in $W^{1,p}(B_1)$, with boundedness in the optimal radial range $n<p+4p/(p-1)$. Castorina--Sanch\'on~\cite{castorina-sanchon} used geometric Sobolev and Morrey inequalities to obtain a priori estimates and regularity of extremal solutions in smooth bounded strictly convex domains, under additional restrictions on $p$ and the nonlinearity. Farina--Sciunzi--Valdinoci~\cite{farina-sciunzi-valdinoci} developed a geometric approach to symmetry for stable solutions of autonomous degenerate or singular elliptic equations. For general nonradial solutions, Cabr\'e, Miraglio, and Sanch\'on~\cite{cms} established the interior estimates recalled below, with different dimensional restrictions for $p\ge2$ and $1<p<2$.

The two results most directly related to \eqref{intro:pde} concern the constant-coefficient $p$-Laplace equation and variable-coefficient semilinear equations.

\medskip
\noindent\textbf{Theorem (Cabr\'e--Miraglio--Sanch\'on~\cite{cms}).}
\begingroup\itshape
Let $u$ be a stable regular solution of $-\Delta_pu=f(u)$ in $B_1$, where $p>1$, $f\in C^1(\R)$, and $f\ge0$. Then there exist $\gamma>0$ and $C>0$, depending only on $n$ and $p$, such that
\[
\norm{\nabla u}_{L^{p+\gamma}(B_{1/2})}
\le
C\norm{u}_{L^1(B_1)},
\]
and
\[
\norm{\Delta_{p+1}u}_{L^1(B_{1/2})}
\le
C\norm{u}_{L^1(B_1)}^p.
\]
If, in addition,
\begin{equation}
\begin{cases}
n<p+\dfrac{4p}{p-1}, & p\ge2,\\[2mm]
n<5p, & 1<p<2,
\end{cases}
\label{intro:cms-dimension}
\end{equation}
then there exist $\alpha>0$ and $C>0$, depending only on $n$ and $p$, such that
\begin{equation}
\norm{u}_{C^\alpha(B_{1/2})}
\le
C\norm{u}_{L^1(B_1)}.
\label{intro:cms-holder}
\end{equation}
\par\endgroup
For $p\ge2$, the dimensional range in \eqref{intro:cms-dimension} is optimal; for $1<p<2$, the condition $n<5p$ is not sharp. The gradient and $(p+1)$-Laplacian estimates hold in every dimension. The latter has the $p$-homogeneous power of $\norm{u}_{L^1}$, obtained by rescaling the normalized estimate in~\cite[Lemma~3.2]{cms} and using the energy bound in the proof of its Theorem~1.1.

Erneta~\cite{erneta} treats the nondivergence-form operator $Lu=a_{ij}(x)u_{ij}+b_i(x)u_i$, with smooth coefficients and a real symmetric uniformly elliptic principal matrix. Its divergence-form representation has drift $\widehat b_i=b_i-\partial_k a_{ki}$, so it need not have a variational structure. Erneta defines stability using a supersolution of the linearized equation in $C^2(B_1)\cap C^0(\overline{B_1})$, positive in $B_1$ and zero on $\partial B_1$.

\medskip
\noindent\textbf{Theorem (Erneta \cite{erneta}).}
\begingroup\itshape
Let $u\in C^\infty(\overline{B_1})$ be a stable solution of
$-Lu=f(u)$, where $f\in C^1(\mathbb R)$ and $f\ge0$. Then there exist
$\gamma=\gamma(n)>0$ and
\[
C=C\bigl(n,\lambda,\|A\|_{C^1(B_1)},\|b\|_{C^0(B_1)}\bigr)>0
\]
such that
\[
\|\nabla u\|_{L^{2+\gamma}(B_{1/2})}
\le C\|u\|_{L^1(B_1)}.
\]
If, in addition, $n\le9$, then there exists
$\alpha=\alpha(n,\lambda,\Lambda)>0$ such that
\begin{equation}
\|u\|_{C^\alpha(B_{1/2})}
\le C\|u\|_{L^1(B_1)},
\label{intro:erneta-holder}
\end{equation}
where $C$ has the same dependence as above. In particular, neither $D^2A$ nor $Db$ enters the constants.
\par\endgroup

Equation \eqref{intro:pde} combines the variable coefficients in~\cite{erneta} with the $p$-growth structure in~\cite{cms}. Spatial differentiation introduces errors involving $DA$, while $p\ne2$ also produces degeneracy or singularity at critical points. We establish quantitative interior estimates with constants depending only on $n,p$, ellipticity, and first-order coefficient control, independently of $f$, $u$, and higher derivatives of $A$.

Related regularity theories address different estimates and assumptions. Weighted second-order estimates for positive solutions of autonomous $p$-Laplace equations were obtained by Damascelli--Sciunzi~\cite{damascelli-sciunzi}. For $L^2$ data, Miao--Peng--Zhou~\cite{miao-peng-zhou} study the same principal part as \eqref{intro:pde} and obtain global $W^{1,2}$ regularity of the nonlinear stress under suitable coefficient, domain, and boundary assumptions. For positive solutions of autonomous anisotropic equations satisfying their structural and nonlinearity hypotheses, Castorina--Riey--Sciunzi~\cite{castorina-riey-sciunzi} establish weighted regularity estimates and prove a Hopf lemma with zero Dirichlet data and an interior sphere condition. Antonini--Ciraolo--Farina~\cite{antonini-ciraolo-farina} prove local stress-field Sobolev regularity and weighted Hessian estimates under suitable integrability assumptions on the source. Our focus is on interior regularity estimates supplied by stability in terms of $\norm{u}_{L^1(B_1)}$, without a norm of $f(u)$ on the right-hand side.

\subsection{Main results}

Theorem~1.1 gives the interior quantitative regularity estimates. Theorem~1.2 provides the Hessian bounds used to control coefficient errors and obtain compactness.

\medskip
\noindent\textbf{Theorem 1.1 (Interior quantitative regularity).}
\begingroup\itshape
Let $p>1$ and let $u$ be a stable regular solution of \eqref{intro:pde}, where $f\in C^1(\R)$ and $f\ge0$. Assume that $A\in W^{1,\infty}_{\mathrm{loc}}(B_1;\R^{n\times n})$ is real symmetric and uniformly elliptic, and set
\[
L:=\norm{DA}_{L^\infty(B_{3/4})}.
\]
There exist
\[
\gamma_0=\gamma_0(n,p,\lambda,\Lambda)>0,
\qquad C=C(n,p,\lambda,\Lambda,L)>0,
\]
such that, in every dimension,
\begin{equation}
\left(
\int_{B_{1/2}}
\bigl(A(x)\nabla u\cdot\nabla u\bigr)^{p/2}\,dx
\right)^{1/p}
+
\norm{\nabla u}_{L^{p+\gamma_0}(B_{1/2})}
\le
C\norm{u}_{L^1(B_1)}.
\label{intro:main-energy}
\end{equation}
If, in addition,
\begin{equation}
n<p+\frac{4p}{p-1}
\quad\text{if }p\ge2,
\qquad
n<5p
\quad\text{if }1<p<2,
\label{intro:main-dimension}
\end{equation}
then there exists $\alpha=\alpha(n,p,\lambda,\Lambda)\in(0,1)$ such that,
after enlarging $C$ with the same dependence if necessary, for every
$y\in B_{1/8}$ and $0<r\le1/16$,
\begin{equation}
\int_{B_r(y)}|\nabla u|^p\,dx
\le
Cr^{\,n-p+p\alpha}
\norm{u}_{L^1(B_1)}^p.
\label{intro:main-morrey}
\end{equation}
Therefore,
\begin{equation}
\norm{u}_{C^\alpha(B_{1/16})}
\le
C\norm{u}_{L^1(B_1)}.
\label{intro:main-holder}
\end{equation}
All exponents and constants above are independent of $f$, $u$, and
derivatives of $A$ of order two or higher.
\par\endgroup

To handle the error terms in the weighted stability calculation, we also need the following second-order estimates.

\medskip
\noindent\textbf{Theorem 1.2 (A priori Hessian estimates).}
\begingroup\itshape
Let $p>1$, and let $u$ be a stable regular solution of \eqref{intro:pde}, where $f\in C^1(\R)$, $f\ge0$, and $A\in W^{1,\infty}_{\mathrm{loc}}(B_1;\R^{n\times n})$ is real symmetric and satisfies
\[
\lambda|\xi|^2
\le
A(x)\xi\cdot\xi
\le
\Lambda|\xi|^2,
\qquad
x\in B_1,\ \xi\in\R^n.
\]
The dimensional restriction for the H\"older estimate in Theorem~1.1 is not needed here. There exists
\[
\delta_0=\delta_0(n,p,\lambda,\Lambda)>0
\]
such that, whenever
\[
\norm{DA}_{L^\infty(B_1)}\le\delta_0,
\]
the estimates below hold with $C=C(n,p,\lambda,\Lambda)$. If $p\ge2$, then
\begin{equation}
\int_{B_{1/2}\cap\{|\nabla u|>0\}}
\bigl(A(x)\nabla u\cdot\nabla u\bigr)^{(p-2)/2}
|D^2u|\,dx
\le
C
\left(
\int_{B_1}
\bigl(A(x)\nabla u\cdot\nabla u\bigr)^{p/2}\,dx
\right)^{(p-1)/p}.
\label{intro:hessian-ge2}
\end{equation}
If $1<p<2$, then
\begin{equation}
\int_{B_{1/2}\cap\{|\nabla u|>0\}}
|D^2u|\,dx
\le
C
\left(
\int_{B_1}
\bigl(A(x)\nabla u\cdot\nabla u\bigr)^{p/2}\,dx
\right)^{1/p}.
\label{intro:hessian-lt2}
\end{equation}
In addition, if $B_{2r}(z)\Subset B_1$ and $0<r\le1$, then
\begin{equation}
\int_{B_r(z)\cap\{|\nabla u|>0\}}
\bigl(A(x)\nabla u\cdot\nabla u\bigr)^{(p-1)/2}
|D^2u|\,dx
\le
\frac Cr
\int_{B_{2r}(z)}
\bigl(A(x)\nabla u\cdot\nabla u\bigr)^{p/2}\,dx.
\label{intro:hessian-local}
\end{equation}
\par\endgroup
The smallness condition can be achieved by rescaling on any fixed interior subdomain, with a scale depending on the structural parameters, a first-order coefficient bound on a larger interior subdomain, and the distance from the smaller subdomain to the boundary of the larger one. Estimate \eqref{intro:hessian-local} controls the Hessian errors in the weighted radial calculation. When $A\equiv I$, the H\"older dimension range in Theorem~1.1 agrees with \eqref{intro:cms-dimension}.

\subsection{Main ideas of the proof}

Theorem~1.1 extends the interior estimates of Cabr\'e, Miraglio, and
Sanch\'on~\cite{cms} to locally Lipschitz coefficient fields $A(x)$.
For the Hessian estimate in Section~3, we adapt the geometric argument of
Sternberg and Zumbrun in the quasilinear form of Farina, Sciunzi, and
Valdinoci~\cite[proof of Theorem~2.5]{farina-sciunzi-valdinoci}.
For $\eta\in C_c^\infty(B_1)$ and $\tau>0$, we test the stability
inequality with
\[
 \psi_\tau(x)=|\nabla u(x)|\,\eta(x)\,
 \chi\!\left(\frac{|\nabla u(x)|}{\tau}\right).
\]
Here $\chi:[0,\infty)\to[0,1]$ is smooth and nondecreasing, with
$\chi=0$ on $[0,1]$ and $\chi=1$ on $[2,\infty)$.
Comparing the stability inequality with the summed linearized identities
cancels the $f'(u)$ terms. The preliminary weighted Hessian bound allows
us to let $\tau\downarrow0$. Differentiating the flux
in \eqref{flux-differentiation} introduces an error involving $DA$.
Since $f\ge0$, the equation bounds the positive part of the second
derivative normal to the level sets of $u$. Integration by parts bounds its
signed integral, allowing us to absorb the coefficient error. Localization gives
\eqref{local-pplus1-hessian}.
Higher integrability, interpolation and Simon's absorption
lemma~\cite{simon-scaling}, in the form of~\cite[Lemma~A.4]{cms},
then give \eqref{intro:main-energy}. For $n\le p$, we have $p+\gamma_0>n$, so Morrey's
embedding gives the H\"older estimate and H\"older's inequality gives the local energy decay.

For $n>p$ in the range \eqref{intro:main-dimension}, fix normalized
coordinates as in \eqref{s5:normalized-setting} and set $r=|x|$,
$u_r=(x/|x|)\cdot\nabla u$. The test based on $x\cdot\nabla u=ru_r$
was used by Cabr\'e, Figalli, Ros-Oton, and Serra~\cite{cfrs} for $p=2$.
Cabr\'e, Miraglio, and Sanch\'on~\cite[(1.22)]{cms} added a cutoff
near $\{\nabla u=0\}$ for the $p$-Laplacian. To prove Theorem~\ref{raw},
we test the stability inequality with
\[
 \xi_\tau=(x\cdot\nabla u)\eta\phi_\tau,\qquad
 \phi_\tau(x)=\phi\left(\frac{|\nabla u(x)|}{\tau}\right).
\]
Here $\phi\in C^\infty(\mathbb R)$ satisfies $0\le\phi\le1$,
$\phi=0$ on $(-\infty,1]$ and $\phi=1$ on $[2,\infty)$;
$\zeta$ is the spatial cutoff in Theorem~\ref{raw}.
For $\eta=r^{-a/2}\zeta$, we first use
$\eta_\delta=w_\delta\zeta$, $w_\delta(r)=\min\{r^{-a/2},\delta^{-a/2}\}$,
and let $\tau\downarrow0$ before $\delta\downarrow0$.
Following Erneta's dyadic argument for semilinear equations~\cite{erneta},
we sum \eqref{local-pplus1-hessian} on annuli. By
\eqref{s5:weighted-hessian-margin}, the mixed error involving $DA$ and $D^2u$
is bounded by the gradient energy with weight $r^{1-a}$, an extra factor
of $r$ over the weight $r^{-a}$.

Following Erneta~\cite{erneta}, we first use a less singular weight,
taking $\theta=\tfrac12\min\{1,n-p\}$ and $a=n-p-\theta$. The full energy
on $B_\rho$ has coefficient $n-p-a=\theta>0$ and absorbs the inner error
at sufficiently small scales. Dropping the nonnegative radial term gives
\eqref{s5:subcritical},
\[
 \int_{B_\rho} r^{p+\theta-n}(A\nabla u\cdot\nabla u)^{p/2}\,dx
 \le C_\theta\int_{B_{4\rho}\setminus B_\rho}
              r^{p+\theta-n}(A\nabla u\cdot\nabla u)^{p/2}\,dx .
\]
At the critical exponent $a=n-p$, the coefficient of the full gradient
energy vanishes. For $n>p$, the radial coefficient $c_p(n-p)$ in
\eqref{s5:reduced-weighted} is positive precisely in the dimension range
\eqref{intro:main-dimension}. We use the
subcritical estimate to move the inner error to the annulus, obtaining
\[
\begin{aligned}
 \int_{B_\rho}r^{p+1-n}(A\nabla u\cdot\nabla u)^{p/2}\,dx
 &\le\rho^{1-\theta}
       \int_{B_\rho}r^{p+\theta-n}(A\nabla u\cdot\nabla u)^{p/2}\,dx\\
 &\le C_\theta\rho
       \int_{B_{4\rho}\setminus B_\rho}
       r^{p-n}(A\nabla u\cdot\nabla u)^{p/2}\,dx .
\end{aligned}
\]
This yields the critical radial estimate \eqref{s5:critical-discrete-input}.
For the converse bound, we adapt the compactness argument
in~\cite[Lemma~4.1]{cms} to varying coefficient fields. Under the doubling
hypothesis, Lemma~\ref{thm:anisotropic-lemma41} controls the full energy
in the ball by the radial energy in its outer annulus. If the comparison
failed, normalization and doubling would preserve positive energy in $B_1$
under local strong $L^p$ convergence of the gradients. Since the limiting
supersolution has zero radial derivative a.e. in $B_{3/2}\setminus\overline{B_1}$,
comparison, the weak Harnack inequality and the strong minimum principle
force it to be constant a.e. in $B_{3/2}$, a contradiction.

The two-sequence lemma of~\cite[Lemma~3.2]{cfrs} was used for the
$p$-Laplacian in~\cite[Lemma~5.1]{cms}. Here the recurrences
\eqref{s5:discrete-radial-bound}--\eqref{s5:discrete-recurrence}
yield a contraction for a positive linear combination of the full and
radial energies, giving power decay of the full gradient energy.
Since affine changes generally send balls to ellipsoids and alter the
radial weight, we first establish this decay at a fixed center and then
transfer it to other centers. Campanato's characterization gives the
H\"older estimate: the exponent depends only on $n,p,\lambda,\Lambda$,
while the constant also depends on $\|DA\|_{L^\infty(B_{3/4})}$.

\subsection{Organization of the paper}

Section~2 derives the weighted radial estimate in normalized coordinates, retaining the Hessian errors from the variable coefficients. Section~3 proves Theorem~1.2 and higher integrability of the gradient. Section~4 compares radial and full gradient energies. Section~5 combines these estimates with local affine normalization and scale iteration to prove Theorem~1.1. Appendix~A supplies the difference-quotient regularity needed for weak linearization and critical-set cutoffs.

\section{A weighted energy estimate}

In this section we derive a weighted energy inequality for stable
regular solutions, with explicit error terms arising from the variable
coefficients. The proof uses the stability inequality with a radially
weighted cutoff of $x\cdot\nabla u$. This inequality will be used in
Section~5 to obtain the radial energy estimate.

Let $1<p<n$ and let $u$ be a stable regular solution
of (1.1), with $f\in C^1(\mathbb R)$. We use the solution class,
stability test space $\mathcal T_u$, and critical-set conventions
of Section~1.1; no sign condition on $f$ is imposed in this section.
Thus stability is understood as in (1.4)--(1.5).

Let $A\in W^{1,\infty}_{\rm loc}(B_1;\mathbb R^{n\times n})$
be real symmetric, and use its locally Lipschitz representative. Assume
\begin{equation}
 \lambda|\xi|^2\le A(x)\xi\cdot\xi\le\Lambda|\xi|^2,
 \qquad x\in B_1,\quad \xi\in\mathbb R^n.
\end{equation}
We work in fixed normalized coordinates centered at the selected point,
with
\begin{equation}
 A(0)=I,\qquad \|DA\|_{L^\infty(B_1)}\le\varepsilon,
 \qquad \varepsilon\ge0.
 \label{normass}
\end{equation}
For coefficient matrices and their directional derivatives, we use the operator
norm and set $|DA|=\sup_{|e|=1}|\partial_eA|_{\rm op}$.
Consequently,
\[
 |A(x)-I|\le\varepsilon|x|,\qquad
 |x_kA_{ij,k}(x)|\le\varepsilon|x|,
\]
where the derivative inequality holds almost everywhere.
All radial quantities refer to these same coordinates:
\begin{equation}
 r=|x|,\qquad u_r=\frac{x}{|x|}\cdot\nabla u,
 \qquad x\cdot\nabla u=ru_r.
\end{equation}
Repeated indices are summed. Such coordinates are obtained by translation
and affine normalization, followed by a spatial rescaling when needed.

\begin{theorem}\label{raw}
Let $0<\rho\le 1/4$, let $0<a<n$, and let $\zeta\in C_c^1(B_{4\rho})$ satisfy
\begin{equation}
0\le\zeta\le1,
\qquad \zeta\equiv1\text{ in }B_\rho,
\qquad |\nabla\zeta|\le \frac{C}{\rho}.
\end{equation}
If $p\ge2$ and
\begin{equation}
pa-\frac{(p-1)a^2}{4}>0,
\end{equation}
\begingroup
\allowdisplaybreaks[0]
then
\begin{align}
&\left(pa-\frac{(p-1)a^2}{4}\right)
\int_{B_\rho}r^{-a}(A\nabla u\cdot\nabla u)^{p/2-1}u_r^2\,dx
+(n-p-a)\int_{B_\rho}r^{-a}(A\nabla u\cdot\nabla u)^{p/2}\,dx
\notag\\
&\le C\int_{B_{4\rho}\setminus B_\rho}r^{-a}(A\nabla u\cdot\nabla u)^{p/2}\,dx
+C\varepsilon\int_{B_{4\rho}}r^{-a+1}(A\nabla u\cdot\nabla u)^{p/2}\,dx
\notag\\
&\quad+C\varepsilon\int_{B_{4\rho}}r^{-a+2}(A\nabla u\cdot\nabla u)^{p/2-1}|\nabla u|\,|D^2u|\,dx .
\label{rawge2}
\end{align}
\endgroup
If $1<p<2$ and
\begin{equation}
pa-\frac{a^2}{4}>0,
\end{equation}
then
\begin{align}
&\left(pa-\frac{a^2}{4}\right)
\int_{B_\rho}r^{-a}(A\nabla u\cdot\nabla u)^{p/2-1}u_r^2\,dx
+(n-p-a)\int_{B_\rho}r^{-a}(A\nabla u\cdot\nabla u)^{p/2}\,dx
\notag\\
&\le C\int_{B_{4\rho}\setminus B_\rho}r^{-a}(A\nabla u\cdot\nabla u)^{p/2}\,dx
+C\varepsilon\int_{B_{4\rho}}r^{-a+1}(A\nabla u\cdot\nabla u)^{p/2}\,dx
\notag\\
&\quad+C\varepsilon\int_{B_{4\rho}}r^{-a+2}(A\nabla u\cdot\nabla u)^{p/2-1}|\nabla u|\,|D^2u|\,dx .
\label{rawlt2}
\end{align}
When $\varepsilon=0$, both terms carrying $\varepsilon$ in
\eqref{rawge2} and \eqref{rawlt2} are omitted.
The constant $C$ depends only on $n,p,\lambda,\Lambda$ and on a fixed admissible range of $a$.
\end{theorem}

For a smooth compactly supported $\eta$, choose
$\phi\in C^\infty(\mathbb R)$ with $0\le\phi\le1$, $\phi=0$ on
$(-\infty,1]$ and $\phi=1$ on $[2,\infty)$, and put
\[
 \phi_\tau(x)=\phi\left(\frac{|\nabla u(x)|}{\tau}\right),
 \qquad \xi_\tau=(x\cdot\nabla u)\eta\phi_\tau,\qquad \tau>0.
\]
Then
\[
 \operatorname{supp}\nabla\phi_\tau\subset\{\tau\le|\nabla u|\le2\tau\},
 \qquad
 |\nabla\phi_\tau|\le\frac C\tau|D^2u|
 \mathbf1_{\{\tau<|\nabla u|<2\tau\}}.
\]
For fixed $\tau>0$, the support of $\xi_\tau$ is compactly contained
in $\{|\nabla u|>0\}$. The local $W^{2,2}$ regularity on this set,
established in Appendix~A, gives $\xi_\tau\in W^{1,2}_0(B_1)$.
Since the coefficients of the stability form are bounded on a
neighbourhood of its support, smooth approximation within this
neighbourhood shows that $\xi_\tau\in\mathcal T_u$.

To let $\tau\to0$, we use the weighted integrability from Appendix~A:
\[
 \int_{K\cap\{|\nabla u|>0\}}
 (A\nabla u\cdot\nabla u)^{(p-2)/2}|D^2u|^2\,dx<\infty,
 \qquad K\Subset B_1.
\]
Here $D^2u$ denotes the weak Hessian on the noncritical
set. Expressions on the critical set, including powers at zero,
are interpreted as in the introduction.

\begin{lemma}\label{prod}
For every smooth compactly supported function $\eta$, with $\phi_\tau$ as above,
we first have the cutoff inequality
\begin{align}
&\int \phi_\tau^2(x\cdot\nabla u)^2(A\nabla u\cdot\nabla u)^{p/2-1}A\nabla\eta\cdot\nabla\eta
\notag\\
&\quad+(p-2)\int \phi_\tau^2(x\cdot\nabla u)^2(A\nabla u\cdot\nabla u)^{p/2-2}(A\nabla u\cdot\nabla\eta)^2
\notag\\
&\ge \int \phi_\tau^2\Bigl\{f'(u)(x\cdot\nabla u)^2
+\operatorname{div}\bigl((A\nabla u\cdot\nabla u)^{p/2-1}A\nabla(x\cdot\nabla u)\bigr)(x\cdot\nabla u)
\notag\\
&\qquad +(p-2)\operatorname{div}\bigl((A\nabla u\cdot\nabla u)^{p/2-2}(A\nabla u\cdot\nabla(x\cdot\nabla u))A\nabla u\bigr)(x\cdot\nabla u)\Bigr\}\eta^2
+o_\tau(1).
\label{prodineqcutoff}
\end{align}
Here $o_\tau(1)\to0$ as $\tau\downarrow0$ for the fixed function $\eta$.
\end{lemma}

\begin{proof}
Keep $\tau>0$ fixed and write $c=x\cdot\nabla u$, so
\[
 \nabla\xi_\tau
 =\phi_\tau(c\nabla\eta+\eta\nabla c)+c\eta\nabla\phi_\tau.
\]
Let $K=\operatorname{supp}\eta$ and
$E_\tau=K\cap\{\tau<|\nabla u|<2\tau\}$.
Every term containing $\nabla\phi_\tau$, both in the stability
expansion and when integration by parts differentiates $\phi_\tau^2$,
satisfies
\[
 |E_\tau^{\rm cut}|
 \le C_\eta\int_{E_\tau}
 \left((A\nabla u\cdot\nabla u)^{(p-2)/2}|D^2u|^2
 +(A\nabla u\cdot\nabla u)^{(p-1)/2}|D^2u|\right)dx.
\]
Indeed, on $E_\tau$,
$(A\nabla u\cdot\nabla u)^{1/2}\asymp|\nabla u|\asymp\tau$,
$|c|\le C_K\tau$,
$|\nabla c|\le|\nabla u|+C_K|D^2u|$, and
$|\nabla\phi_\tau|\le C|D^2u|/\tau$.  The first quadratic form
therefore produces precisely these two types of terms; the absolute
value of the second is controlled by the same form by ellipticity,
also when $p-2<0$.  The first density is integrable by the Appendix,
and Cauchy--Schwarz gives
\[
\begin{aligned}
 &\int_K(A\nabla u\cdot\nabla u)^{(p-1)/2}|D^2u|\,dx\\
 &\quad\le\left(\int_{K\cap\{|\nabla u|>0\}}
 (A\nabla u\cdot\nabla u)^{(p-2)/2}|D^2u|^2\,dx\right)^{1/2}
 \left(\int_K(A\nabla u\cdot\nabla u)^{p/2}\,dx\right)^{1/2}.
\end{aligned}
\]
Thus $E_\tau^{\rm cut}=o_\tau(1)$ for fixed $\eta$, by absolute
continuity.  Stability becomes
\begin{align}
&\int\phi_\tau^2(A\nabla u\cdot\nabla u)^{p/2-1}
 A(c\nabla\eta+\eta\nabla c)\cdot(c\nabla\eta+\eta\nabla c)
\notag\\
&\quad+(p-2)\int\phi_\tau^2(A\nabla u\cdot\nabla u)^{p/2-2}
 \bigl(A\nabla u\cdot(c\nabla\eta+\eta\nabla c)\bigr)^2
\notag\\
&\ge\int f'(u)c^2\eta^2\phi_\tau^2+o_\tau(1).
\label{prod-cutoff-stability}
\end{align}
In the two square expansions, integrate the cross terms against
$\nabla(\eta^2)$ by parts.  On the noncritical support all vector
fields below are locally $L^2$ and all pairings are distributional;
the resulting identity is
\begin{align*}
&\int\phi_\tau^2c(A\nabla u\cdot\nabla u)^{p/2-1}
 A\nabla c\cdot\nabla(\eta^2)\\
&\quad+(p-2)\int\phi_\tau^2c(A\nabla u\cdot\nabla u)^{p/2-2}
 (A\nabla u\cdot\nabla c)A\nabla u\cdot\nabla(\eta^2)\\
&=-\int\phi_\tau^2\eta^2
 \left[(A\nabla u\cdot\nabla u)^{p/2-1}A\nabla c\cdot\nabla c
 +(p-2)(A\nabla u\cdot\nabla u)^{p/2-2}
 (A\nabla u\cdot\nabla c)^2\right]\\
&\quad-\int\phi_\tau^2\eta^2c\,
 \operatorname{div}\bigl((A\nabla u\cdot\nabla u)^{p/2-1}A\nabla c\bigr)\\
&\quad-(p-2)\int\phi_\tau^2\eta^2c\,
 \operatorname{div}\bigl((A\nabla u\cdot\nabla u)^{p/2-2}
 (A\nabla u\cdot\nabla c)A\nabla u\bigr)+o_\tau(1).
\end{align*}
The first integral on the right cancels exactly the two
$\eta^2$ square terms in \eqref{prod-cutoff-stability}; derivatives
of $\phi_\tau^2$ are covered by the cutoff bound above.  Moving the
remaining divergence pairings to the right proves
\eqref{prodineqcutoff}.
\end{proof}

\begin{lemma}\label{dirid}
Let $U'\Subset\{|\nabla u|>0\}$. Then
\begin{equation}
\nabla(x\cdot\nabla u)=\nabla u+(D^2u)x
\end{equation}
holds almost everywhere in $U'$, and the following identity holds
in the sense of distributions in $U'$:
\begin{align}
&f'(u)(x\cdot\nabla u)+\operatorname{div}\bigl((A\nabla u\cdot\nabla u)^{p/2-1}A\nabla(x\cdot\nabla u)\bigr)
\notag\\
&\quad +(p-2)\operatorname{div}\bigl((A\nabla u\cdot\nabla u)^{p/2-2}(A\nabla u\cdot\nabla(x\cdot\nabla u))A\nabla u\bigr)
\notag\\
&=\partial_i\Bigl[(A\nabla u\cdot\nabla u)^{p/2-1}(pA_{ij}-x_kA_{ij,k})u_j
\notag\\
&\qquad\qquad-\left(\frac p2-1\right)(A\nabla u\cdot\nabla u)^{p/2-2}(u_mx_kA_{m\ell,k}u_\ell)A_{ij}u_j\Bigr].
\label{diridentity}
\end{align}
\end{lemma}

\begin{proof}
By Appendix~A, $u\in W^{2,2}(U')$, so the first identity follows
from the Sobolev product rule. Set
\[
 W_i=(A\nabla u\cdot\nabla u)^{p/2-1}A_{ij}u_j.
\]
Since $|\nabla u|$ is bounded away from zero on $\overline{U'}$,
the Sobolev chain rule gives $W_i,f(u)\in W^{1,2}(U')$.
Differentiating $-\partial_iW_i=f(u)$ in distributions yields
\[
 f'(u)(x\cdot\nabla u)
 =\partial_i\bigl(W_i-x_k\partial_kW_i\bigr).
\]
Using $\partial_j(x_ku_k)=u_j+x_ku_{kj}$, the symmetries of $A$
and $D^2u$, and $2(p/2-1)=p-2$, the chain rule gives
\begin{align*}
 x_k\partial_kW_i
 &=(A\nabla u\cdot\nabla u)^{p/2-1}
   A_{ij}\partial_j(x_ku_k)\\
 &\quad +(p-2)(A\nabla u\cdot\nabla u)^{p/2-2}
   (A\nabla u\cdot\nabla(x\cdot\nabla u))A_{ij}u_j\\
 &\quad -(p-1)W_i
   +(A\nabla u\cdot\nabla u)^{p/2-1}x_kA_{ij,k}u_j\\
 &\quad +\left(\frac p2-1\right)
   (A\nabla u\cdot\nabla u)^{p/2-2}
   (u_mx_kA_{m\ell,k}u_\ell)A_{ij}u_j
\end{align*}
almost everywhere in $U'$. Substituting this expression into the
preceding distributional identity gives \eqref{diridentity}.
All vector fields involved are in $L^2(U')$; the final divergence
is taken in distributions, so no second derivative of $A$ is needed.
\end{proof}

\begin{lemma}\label{explicit}
Let $\eta\in C_c^\infty(B_1)$. For $\tau>0$, one has
\begin{align}
&\int \phi_\tau^2(x\cdot\nabla u)^2
(A\nabla u\cdot\nabla u)^{p/2-1}
A\nabla\eta\cdot\nabla\eta\,dx
\notag\\
&\quad +(p-2)\int \phi_\tau^2(x\cdot\nabla u)^2
(A\nabla u\cdot\nabla u)^{p/2-2}
(A\nabla u\cdot\nabla\eta)^2\,dx
\notag\\
&\ge \int \partial_i\Bigl[
(A\nabla u\cdot\nabla u)^{p/2-1}
(pA_{ij}-x_kA_{ij,k})u_j
\notag\\
&\qquad\qquad
-\left(\frac p2-1\right)
(A\nabla u\cdot\nabla u)^{p/2-2}
(u_mx_kA_{m\ell,k}u_\ell)A_{ij}u_j
\Bigr](x\cdot\nabla u)\eta^2\phi_\tau^2\,dx
+o_\tau(1),
\label{explicitineqcutoff}
\end{align}
where $o_\tau(1)\to0$ as $\tau\downarrow0$ for fixed $\eta$.
Passing to the limit gives
\begin{align}
&\int (x\cdot\nabla u)^2
(A\nabla u\cdot\nabla u)^{p/2-1}
A\nabla\eta\cdot\nabla\eta\,dx
\notag\\
&\quad +(p-2)\int (x\cdot\nabla u)^2
(A\nabla u\cdot\nabla u)^{p/2-2}
(A\nabla u\cdot\nabla\eta)^2\,dx
\notag\\
&\ge \int \partial_i\Bigl[
(A\nabla u\cdot\nabla u)^{p/2-1}
(pA_{ij}-x_kA_{ij,k})u_j
\notag\\
&\qquad\qquad
-\left(\frac p2-1\right)
(A\nabla u\cdot\nabla u)^{p/2-2}
(u_mx_kA_{m\ell,k}u_\ell)A_{ij}u_j
\Bigr](x\cdot\nabla u)\eta^2\,dx.
\label{explicitineq}
\end{align}
\end{lemma}

\begin{proof}
Fix $\tau>0$ and put $K=\operatorname{supp}\eta$.
Since $\phi_\tau=0$ where $|\nabla u|\le\tau$, choose
$U_\tau\Subset\{|\nabla u|>0\}$ containing
$K\cap\operatorname{supp}\phi_\tau$.
By Appendix~A, $u\in W^{2,2}(U_\tau)$.
The coefficients and gradient are bounded there, and
$|\nabla u|$ is bounded away from zero. Hence the vector fields
in \eqref{diridentity} and $f'(u)(x\cdot\nabla u)$ belong to
$L^2(U_\tau)$, while
\[
 (x\cdot\nabla u)\eta^2\phi_\tau^2
 \in W_0^{1,2}(U_\tau).
\]
Smooth approximation therefore allows us to pair
\eqref{diridentity} with this function. Substitution into
\eqref{prodineqcutoff} gives \eqref{explicitineqcutoff}.

After distributional integration by parts on the right-hand side,
the terms containing $\nabla\phi_\tau$ are bounded in absolute
value by
\[
 C_\eta
 \int_{K\cap\{\tau<|\nabla u|<2\tau\}}
 (A\nabla u\cdot\nabla u)^{(p-1)/2}|D^2u|\,dx
 \longrightarrow0.
\]
This bound follows from ellipticity, the local bound on $DA$,
and the cutoff estimates in Lemma~\ref{prod}, with $C_\eta$
independent of $\tau$.
The convergence follows from the weighted Hessian integrability
and Cauchy--Schwarz, as in Lemma~\ref{prod}.

The remaining integrands are bounded by
\[
 C_\eta\left(
 (A\nabla u\cdot\nabla u)^{p/2}
 +(A\nabla u\cdot\nabla u)^{(p-1)/2}|D^2u|
 \right),
\]
which is integrable on $K\cap\{|\nabla u|>0\}$.
Dominated convergence now yields \eqref{explicitineq},
with its divergence pairing defined by this limit.
\end{proof}

\begin{remark}\label{radial-weight-regularization}
For $0<\delta<\rho$, put
\[
 w_\delta(r)=\min\{r^{-a/2},\delta^{-a/2}\},\qquad
 \eta_\delta=w_\delta\zeta.
\]
At fixed $\delta$, use compactly supported smooth mollifications of
$\eta_\delta$ with uniform convergence, almost-everywhere convergence
of the gradients, and uniformly bounded Lipschitz constants.  Their
supports lie in a common compact subset of $B_{4\rho}$.  Apply the
critical-set cutoff argument of Lemma~\ref{explicit} at this fixed
radial truncation and first let $\tau\downarrow0$.  Dominated
convergence then justifies $\eta_\delta$ in its distributional pairings.
Only afterwards do we let $\delta\downarrow0$ to use
$\eta=r^{-a/2}\zeta$.

If $\varepsilon=0$, \eqref{normass} gives $A\equiv I$.  Both $DA$
fields in \eqref{explicitineq} and all errors involving $A-I$ vanish.
At fixed $\delta$, perform the cutoff integrations and the Pohozaev
cancellations below before letting $\delta\downarrow0$.  The remaining
terms are controlled by $r^{-a}(A\nabla u\cdot\nabla u)^{p/2}$ on
$\operatorname{supp}\zeta$, with the small-ball bound below; thus no
Hessian-error majorant is needed in this case.
For $\varepsilon>0$, the raw inequality is immediate if
\[
 \int_{B_{4\rho}}r^{-a+2}(A\nabla u\cdot\nabla u)^{p/2-1}
          |\nabla u|\,|D^2u|\,dx=+\infty,
\]
since its right-hand side is then infinite; otherwise this integral
dominates the truncated Hessian errors.  The remaining energy and radial terms, on the compact support
of $\zeta$, are bounded by a constant times
$r^{-a}(A\nabla u\cdot\nabla u)^{p/2}$, which is integrable because
$u\in C^1$ and $0<a<n$.  Inside $B_\delta$, the radial derivative of
$w_\delta$ vanishes, and the missing limiting terms are controlled by
\[
 \int_{B_\delta}r^{-a}(A\nabla u\cdot\nabla u)^{p/2}\,dx
 \le C\|\nabla u\|_{L^\infty(B_\rho)}^p\delta^{n-a}
 \longrightarrow0.
\]
Together with dominated convergence outside $B_\delta$, this proves
the radial limiting formulas below.  For $1<p<2$, the Hessian density
is understood on $\{|\nabla u|>0\}$ and extended by zero.
\end{remark}

\begin{proof}[Proof of Theorem~\ref{raw}]
\normalfont
We take $\eta=r^{-a/2}\zeta$ in \eqref{explicitineq}, integrating by parts
and cancelling terms before removing the radial truncation, as justified
in Remark~\ref{radial-weight-regularization}. Then
\begin{align*}
 \eta^2&=r^{-a}\zeta^2,\qquad
 \nabla\eta=r^{-a/2}\nabla\zeta-\frac a2r^{-a/2-2}x\zeta,\\
 \nabla(\eta^2)&=-ar^{-a-2}x\zeta^2+2r^{-a}\zeta\nabla\zeta.
\end{align*}
The bound $|A-I|\le\varepsilon r$ gives
\begin{equation}
 |Ax\cdot x-r^2|\le C\varepsilon r^3,\qquad
 |A\nabla u\cdot x-x\cdot\nabla u|\le C\varepsilon r^2|\nabla u|.
 \label{comp}
\end{equation}

On $B_\rho$, where $\zeta=1$, \eqref{comp} gives the first quadratic
term in \eqref{explicitineq} the radial coefficient $a^2/4$, with error
bounded by
$C\varepsilon\int_{B_\rho}r^{1-a}(A\nabla u\cdot\nabla u)^{p/2}$.
For $p\ge2$,
\[
 (A\nabla u\cdot x)^2
 \le(A\nabla u\cdot\nabla u)(Ax\cdot x)
\]
bounds the second quadratic term by $p-2$ times the first.
On $B_{4\rho}\setminus B_\rho$, $r\asymp\rho$ and
$|\nabla\zeta|\le C/\rho$ control both terms by the annular energy.
Thus the left side of \eqref{explicitineq} is at most
\begin{align}
 &\frac{(p-1)a^2}{4}\int_{B_\rho}r^{-a}
       (A\nabla u\cdot\nabla u)^{p/2-1}u_r^2
 +C\int_{B_{4\rho}\setminus B_\rho}r^{-a}(A\nabla u\cdot\nabla u)^{p/2}
 \notag\\
 &\quad+C\varepsilon\int_{B_{4\rho}}r^{1-a}(A\nabla u\cdot\nabla u)^{p/2}.
 \label{upperge2}
\end{align}
For $1<p<2$, dropping the nonpositive second term gives the upper bound
\begin{align}
 &\frac{a^2}{4}\int_{B_\rho}r^{-a}
       (A\nabla u\cdot\nabla u)^{p/2-1}u_r^2
 +C\int_{B_{4\rho}\setminus B_\rho}r^{-a}(A\nabla u\cdot\nabla u)^{p/2}
 \notag\\
 &\quad+C\varepsilon\int_{B_{4\rho}}r^{1-a}(A\nabla u\cdot\nabla u)^{p/2}.
 \label{upperlt2}
\end{align}

For the lower bound, pair the divergence in \eqref{explicitineq}
with $(x\cdot\nabla u)\eta^2$. At fixed $\delta$, retain $\phi_\tau$
in the integrations by parts. The extra cutoff terms are bounded by
\[
 C_\delta\int_{\operatorname{supp}\eta_\delta\cap\{\tau<|\nabla u|<2\tau\}}
 (A\nabla u\cdot\nabla u)^{(p-1)/2}|D^2u|\,dx\longrightarrow0
\]
by Cauchy--Schwarz and the preliminary weighted Hessian bound,
as in Lemma~\ref{prod}; the coefficient-derivative energy terms
are locally integrable.

On the noncritical set, the chain rule and Lemma~\ref{dirid} give
\begin{align*}
 &p(A\nabla u\cdot\nabla u)^{p/2-1}
 A\nabla u\cdot\nabla(x\cdot\nabla u)\\
 &=p(A\nabla u\cdot\nabla u)^{p/2}
   +x\cdot\nabla\bigl((A\nabla u\cdot\nabla u)^{p/2}\bigr)\\
 &\quad-\frac p2(A\nabla u\cdot\nabla u)^{p/2-1}
       u_i x_kA_{ij,k}u_j.
\end{align*}
Using $\operatorname{div}(x\eta^2)=n\eta^2+x\cdot\nabla(\eta^2)$,
the contribution of the principal field is therefore
\begin{equation}
\begin{aligned}
 &-p\int_{B_{4\rho}}(A\nabla u\cdot\nabla u)^{p/2-1}
       A\nabla u\cdot\nabla\bigl((x\cdot\nabla u)\eta^2\bigr)\\
 &=(n-p)\int_{B_{4\rho}}(A\nabla u\cdot\nabla u)^{p/2}\eta^2
   +\int_{B_{4\rho}}(A\nabla u\cdot\nabla u)^{p/2}x\cdot\nabla(\eta^2)\\
 &\quad+\frac p2\int_{B_{4\rho}}(A\nabla u\cdot\nabla u)^{p/2-1}
       (u_i x_kA_{ij,k}u_j)\eta^2\\
 &\quad-p\int_{B_{4\rho}}(x\cdot\nabla u)
       (A\nabla u\cdot\nabla u)^{p/2-1}A\nabla u\cdot\nabla(\eta^2).
\end{aligned}
 \label{poho}
\end{equation}
On $B_\rho$, the first two terms give the energy coefficient
$n-p-a$. The last term, by \eqref{comp}, satisfies
\begin{align}
 &pa\int_{B_\rho}r^{-a-2}(x\cdot\nabla u)(A\nabla u\cdot x)
       (A\nabla u\cdot\nabla u)^{p/2-1}\notag\\
 &\quad\ge pa\int_{B_\rho}r^{-a}(A\nabla u\cdot\nabla u)^{p/2-1}u_r^2
       -C\varepsilon\int_{B_\rho}r^{1-a}(A\nabla u\cdot\nabla u)^{p/2}.
 \label{mainradial}
\end{align}

Both $DA$ fields in \eqref{explicitineq} have magnitude at most
\[
 C\varepsilon r(A\nabla u\cdot\nabla u)^{p/2-1}|\nabla u|.
\]
For the second, use $|u_mx_kA_{m\ell,k}u_\ell|\le C\varepsilon r
(A\nabla u\cdot\nabla u)$. After integration by parts, estimate
both pairings with $\nabla((x\cdot\nabla u)\eta^2)$ in absolute
value, also for $p<2$. In the decomposition
\begin{equation}
 \partial_i((x\cdot\nabla u)\eta^2)
 =u_i\eta^2+x_ku_{ki}\eta^2+(x\cdot\nabla u)\partial_i(\eta^2),
 \label{derivativeceta}
\end{equation}
the first and third pieces, together with the $DA$ term in
\eqref{poho}, are bounded by
\begin{equation}
\begin{aligned}
 &C\varepsilon\int_{B_{4\rho}}
  \bigl(r^{1-a}\zeta^2+r^{2-a}|\zeta|\,|\nabla\zeta|\bigr)
  (A\nabla u\cdot\nabla u)^{p/2}\\
 &\qquad\le C\varepsilon\int_{B_{4\rho}}r^{1-a}
  (A\nabla u\cdot\nabla u)^{p/2}.
\end{aligned}
 \label{BnoHess}
\end{equation}
The only Hessian contribution, from $x_ku_{ki}\eta^2$, has absolute value at most
\begin{equation}
 C\varepsilon\int_{B_{4\rho}}r^{2-a}\zeta^2
 (A\nabla u\cdot\nabla u)^{p/2-1}|\nabla u|\,|D^2u|.
 \label{BHess}
\end{equation}
Before enlargement, coefficient errors are supported in
$\operatorname{supp}\zeta$ and $\nabla\zeta$ terms in
$B_{4\rho}\setminus B_\rho$. If $\operatorname{supp}\zeta\subset B_{2\rho}$,
these supports lie in $B_{2\rho}$ and $B_{2\rho}\setminus B_\rho$,
respectively.

By \eqref{poho}--\eqref{BHess}, $0\le\zeta\le1$, and the annular
energy bound, the right side of \eqref{explicitineq} is at least
\begin{align}
 &pa\int_{B_\rho}r^{-a}(A\nabla u\cdot\nabla u)^{p/2-1}u_r^2
 +(n-p-a)\int_{B_\rho}r^{-a}(A\nabla u\cdot\nabla u)^{p/2}\notag\\
 &\quad-C\int_{B_{4\rho}\setminus B_\rho}r^{-a}(A\nabla u\cdot\nabla u)^{p/2}
       -C\varepsilon\int_{B_{4\rho}}r^{1-a}(A\nabla u\cdot\nabla u)^{p/2}\notag\\
 &\quad-C\varepsilon\int_{B_{4\rho}}r^{2-a}
       (A\nabla u\cdot\nabla u)^{p/2-1}|\nabla u|\,|D^2u|.
 \label{rhsbelow}
\end{align}
Comparing with \eqref{upperge2} or \eqref{upperlt2} proves
\eqref{rawge2} and \eqref{rawlt2}.
\end{proof}

\section{Anisotropic Hessian estimates and higher integrability}\label{sec:hessian-estimate}

This section proves the Hessian estimates in Theorem~1.2 and the
higher integrability of the gradient. We first combine stability with
the equation to estimate the Hessian, and then use the coarea formula
and the Sobolev--Poincar\'e inequality to obtain higher integrability.
We use the definitions of regular solutions and stability from
Section~1.1; admissible stability tests belong to $\mathcal T_u$ in
(1.4). Repeated indices are summed.

For the solutions considered below, $D^2u$ denotes the weak Hessian on
the noncritical set, supplied by Appendix~A. On this set we write
\begin{equation}\label{c-e-N-M-def}
 e=\frac{\nabla u}{|\nabla u|},\qquad
 N=D^2u(e,e),\qquad M=D^2u-Ne\otimes e.
\end{equation}
The quantities $|D^2u|$ and $|M|$ denote the Frobenius norms.

\begin{thm}\label{thm:hessian-estimate}
Let $p>1$ and let $u$ be a stable regular solution
of (1.1), with $f\in C^1(\mathbb R)$ nonnegative. Let
$A\in W^{1,\infty}_{\rm loc}(\Ball_1;\mathbb R^{n\times n})$ be real symmetric.
We use its locally Lipschitz representative and assume
\begin{equation}\label{A-elliptic}
 \lambda|\xi|^2\le A(x)\xi\cdot\xi\le\Lambda|\xi|^2
 \qquad\text{for }x\in\Ball_1,\ \xi\in\mathbb R^n.
\end{equation}
There exists $\delta_0=\delta_0(n,p,\lambda,\Lambda)\in(0,1]$
such that, if
\begin{equation}\label{A-small-gradient}
 \max_{i,j,k}\|A_{ij,k}\|_{L^\infty(\Ball_1)}\le\delta_0,
\end{equation}
the following estimates hold with $C=C(n,p,\lambda,\Lambda)$.

If \(p\ge 2\), then
\begin{equation}\label{hessian-p-ge-2}
 \int_{\Ball_{1/2}\cap\{|\nabla u|>0\}}(A\nabla u\cdot\nabla u)^{(p-2)/2}|D^2u|\dd
 \le C\left(\int_{\Ball_1}(A\nabla u\cdot\nabla u)^{p/2}\dd\right)^{(p-1)/p}.
\end{equation}
If \(1<p<2\), then
\begin{equation}\label{hessian-p-lt-2}
 \int_{\Ball_{1/2}\cap\{|\nabla u|>0\}}|D^2u|\dd
 \le C\left(\int_{\Ball_1}(A\nabla u\cdot\nabla u)^{p/2}\dd\right)^{1/p}.
\end{equation}
Moreover, for every \(\eta\in C_c^\infty(\Ball_1)\),
\begin{equation}\label{M-direct-estimate}
 \int_{\Ball_1\cap\{|\nabla u|>0\}}(A\nabla u\cdot\nabla u)^{(p-2)/2}|M|^2\eta^2\dd
 \le C\int_{\Ball_1}(A\nabla u\cdot\nabla u)^{p/2}
        (|\nabla\eta|^2+\eta^2)\dd .
\end{equation}
\end{thm}

No normalization $A(0)=I$ is required. We write $A_{ij,k}$ for the weak
derivative $\partial_k A_{ij}$; pointwise identities involving these
derivatives hold almost everywhere. The estimates use only the
ellipticity bounds and the first weak derivatives of $A$.

\begin{proof}
By Theorem~\ref{dq:main-theorem}, for every $K\Subset\Ball_1$,
\[
 \int_{K\cap\{|\nabla u|>0\}}
 (A\nabla u\cdot\nabla u)^{(p-2)/2}|D^2u|^2\dd<\infty,
\]
and $u\in W^{2,2}_{\rm loc}(\{|\nabla u|>0\})$.
For all $p>1$, this preliminary input follows from the weak equation
and the local $C^1$ regularity supplied by
Section~\ref{intro:solution-class}.

\medskip
\noindent\textbf{Step 1. }
  Define
\begin{equation}\label{Q-definition}
 Q_{ij}
 =(p-2)(A_{ab}u_au_b)^{p/2-2}(A_{i\ell}u_\ell)(A_{jm}u_m)
 +(A_{ab}u_au_b)^{p/2-1}A_{ij}.
\end{equation}
For admissible tests supported away from the critical set, stability reads
\begin{equation}\label{stability-Q-form}
 \int_{\Ball_1}f'(u)\xi^2\dd
 \le\int_{\Ball_1\cap\{|\nabla u|>0\}}Q_{ij}\xi_i\xi_j\dd.
\end{equation}
Cauchy's inequality in the $A$ scalar product gives
$(A\nabla u\cdot\zeta)^2\le
(A\nabla u\cdot\nabla u)(A\zeta\cdot\zeta)$.
Accounting for the sign of $p-2$ and using ellipticity, we obtain
\begin{align}
 Q_{ij}\zeta_i\zeta_j
 &\ge\mu(A\nabla u\cdot\nabla u)^{(p-2)/2}|\zeta|^2,
 \label{Q-lower}\\*
 |Q_{ij}\zeta_i\zeta_j|
 &\le C(A\nabla u\cdot\nabla u)^{(p-2)/2}|\zeta|^2,
 \label{Q-upper}
\end{align}
where $\mu=\lambda\min\{1,p-1\}>0$ and
$C=\Lambda\max\{1,p-1\}$.

Fix $U\Subset U'\Subset\{|\nabla u|>0\}$ and set
\[
 F_i(x,\xi)=(A(x)\xi\cdot\xi)^{p/2-1}A_{ij}(x)\xi_j.
\]
On $U'$, $u\in W^{2,2}$, $A\in W^{1,\infty}$, and $\nabla u$ is
bounded and bounded away from zero. The Sobolev chain rule therefore
gives $F_i(x,\nabla u),f(u)\in W^{1,2}(U')$ and
\begin{equation}\label{flux-differentiation}
 \partial_k\left((A_{ab}u_au_b)^{p/2-1}A_{ij}u_j\right)
 =Q_{ij}u_{jk}+R_i^{(k)},
\end{equation}
where
\begin{align}\label{R-definition}
 R_i^{(k)}
 &=\left(\frac p2-1\right)(A_{ab}u_au_b)^{p/2-2}
       A_{mn,k}u_mu_n(A_{i\ell}u_\ell)\notag\\*
 &\quad+(A_{ab}u_au_b)^{p/2-1}A_{ij,k}u_j.
\end{align}
Differentiating (1.1) in distributions and extending by density yields
\begin{equation}\label{differentiated-equation}
 \int_U\bigl(Q_{ij}u_{jk}+R_i^{(k)}\bigr)\partial_i\varphi\dd
 =\int_U f'(u)u_k\varphi\dd
 \qquad\text{for }\varphi\in W_0^{1,2}(U).
\end{equation}
Both sides are continuous in this test space since the vector field
and $f'(u)u_k$ belong to $L^2(U)$. Ellipticity and \eqref{A-small-gradient} bound
each term in \eqref{R-definition}, giving
\begin{equation}\label{R-pointwise-bound}
 |R_i^{(k)}|\le C\delta_0(A\nabla u\cdot\nabla u)^{(p-1)/2},
\end{equation}
with $C=C(n,p,\lambda,\Lambda)$.

\medskip
\noindent\textbf{Step 2. }
Let $\eta\in C_c^\infty(\Ball_1)$ and choose a nondecreasing
$\chi\in C^\infty([0,\infty))$ with $0\le\chi\le1$,
$\chi=0$ on $[0,1]$, $\chi=1$ on $[2,\infty)$, and $|\chi'|\le C$. Set
\[
 \chi_\tau=\chi(|\nabla u|/\tau),\qquad \psi_\tau=|\nabla u|\eta\chi_\tau,
 \qquad E_\tau=\{\tau<|\nabla u|<2\tau\}.
\]
For fixed $\tau>0$, noncritical $W^{2,2}$ regularity gives
$\psi_\tau,u_k\eta^2\chi_\tau^2\in W_0^{1,2}$ with compact support
in $\{|\nabla u|>0\}$. Since the weight is bounded on a neighbourhood of this
support, smooth approximation makes them admissible in
\eqref{stability-Q-form} and \eqref{differentiated-equation}, respectively.
There,
\[
 \partial_i|\nabla u|=\frac{u_ku_{ki}}{|\nabla u|},\qquad
 |\nabla\chi_\tau|\le C\tau^{-1}|D^2u|\mathbf1_{E_\tau}.
\]
Test the linearized equations by $u_k\eta^2\chi_\tau^2$, sum over $k$,
and compare with stability tested by $\psi_\tau$. The identities
$\sum_k u_k^2=|\nabla u|^2$ and $\sum_k u_k\nabla u_k=|\nabla u|\nabla|\nabla u|$ cancel the
reaction and mixed quadratic terms, giving
\begin{align*}
 &\int_{\Ball_1\cap\{|\nabla u|>0\}}
 (Q_{ij}u_{ki}u_{kj}-Q_{ij}(\partial_i|\nabla u|)(\partial_j|\nabla u|))\eta^2\chi_\tau^2\dd\\*
 &\quad+\int_{\Ball_1\cap\{|\nabla u|>0\}}
 R_i^{(k)}\partial_i(u_k\eta^2\chi_\tau^2)\dd
 \le\int_{\Ball_1\cap\{|\nabla u|>0\}}
 |\nabla u|^2Q\nabla(\eta\chi_\tau)\cdot\nabla(\eta\chi_\tau)\dd.
\end{align*}
After expansion, the terms containing $\nabla\chi_\tau$ are supported
in $E_\tau$. There, \eqref{Q-upper}, \eqref{R-pointwise-bound}, and
$|\nabla u|/\tau\le2$ give
\begin{align*}
 |\nabla u|^2|\eta Q\nabla\eta\cdot\nabla\chi_\tau|
 +|R_i^{(k)}u_k\eta^2\partial_i\chi_\tau|
 &\le C|\nabla u|^{p-1}|D^2u||\eta|(|\nabla\eta|+\delta_0|\eta|),\\*
 |\nabla u|^2\eta^2|Q\nabla\chi_\tau\cdot\nabla\chi_\tau|
 &\le C|\nabla u|^{p-2}|D^2u|^2\eta^2.
\end{align*}
Their integrals tend to zero by Cauchy--Schwarz with weights
$|\nabla u|^{p-2}$ and $|\nabla u|^p$, and absolute continuity of the preliminary
weighted Hessian integral. The same bounds give integrable majorants
for the remaining terms. Letting $\tau\downarrow0$, we obtain
\begin{align}\label{main-before-lower-bound}
 &\int_{\Ball_1\cap\{|\nabla u|>0\}}
 (Q_{ij}u_{ki}u_{kj}-Q_{ij}(\partial_i|\nabla u|)(\partial_j|\nabla u|))\eta^2\dd
 +\left\langle-\partial_iR_i^{(k)},u_k\eta^2\right\rangle\notag\\*
 &\hspace{30mm}\le
 \int_{\Ball_1\cap\{|\nabla u|>0\}}|\nabla u|^2Q\nabla\eta\cdot\nabla\eta\dd.
\end{align}
Here the pairing is defined by this cutoff limit, and
\eqref{R-pointwise-bound} gives
\begin{align}\label{R-error-integrate}
 \left|\left\langle-\partial_iR_i^{(k)},u_k\eta^2\right\rangle\right|
 &=\left|\int_{\Ball_1\cap\{|\nabla u|>0\}}
 R_i^{(k)}(u_{ki}\eta^2+2u_k\eta\eta_i)\dd\right|\notag\\*
 &\le C\delta_0\int_{\Ball_1\cap\{|\nabla u|>0\}}
 (A\nabla u\cdot\nabla u)^{(p-1)/2}|D^2u|\eta^2\dd\notag\\*
 &\quad+C\delta_0\int_{\Ball_1}
 (A\nabla u\cdot\nabla u)^{p/2}|\eta||\nabla\eta|\dd.
\end{align}

At a noncritical point, choose an orthonormal basis with $e=e_n$.
This pointwise change of basis gives $\partial_i|\nabla u|=u_{ni}$, and Hessian symmetry yields
\begin{align*}
 Q_{ij}u_{ki}u_{kj}-Q_{ij}(\partial_i|\nabla u|)(\partial_j|\nabla u|)
 &=\sum_{\alpha=1}^{n-1}Q(D^2u\,e_\alpha)\cdot(D^2u\,e_\alpha),\\*
 |M|^2
 &=\sum_{\alpha,\beta<n}u_{\alpha\beta}^2
   +2\sum_{\alpha<n}u_{\alpha n}^2
 \le2\sum_{\alpha=1}^{n-1}|D^2u\,e_\alpha|^2.
\end{align*}
Thus \eqref{Q-lower} gives
\begin{equation}\label{M-positive-lower}
 Q_{ij}u_{ki}u_{kj}-Q_{ij}(\partial_i|\nabla u|)(\partial_j|\nabla u|)
 \ge\frac\mu2(A\nabla u\cdot\nabla u)^{(p-2)/2}|M|^2,
\end{equation}
while \eqref{Q-upper} bounds the right-hand density of
\eqref{main-before-lower-bound} by
$C(A\nabla u\cdot\nabla u)^{p/2}|\nabla\eta|^2$.

\medskip
\noindent\textbf{Step 3. }
On the noncritical set, (1.1) gives
\begin{equation}\label{equation-expanded}
 Q_{ij}u_{ij}+E_0=-f(u),
\end{equation}
where
\begin{equation}\label{E0-def}
 E_0=\left(\frac p2-1\right)(A_{ab}u_au_b)^{p/2-2}
 A_{mn,i}u_mu_n(A_{ij}u_j)
 +(A_{ab}u_au_b)^{p/2-1}A_{ij,i}u_j.
\end{equation}
As in \eqref{R-pointwise-bound},
$|E_0|\le C\delta_0(A\nabla u\cdot\nabla u)^{(p-1)/2}$.
Since $D^2u=M+Ne\otimes e$,
\begin{equation}\label{normal-equation}
 (Qe\cdot e)N=-f(u)-E_0-Q:M.
\end{equation}
Thus $f\ge0$ and \eqref{Q-lower}--\eqref{Q-upper} imply
\begin{equation}\label{N-positive-bound}
 N_+\le C|M|+C\delta_0(A\nabla u\cdot\nabla u)^{1/2},
 \qquad N_+=\max\{N,0\}.
\end{equation}

Fix $\beta\ge\max\{0,(p-2)/2\}$. On $\{|\nabla u|>0\}$,
$N=\nabla|\nabla u|\cdot e$ and $|\nabla u|\operatorname{div}e=\operatorname{tr}M$, so
\[
 \operatorname{div}(|\nabla u|^{\beta+1}e)
 =|\nabla u|^\beta\bigl((\beta+1)N+\operatorname{tr}M\bigr).
\]
Test this identity with $\eta^2\chi_\tau^2$ and integrate by parts.
The term differentiating $\chi_\tau^2$ satisfies
\begin{align*}
 &\left|\frac2{\beta+1}\int_{E_\tau}
 |\nabla u|^{\beta+1}e_i\eta^2\chi_\tau\partial_i\chi_\tau\dd\right|\\*
 &\quad\le C
 \left(\int_{E_\tau}|\nabla u|^{p-2}|D^2u|^2\eta^2\dd\right)^{1/2}
 \left(\int_{E_\tau}|\nabla u|^{2\beta-(p-2)}\eta^2\dd\right)^{1/2}
 \longrightarrow0.
\end{align*}
Here $2\beta-(p-2)\ge0$, so local boundedness of $|\nabla u|$ and the preliminary
weighted Hessian integrability justify the limit. The same estimate on
$\operatorname{supp}\eta$ gives $|\nabla u|^\beta|D^2u|\eta^2\in L^1$, while
$|\nabla u|^{\beta+1}|\eta||\nabla\eta|$ is locally integrable. Dominated convergence yields
\begin{align}\label{signed-N-identity}
 \int_{\Ball_1\cap\{|\nabla u|>0\}}|\nabla u|^\beta N\eta^2\dd
 &=-\frac1{\beta+1}\int_{\Ball_1\cap\{|\nabla u|>0\}}
 |\nabla u|^\beta\operatorname{tr}M\,\eta^2\dd\notag\\*
 &\quad-\frac2{\beta+1}\int_{\Ball_1\cap\{|\nabla u|>0\}}
 |\nabla u|^{\beta+1}e_i\eta\eta_i\dd.
\end{align}
All integrals are finite. Using $|N|=2N_+-N$ and
\eqref{N-positive-bound} first bounds the absolute-value integral with weight $|\nabla u|^\beta$;
ellipticity then compares the nonnegative weighted integrals, giving
\begin{align}\label{N-absolute-aniso-beta}
 &\int_{\Ball_1\cap\{|\nabla u|>0\}}
 (A\nabla u\cdot\nabla u)^{\beta/2}|N|\eta^2\dd\notag\\*
 &\quad\le C\int_{\Ball_1\cap\{|\nabla u|>0\}}
 (A\nabla u\cdot\nabla u)^{\beta/2}|M|\eta^2\dd\notag\\*
 &\qquad+C\delta_0\int_{\Ball_1}
 (A\nabla u\cdot\nabla u)^{(\beta+1)/2}\eta^2\dd\notag\\*
 &\qquad+C\int_{\Ball_1}
 (A\nabla u\cdot\nabla u)^{(\beta+1)/2}|\eta||\nabla\eta|\dd.
\end{align}
Here $C$ may also depend on $\beta$, whose later choices depend only on $p$;
$f\ge0$ is used only in \eqref{N-positive-bound}.

Since $p-1\ge\max\{0,(p-2)/2\}$, use \eqref{N-absolute-aniso-beta}
with $\beta=p-1$ and $|D^2u|\le|M|+|N|$ in
\eqref{R-error-integrate} to obtain
\begin{align}\label{R-error-before-young}
 \left|\left\langle-\partial_iR_i^{(k)},u_k\eta^2\right\rangle\right|
 &\le C\delta_0\int_{\Ball_1\cap\{|\nabla u|>0\}}
 (A\nabla u\cdot\nabla u)^{(p-1)/2}|M|\eta^2\dd\notag\\*
 &\quad+C\int_{\Ball_1}(A\nabla u\cdot\nabla u)^{p/2}
 (\delta_0^2\eta^2+\delta_0|\eta||\nabla\eta|)\dd.
\end{align}
By Young's inequality, for every $\varepsilon>0$ we have
\begin{align*}
 \left|\left\langle-\partial_iR_i^{(k)},u_k\eta^2\right\rangle\right|
 &\le\varepsilon\int_{\Ball_1\cap\{|\nabla u|>0\}}
 (A\nabla u\cdot\nabla u)^{(p-2)/2}|M|^2\eta^2\dd\\*
 &\quad+\int_{\Ball_1}(A\nabla u\cdot\nabla u)^{p/2}
 (C_\varepsilon\delta_0^2\eta^2+C|\nabla\eta|^2)\dd.
\end{align*}
Insert this in \eqref{main-before-lower-bound}. By \eqref{M-positive-lower},
$\varepsilon=\mu/4$ absorbs the first term into the $\mu/2$ lower bound.
Since $\delta_0\le1$, this proves \eqref{M-direct-estimate} with
$C=C(n,p,\lambda,\Lambda)$.

\medskip
\noindent\textbf{Step 4. }
Using a cutoff equal to one on \(\Ball_{3/4}\), supported in
\(\Ball_1\), with bounded gradient, \eqref{M-direct-estimate} gives
\begin{equation}\label{M-L2-global}
 \int_{\Ball_{3/4}\cap\{|\nabla u|>0\}}
 (A\nabla u\cdot\nabla u)^{(p-2)/2}|M|^2\dd
 \le C\int_{\Ball_1}(A\nabla u\cdot\nabla u)^{p/2}\dd .
\end{equation}
For \(p\ge2\), Cauchy's and H\"older's inequalities give
\begin{align}\label{M-L1-pge2}
 &\int_{\Ball_{3/4}\cap\{|\nabla u|>0\}}
 (A\nabla u\cdot\nabla u)^{(p-2)/2}|M|\dd\notag\\*
 &\quad\le C\left(\int_{\Ball_1}
 (A\nabla u\cdot\nabla u)^{p/2}\dd\right)^{1/2}
 \left(\int_{\Ball_{3/4}}
 (A\nabla u\cdot\nabla u)^{(p-2)/2}\dd\right)^{1/2}\notag\\*
 &\quad\le C\left(\int_{\Ball_1}
 (A\nabla u\cdot\nabla u)^{p/2}\dd\right)^{(p-1)/p}.
\end{align}
For \(p=2\), the second factor is simply \(|\Ball_{3/4}|^{1/2}\).
Take \(\eta\in C_c^\infty(\Ball_{3/4})\), equal to one on
\(\Ball_{1/2}\), with \(0\le\eta\le1\) and \(|\nabla\eta|\le C\).
Applying \eqref{N-absolute-aniso-beta} with \(\beta=p-2\),
\eqref{M-L1-pge2}, and H\"older's inequality to the remaining integral of
\((A\nabla u\cdot\nabla u)^{(p-1)/2}\) gives the same bound for the
weighted integral of \(|N|\) on \(\Ball_{1/2}\).
Thus \(|D^2u|\le|M|+|N|\) proves \eqref{hessian-p-ge-2}.

For \(1<p<2\), the same inequalities and \(0<2-p<p\) give
\begin{align}\label{M-L1-plt2}
 &\int_{\Ball_{3/4}\cap\{|\nabla u|>0\}}|M|\dd\notag\\*
 &\quad\le C\left(\int_{\Ball_1}
 (A\nabla u\cdot\nabla u)^{p/2}\dd\right)^{1/2}
 \left(\int_{\Ball_{3/4}}
 (A\nabla u\cdot\nabla u)^{(2-p)/2}\dd\right)^{1/2}\notag\\*
 &\quad\le C\left(\int_{\Ball_1}
 (A\nabla u\cdot\nabla u)^{p/2}\dd\right)^{1/p}.
\end{align}
Using the same \(\eta\) in \eqref{N-absolute-aniso-beta} with \(\beta=0\)
and H\"older's inequality for the remaining integral of
\((A\nabla u\cdot\nabla u)^{1/2}\) gives the same bound for
\(\int_{\Ball_{1/2}\cap\{|\nabla u|>0\}}|N|\dd\), proving
\eqref{hessian-p-lt-2}.
The Appendix gives \(u\in W^{2,p}_{\rm loc}\) for \(1<p<2\) and
\(u\in W^{2,2}_{\rm loc}\) for \(p=2\), so Sobolev--Stampacchia
\cite[Lemma~7.7]{gilbarg-trudinger} gives
\(D^2u=0\) a.e. on \(\{\nabla u=0\}\).
For \(p>2\), only the noncritical Hessian is used.
If the outer energy is zero, \(\nabla u\equiv0\) and the conclusions
hold directly.
\end{proof}

\begin{cor}[]\label{cor:local-pplus1-hessian}\label{lem:pplus-hessian}
Under the assumptions of Theorem~\ref{thm:hessian-estimate}, let
\(B_{2r}(z)\Subset\Ball_1\) and \(0<r\le1\). Then
\begin{equation}\label{local-pplus1-hessian}
 \int_{\Ball_r(z)\cap\{|\nabla u|>0\}}
 (A\nabla u\cdot\nabla u)^{(p-1)/2}|D^2u|\dd
 \le \frac{C}{r}\int_{\Ball_{2r}(z)}(A\nabla u\cdot\nabla u)^{p/2}\dd .
\end{equation}
Here \(C\) depends only on \(n,p,\lambda,\Lambda\).
\end{cor}

\begin{proof}
Choose \(\varphi\in C_c^\infty(\Ball_{2r}(z))\) with
\(0\le\varphi\le1\), \(\varphi\equiv1\) on \(\Ball_r(z)\), and
\(|\nabla\varphi|\le C/r\). Since \(r\le1\),
\eqref{M-direct-estimate} and Cauchy's inequality give
\begin{align}
 &\int_{\Ball_{2r}(z)\cap\{|\nabla u|>0\}}
 (A\nabla u\cdot\nabla u)^{(p-2)/2}|M|^2\varphi^2\dd\notag\\*
 &\quad\le\frac{C}{r^2}\int_{\Ball_{2r}(z)}
 (A\nabla u\cdot\nabla u)^{p/2}\dd ,
 \label{local-M-L2}\\*
 &\int_{\Ball_{2r}(z)\cap\{|\nabla u|>0\}}
 (A\nabla u\cdot\nabla u)^{(p-1)/2}|M|\varphi^2\dd\notag\\*
 &\quad\le\frac{C}{r}\int_{\Ball_{2r}(z)}
 (A\nabla u\cdot\nabla u)^{p/2}\dd .
 \label{local-M-L1}
\end{align}
The second estimate pairs the square root of the first weighted Hessian
integral with
\[
 \left(\int_{\Ball_{2r}(z)}
 (A\nabla u\cdot\nabla u)^{p/2}\varphi^2\dd\right)^{1/2}.
\]
Now \(p-1\ge\max\{0,(p-2)/2\}\), so
\eqref{N-absolute-aniso-beta} applies with \(\beta=p-1\) and
\(\eta=\varphi\), with all its critical-set limits already justified.
Using \eqref{local-M-L1}, \(\delta_0\le1\), and \(r\le1\), it yields
\begin{equation}\label{local-N-L1}
 \int_{\Ball_{2r}(z)\cap\{|\nabla u|>0\}}
 (A\nabla u\cdot\nabla u)^{(p-1)/2}|N|\varphi^2\dd
 \le\frac{C}{r}\int_{\Ball_{2r}(z)}
 (A\nabla u\cdot\nabla u)^{p/2}\dd .
\end{equation}
Indeed, its remaining terms are bounded by
\(C(\delta_0+r^{-1})\int_{\Ball_{2r}(z)}
(A\nabla u\cdot\nabla u)^{p/2}\dd\).
Combining the last two bounds and using \(\varphi=1\) on \(\Ball_r(z)\)
proves \eqref{local-pplus1-hessian}. All constants are independent of
\(z,r,u\).
\end{proof}

\begin{lem}[]\label{lem:direct-anisotropic-higher}
Under the assumptions of Theorem \ref{thm:hessian-estimate}, there exist
\(\gamma>0\) and \(C>0\), depending only on
\(n,p,\lambda,\Lambda\) and on an upper bound for
\(\|A\|_{W^{1,\infty}(\Ball_1)}\), such that
\begin{equation}\label{higher-final}
 \left(\int_{\Ball_{1/2}}(A\nabla u\cdot\nabla u)^{(p+\gamma)/2}\dd\right)^{1/(p+\gamma)}
 \le C\left(\int_{\Ball_1}(A\nabla u\cdot\nabla u)^{p/2}\dd\right)^{1/p}.
\end{equation}
\end{lem}

\begin{proof}
Let
\[
 E=\left(\int_{\Ball_1}(A\nabla u\cdot\nabla u)^{p/2}\dd\right)^{1/p}<\infty,
\]
where finiteness follows from the energy-solution class.  If \(E=0\),
then \(\nabla u=0\) a.e. in \(\Ball_1\) and the conclusion is trivial.  Otherwise,
replacing \(u\) by \(E^{-1}u\) replaces \(f\) by
\(\widetilde f(t)=E^{1-p}f(Et)\ge0\), with
\(\widetilde f'(t)=E^{2-p}f'(Et)\).  Both sides of the stability inequality
are multiplied by \(E^{2-p}\), and the critical set and admissible test
spaces are unchanged.  It therefore suffices to work under
\begin{equation}\label{normalization}
 \int_{\Ball_1}(A\nabla u\cdot\nabla u)^{p/2}\dd\le1.
\end{equation}

\medskip
\noindent
Set \(G=(A\nabla u\cdot\nabla u)^{(p-1)/2}A\nabla u\).
On \(\{|\nabla u|>0\}\), noncritical \(W^{2,2}_{\rm loc}\) regularity,
the Sobolev chain and product rules, and ellipticity give
\[
 |DG|\le C(A\nabla u\cdot\nabla u)^{(p-1)/2}|D^2u|
       +C|DA|(A\nabla u\cdot\nabla u)^{p/2}.
\]
For every compact \(K\Subset\Ball_1\), Cauchy's inequality and the
preliminary weighted Hessian regularity yield
\begin{align*}
 &\int_{K\cap\{|\nabla u|>0\}}
 (A\nabla u\cdot\nabla u)^{(p-1)/2}|D^2u|\dd\\*
 &\quad\le
 \left(\int_{K\cap\{|\nabla u|>0\}}
 (A\nabla u\cdot\nabla u)^{(p-2)/2}|D^2u|^2\dd\right)^{1/2}
 \left(\int_K(A\nabla u\cdot\nabla u)^{p/2}\dd\right)^{1/2}<\infty.
\end{align*}
The coefficient term is also locally integrable.
For the cutoff \(\chi_\tau\) from Theorem~\ref{thm:hessian-estimate},
\(G\chi_\tau\) vanishes near \(\{\nabla u=0\}\), so gluing the
noncritical calculation with zero gives its local weak derivatives.
On \(E_\tau=\{\tau<|\nabla u|<2\tau\}\),
\[
 |G|\,|\nabla\chi_\tau|
 \le \frac C\tau(A\nabla u\cdot\nabla u)^{p/2}|D^2u|
 \le C(A\nabla u\cdot\nabla u)^{(p-1)/2}|D^2u|.
\]
Thus the additional derivative term tends to zero in \(L^1_{\rm loc}\).
Dominated convergence gives \(G\chi_\tau\to G\) and convergence of
\(\chi_\tau DG\) to the zero extension of the noncritical derivative,
both in \(L^1_{\rm loc}\). Passing to the limit in the weak derivative
identity therefore gives
\[
 G\in W^{1,1}_{\rm loc}(\Ball_1;\mathbb R^n),\qquad
 \operatorname{div}G\in L^1_{\rm loc}(\Ball_1).
\]
Stampacchia's property now gives \(DG=0\) a.e. on \(\{\nabla u=0\}\).
In the identities below, extend each weighted Hessian term and each
coefficient-derivative term, taken as a whole, by zero on the critical set.
These extended terms satisfy the preceding \(L^1_{\rm loc}\) bounds, while
\((A\nabla u\cdot\nabla u)^{1/2}f(u)\) is locally bounded and vanishes
on the critical set. The noncritical identities thus hold a.e. and
in distributions on \(\Ball_1\).

For \(\eta\in C_c^\infty(\Ball_1)\), we prove
\begin{equation}\label{div-flux-target}
\begin{aligned}
 &\int_{\Ball_1}\left|
 \operatorname{div}\left((A\nabla u\cdot\nabla u)^{(p-1)/2}A\nabla u\right)
 \right|\eta^2\dd \\*
 &\hspace{25mm}\le
 C\int_{\Ball_1}(A\nabla u\cdot\nabla u)^{p/2}
       (|\nabla\eta|^2+\eta^2)\dd .
\end{aligned}
\end{equation}
The product rule and the original equation give
\begin{align}\label{first-div-identity}
 &\operatorname{div}\left((A\nabla u\cdot\nabla u)^{(p-1)/2}A\nabla u\right)\notag\\*
 &=-(A\nabla u\cdot\nabla u)^{1/2}f(u)
   +(A\nabla u\cdot\nabla u)^{(p-1)/2}
     \frac{(A_{i\ell}u_\ell)(A_{jm}u_m)u_{ij}}{A\nabla u\cdot\nabla u}\notag\\*
 &\quad+\frac12(A\nabla u\cdot\nabla u)^{(p-3)/2}
     A_{ij}u_jA_{mn,i}u_mu_n .
\end{align}
Expanding the original equation and isolating the normal contraction yields
\begin{align}\label{normal-aniso-isolation}
 &(p-1)(A\nabla u\cdot\nabla u)^{(p-1)/2}
 \frac{(A_{i\ell}u_\ell)(A_{jm}u_m)u_{ij}}{A\nabla u\cdot\nabla u}\notag\\*
 &=-(A\nabla u\cdot\nabla u)^{1/2}f(u)\notag\\*
 &\quad-(A\nabla u\cdot\nabla u)^{(p-1)/2}
 \left(A_{ij}u_{ij}-\frac{(A_{i\ell}u_\ell)(A_{jm}u_m)u_{ij}}{A\nabla u\cdot\nabla u}\right)\notag\\*
 &\quad-\frac{p-2}{2}(A\nabla u\cdot\nabla u)^{(p-3)/2}
       A_{ij}u_jA_{mn,i}u_mu_n
 -(A\nabla u\cdot\nabla u)^{(p-1)/2}A_{ij,i}u_j .
\end{align}
Combining these identities, including both coefficient terms, gives
\begin{align}\label{div-flux-identity}
 &(p-1)\operatorname{div}\left((A\nabla u\cdot\nabla u)^{(p-1)/2}A\nabla u\right)\notag\\*
 &=-p(A\nabla u\cdot\nabla u)^{1/2}f(u)\notag\\*
 &\quad-(A\nabla u\cdot\nabla u)^{(p-1)/2}
 \left(A_{ij}u_{ij}-\frac{(A_{i\ell}u_\ell)(A_{jm}u_m)u_{ij}}{A\nabla u\cdot\nabla u}\right)\notag\\*
 &\quad+\frac12(A\nabla u\cdot\nabla u)^{(p-3)/2}
       A_{ij}u_jA_{mn,i}u_mu_n
 -(A\nabla u\cdot\nabla u)^{(p-1)/2}A_{ij,i}u_j .
\end{align}
By ellipticity, the absolute value of each coefficient term here, and of
the corresponding terms in \eqref{first-div-identity}--\eqref{normal-aniso-isolation},
is at most
\(C\|DA\|_{L^\infty(\Ball_1)}(A\nabla u\cdot\nabla u)^{p/2}\).

Integrating \eqref{first-div-identity} against \(\eta^2\) and using
\(G\in W^{1,1}_{\rm loc}\) gives
\begin{align}\label{f-term-start}
 &\int_{\Ball_1}(A\nabla u\cdot\nabla u)^{1/2}f(u)\eta^2\dd\notag\\*
 &=\int_{\Ball_1}(A\nabla u\cdot\nabla u)^{(p-1)/2}
       A\nabla u\cdot\nabla(\eta^2)\dd\notag\\*
 &\quad+\int_{\Ball_1}(A\nabla u\cdot\nabla u)^{(p-1)/2}
       \frac{(A_{i\ell}u_\ell)(A_{jm}u_m)u_{ij}}{A\nabla u\cdot\nabla u}\eta^2\dd\notag\\*
 &\quad+\frac12\int_{\Ball_1}(A\nabla u\cdot\nabla u)^{(p-3)/2}
       A_{ij}u_jA_{mn,i}u_mu_n\eta^2\dd .
\end{align}
Since \(f\ge0\) and \(p-1>0\), \eqref{normal-aniso-isolation} bounds
the signed normal integral from above by the absolute tangential
integral and the coefficient error. Using the preceding coefficient
bounds and \(2|\eta||\nabla\eta|\le\eta^2+|\nabla\eta|^2\), we obtain
\begin{align}\label{f-term-final}
 &\int_{\Ball_1}(A\nabla u\cdot\nabla u)^{1/2}f(u)\eta^2\dd\notag\\*
 &\le C\int_{\Ball_1}(A\nabla u\cdot\nabla u)^{(p-1)/2}
 \left|A_{ij}u_{ij}-\frac{(A_{i\ell}u_\ell)(A_{jm}u_m)u_{ij}}{A\nabla u\cdot\nabla u}\right|\eta^2\dd\notag\\*
 &\quad+C\int_{\Ball_1}(A\nabla u\cdot\nabla u)^{p/2}
       (|\nabla\eta|^2+\eta^2)\dd .
\end{align}
Taking absolute values in \eqref{div-flux-identity} and using
\eqref{f-term-final} leaves only the tangential integral to control.
Since
\[
 \left(A_{ij}-\frac{(A_{i\ell}u_\ell)(A_{jm}u_m)}{A\nabla u\cdot\nabla u}\right)e_i e_j
 =Ae\cdot e-\frac{(A\nabla u\cdot e)^2}{A\nabla u\cdot\nabla u}=0,
\]
the decomposition \(D^2u=M+Ne\otimes e\) and ellipticity give
\begin{equation}\label{tangential-Hessian-bound}
 \left|A_{ij}u_{ij}-\frac{(A_{i\ell}u_\ell)(A_{jm}u_m)u_{ij}}{A\nabla u\cdot\nabla u}\right|
 \le C|M|.
\end{equation}
Cauchy's inequality and \eqref{M-direct-estimate} yield
\begin{align}\label{M-close-div}
 &\int_{\Ball_1}(A\nabla u\cdot\nabla u)^{(p-1)/2}|M|\eta^2\dd\notag\\*
 &\le \left(\int_{\Ball_1}(A\nabla u\cdot\nabla u)^{(p-2)/2}|M|^2\eta^2\dd\right)^{1/2}
       \left(\int_{\Ball_1}(A\nabla u\cdot\nabla u)^{p/2}\eta^2\dd\right)^{1/2}\notag\\*
 &\le C\int_{\Ball_1}(A\nabla u\cdot\nabla u)^{p/2}
       (|\nabla\eta|^2+\eta^2)\dd .
\end{align}
Together with \eqref{tangential-Hessian-bound}, this proves the uniform
reaction bound and \eqref{div-flux-target}.

\medskip
\noindent
Fix \(\eta\in C_c^\infty(\Ball_1)\) with \(0\le\eta\le1\),
\(\eta\equiv1\) on \(\Ball_{1/2}\), and \(|\nabla\eta|\le C\).  For \(t\in\R\), define
\begin{equation}\label{Phi-def}
 \Phi(t)=\int_{\{u=t\}\cap\Ball_1\cap\{|\nabla u|>0\}}
 \frac{(A\nabla u\cdot\nabla u)^{(p+1)/2}}{|\nabla u|}
 \eta^2\,d{\mathcal H}^{n-1}.
\end{equation}
Coarea and local \(C^1\) regularity give the qualitative integrability
\[
 \int_\R\Phi(t)\,dt
 =\int_{\Ball_1}(A\nabla u\cdot\nabla u)^{(p+1)/2}\eta^2\dd<\infty.
\]
Choose smooth nondecreasing cutoffs \(K_\varepsilon\) satisfying
\[
 \begin{gathered}
 K_\varepsilon(s)=0\ (s\le0),\qquad
 K_\varepsilon(s)=1\ (s\ge\varepsilon),\\
 0\le K'_\varepsilon\le C/\varepsilon,\qquad
 \int_\R K'_\varepsilon=1.
 \end{gathered}
\]
For fixed \(\varepsilon\), \(K_\varepsilon(u-t)\) is bounded and locally
Lipschitz near \(\operatorname{supp}\eta\).  Since
\(G\eta^2\in W^{1,1}_0(\Ball_1;\mathbb R^n)\), integration by parts and coarea give
\begin{align}\label{heaviside-identity}
 &-\int_{\Ball_1}K_\varepsilon(u-t)
 \operatorname{div}\left((A\nabla u\cdot\nabla u)^{(p-1)/2}A\nabla u\,\eta^2\right)\dd\notag\\*
 &=\int_{\Ball_1}K'_\varepsilon(u-t)
       (A\nabla u\cdot\nabla u)^{(p+1)/2}\eta^2\dd
 =\int_\R K'_\varepsilon(s-t)\Phi(s)\,ds.
\end{align}
As \(\varepsilon\downarrow0\), the right-hand side tends to \(\Phi(t)\)
at its Lebesgue points, while the left-hand side converges by dominated
convergence, using the integrable absolute divergence and \(0\le K_\varepsilon\le1\).
Thus, for a.e. \(t\),
\begin{equation}\label{Phi-div-identity}
 \Phi(t)=-\int_{\{u>t\}\cap\Ball_1}
 \operatorname{div}\left((A\nabla u\cdot\nabla u)^{(p-1)/2}A\nabla u\,\eta^2\right)\dd .
\end{equation}
Expanding this divergence and using ellipticity,
\eqref{div-flux-target}, and \eqref{normalization}, we obtain
\begin{equation}\label{Phi-bound}
 \Phi(t)\le C\qquad\hbox{for a.e. }t\in\R.
\end{equation}

Ellipticity, \eqref{normalization}, and Sobolev--Poincar\'e give
\begin{equation}\label{sobolev-poincare}
 \|u-\bar u\|_{L^m(\Ball_1)}\le C,\qquad
 \bar u=|\Ball_1|^{-1}\int_{\Ball_1}u\dd,
\end{equation}
where \(m=np/(n-p)\) if \(p<n\), and \(m=2p\) if \(p\ge n\). Set
\begin{equation}\label{h-def}
 h(t)=\max\{1,|t-\bar u|\}.
\end{equation}
Coarea and \eqref{sobolev-poincare} give
\begin{equation}\label{coarea-hm}
 \int_\R\int_{\{u=t\}\cap\Ball_1\cap\{|\nabla u|>0\}}
 \frac{h(t)^m}{|\nabla u|}\,d{\mathcal H}^{n-1}dt
 =\int_{\Ball_1\cap\{|\nabla u|>0\}}h(u(x))^m\dd\le C.
\end{equation}
Fix \(q\in(1,m/p)\) depending only on \(n,p\), and set
\begin{equation}\label{theta-beta-def}
 \theta=\frac{q}{m+q},\qquad
 \beta=(p+1)(1-\theta)=\frac{(p+1)m}{m+q}.
\end{equation}
Then \(\gamma:=\beta-p=(m-pq)/(m+q)>0\).
Since \(m\theta=q(1-\theta)\), on \(\{u=t\}\cap\{|\nabla u|>0\}\) we have
\begin{equation}\label{factorization}
 \frac{(A\nabla u\cdot\nabla u)^{\beta/2}}{|\nabla u|}
 =\left(\frac{h(t)^m}{|\nabla u|}\right)^\theta
  \left(h(t)^{-q}\frac{(A\nabla u\cdot\nabla u)^{(p+1)/2}}{|\nabla u|}\right)^{1-\theta}.
\end{equation}
Coarea and H\"older's inequality with exponents \(1/\theta\), \(1/(1-\theta)\) give
\begin{align}\label{higher-coarea-holder}
 &\int_{\Ball_{1/2}}(A\nabla u\cdot\nabla u)^{\beta/2}\dd\notag\\*
 &\le\left(\int_\R\int_{\{u=t\}\cap\Ball_1\cap\{|\nabla u|>0\}}
       \frac{h(t)^m}{|\nabla u|}\,d{\mathcal H}^{n-1}dt\right)^\theta\notag\\*
 &\quad\times\left(\int_\R h(t)^{-q}
       \int_{\{u=t\}\cap\Ball_{1/2}\cap\{|\nabla u|>0\}}
       \frac{(A\nabla u\cdot\nabla u)^{(p+1)/2}}{|\nabla u|}
       \,d{\mathcal H}^{n-1}dt\right)^{1-\theta}.
\end{align}
Both factors are bounded by \eqref{coarea-hm}, \eqref{Phi-bound},
\(\eta=1\) on \(\Ball_{1/2}\), and
\mbox{\(\int_\R h(t)^{-q}dt=2+2/(q-1)\)}, independently of \(\bar u\).
Truncation and Tonelli justify the nonnegative integrals without assuming
the left-hand side finite. Thus
\begin{equation}\label{higher-normalized}
 \int_{\Ball_{1/2}}(A\nabla u\cdot\nabla u)^{(p+\gamma)/2}\dd\le C.
\end{equation}
Rescaling by \(E\) gives \eqref{higher-final}.
\end{proof}

\section{Control of the full gradient by the weighted radial derivative}\label{sec:anisotropic-lemma41}

In this section we show that, under a doubling assumption on the
gradient energy, the full gradient energy in a ball is controlled by
the weighted radial energy in its outer annulus.

Let $n\ge2$ and $p>1$. Throughout this section, $A$ is symmetric and satisfies
\[
 A\in W^{1,\infty}_{\rm loc}(B_2;\mathbb R^{n\times n}),
 \qquad
 \lambda |\xi|^2\le A(x)\xi\cdot\xi\le\Lambda |\xi|^2,
 \qquad
 \|DA\|_{L^\infty(B_2)}\le\delta_0,
\]
where $0<\lambda\le\Lambda$ and $\delta_0$ is the small constant fixed in
Section~\ref{sec:hessian-estimate}. We use locally Lipschitz representatives
and weak derivatives, with the norm from Section~2.

For a stable regular solution $v$ of (1.1) in $B_2$ with nonnegative
$C^1$ nonlinearity, fixed-factor rescaling and a finite covering of $B_{3/2}$
by balls whose doubles lie compactly in $B_2$ allow us to apply
Lemma~\ref{lem:direct-anisotropic-higher} and Theorem~\ref{thm:hessian-estimate}:
\begin{equation}\label{sec4:higher-input}
 \|\nabla v\|_{L^{p+\gamma}(B_{3/2})}
 \le C\|\nabla v\|_{L^p(B_2)}
\end{equation}
\begin{equation}\label{sec4:hessian-input-ge2}
 \int_{B_{3/2}\cap\{|\nabla v|>0\}}
 (A\nabla v\cdot\nabla v)^{(p-2)/2}|D^2v|\,dx
 \le C\|\nabla v\|_{L^p(B_2)}^{p-1}
 \qquad (p\ge2)
\end{equation}
\begin{equation}\label{sec4:hessian-input-lt2}
 \int_{B_{3/2}\cap\{|\nabla v|>0\}}|D^2v|\,dx
 \le C\|\nabla v\|_{L^p(B_2)}
 \qquad (1<p<2).
\end{equation}
Here $C$ and $\gamma>0$ depend only on $n,p,\lambda,\Lambda$.
For $p>2$, we use the noncritical Hessian with its weighted zero extension
from Corollary~\ref{cor:dq-weighted-sobolev}. For $1<p\le2$,
Theorem~\ref{dq:main-theorem} and Sobolev--Stampacchia identify the zero
extension with the ordinary weak Hessian.

For symmetric uniformly elliptic $B\in W^{1,\infty}_{\rm loc}$, the field
\[
 \mathcal A_B(x,\xi)
 :=(B(x)\xi\cdot\xi)^{p/2-1}B(x)\xi
\]
is Carath\'eodory, with $p$-growth, coercivity, and strict monotonicity
controlled by $p,\lambda,\Lambda$. It also satisfies the homogeneity condition
\[
 \mathcal A_B(x,t\xi)
 =|t|^{p-2}t\,\mathcal A_B(x,\xi)
 \qquad (t\in\mathbb R,\ \xi\in\mathbb R^n),
\]
where the scalar factor is defined to be zero at $t=0$.
Weak comparison, Dirichlet existence
and uniqueness for the homogeneous equation, weak Harnack for nonnegative
supersolutions, and the strong minimum principle therefore apply to the
limiting operator with locally Lipschitz coefficient $A_\infty$; we use the
unweighted theory of \cite[Chapters~3 and~7]{hkm}, with the specific results
identified below.

The next lemma uses the compactness argument in \cite[Lemma~A.2]{cms}
for the transformed variables.

\begin{lemma}[]\label{lem:sec4-nonlinear-rellich}
Let $p\ge2$ and let $V\subset\mathbb R^n$ be a bounded Lipschitz domain.
Suppose that $w_k\in L^{p-1}(V)$ and
\[
 z_k:=|w_k|^{p-2}w_k\in W^{1,1}(V),
 \qquad
 \sup_k\bigl(\|z_k\|_{L^1(V)}+\|Dz_k\|_{L^1(V)}\bigr)<\infty.
\]
Then, after passing to a subsequence, there is $w\in L^{p-1}(V)$ such that
\[
 w_k\longrightarrow w
 \qquad\hbox{strongly in }L^{p-1}(V).
\]
If, in addition, $w_k\to w_0$ in the sense of distributions, then $w=w_0$.
\end{lemma}

\begin{proof}
The compact embedding $W^{1,1}(V)\Subset L^1(V)$ yields, after passing to a
subsequence,
\[
 z_k\longrightarrow z
 \qquad\hbox{strongly in }L^1(V).
\]
Let $p'=p/(p-1)$ and set
\[
 w:=|z|^{p'-2}z,
\]
with the value zero at $z=0$.  The inverse-monotonicity estimate
\begin{equation}\label{sec4:inverse-monotonicity}
 |a-b|^{p-1}
 \le C_p\left||a|^{p-2}a-|b|^{p-2}b\right|
 \qquad (a,b\in\mathbb R)
\end{equation}
gives
\[
 \|w_k-w\|_{L^{p-1}(V)}^{p-1}
 \le C_p\|z_k-z\|_{L^1(V)}\longrightarrow0.
\]
The last assertion follows because strong $L^{p-1}$ convergence implies
convergence in distributions and distributional limits are unique.
\end{proof}

\begin{lemma}[]
\label{thm:anisotropic-lemma41}
Let $n\ge2$, $p>1$, and $\delta\in(0,1)$.  Let
$A=A(x)\in W^{1,\infty}_{\rm loc}(B_2;\mathbb R^{n\times n})$ be symmetric and satisfy
\[
 \lambda |\xi|^2\le A(x)\xi\cdot\xi\le\Lambda |\xi|^2,
 \qquad
 \|DA\|_{L^\infty(B_2)}\le\delta_0
\]
for every $x\in B_2$ and every $\xi\in\mathbb R^n$.  Let $u$ be a stable
regular solution of
\[
 -\operatorname{div}\left(
 (A(x)\nabla u\cdot\nabla u)^{p/2-1}A(x)\nabla u\right)=f(u)
 \quad\hbox{in }B_2,
\]
where $f\in C^1(\mathbb R)$ and $f\ge0$.  Assume that
\[
 \int_{B_1}(A(x)\nabla u\cdot\nabla u)^{p/2}\,dx
 \ge\delta
 \int_{B_2}(A(x)\nabla u\cdot\nabla u)^{p/2}\,dx.
\]
Then there exists $C_\delta=C(n,p,\lambda,\Lambda,\delta)$ such that
\[
 \int_{B_{3/2}}(A(x)\nabla u\cdot\nabla u)^{p/2}\,dx
 \le C_\delta
 \int_{B_{3/2}\setminus B_1}
 (A(x)\nabla u\cdot\nabla u)^{p/2-1}
 \left(\frac{x}{|x|}\cdot\nabla u\right)^2\,dx.
\]
When $1<p<2$, the integrand on the right is defined to be zero at
$\nabla u=0$.  This is its continuous extension, since ellipticity gives
\[
 0\le (A\nabla u\cdot\nabla u)^{p/2-1}u_r^2
 \le C|\nabla u|^p.
\]
\end{lemma}

\begin{proof}
If the conclusion failed, there would be sequences $A_k,u_k,f_k$
satisfying all the hypotheses with the same structural constants and the
same $\delta$, such that
\begin{equation}\label{sec4:counterexample-ratio}
 \frac{
 \displaystyle\int_{B_{3/2}\setminus B_1}
 (A_k\nabla u_k\cdot\nabla u_k)^{p/2-1}
 \left(\frac{x}{|x|}\cdot\nabla u_k\right)^2\,dx
 }{
 \displaystyle\int_{B_{3/2}}
 (A_k\nabla u_k\cdot\nabla u_k)^{p/2}\,dx
 }
 \longrightarrow0.
\end{equation}
The denominator is positive, since otherwise ellipticity would give
$\nabla u_k=0$ a.e. in $B_{3/2}$ and the desired inequality would be trivial.
Set
\[
 m_k=(u_k)_{B_2},
 \qquad
 E_k=\left(\int_{B_{3/2}}
 (A_k\nabla u_k\cdot\nabla u_k)^{p/2}\,dx\right)^{1/p}.
\]
Then $v_k=(u_k-m_k)/E_k$ solves the same equation with coefficient $A_k$
and nonlinearity
\[
 \widetilde f_k(t)=E_k^{1-p}f_k(m_k+E_kt)\in C^1(\mathbb R),
 \qquad
 \widetilde f_k\ge0,
 \qquad
 \widetilde f_k'(t)=E_k^{2-p}f_k'(m_k+E_kt).
\]
The critical set is unchanged, and
\[
 (A_k\nabla v_k\cdot\nabla v_k)^{p/2-1}
 =E_k^{2-p}(A_k\nabla u_k\cdot\nabla u_k)^{p/2-1}.
\]
Thus, for each $k$, the stability test spaces coincide and have equivalent
norms. Both principal terms in the stability inequality and the reaction
term are multiplied by $E_k^{2-p}$, so $v_k$ is again a stable regular solution.

By construction,
\begin{equation}\label{sec4:normalized-energy}
 \int_{B_{3/2}}(A_k\nabla v_k\cdot\nabla v_k)^{p/2}\,dx=1
\end{equation}
and
\begin{equation}\label{sec4:normalized-radial}
 \int_{B_{3/2}\setminus B_1}
 (A_k\nabla v_k\cdot\nabla v_k)^{p/2-1}
 \left(\frac{x}{|x|}\cdot\nabla v_k\right)^2\,dx
 \longrightarrow0.
\end{equation}
The doubling condition, \eqref{sec4:normalized-energy}, and ellipticity give
\[
 \lambda^{p/2}\int_{B_2}|\nabla v_k|^p\,dx
 \le\int_{B_2}(A_k\nabla v_k\cdot\nabla v_k)^{p/2}\,dx
 \le\delta^{-1}.
\]
Since $(v_k)_{B_2}=0$, Poincar\'e's inequality
\cite[(7.45)]{gilbarg-trudinger} yields
\begin{equation}\label{sec4:W1p-bound}
 \|v_k\|_{W^{1,p}(B_2)}\le C.
\end{equation}
After passing to a subsequence,
\begin{equation}\label{sec4:weak-convergence}
 v_k\rightharpoonup v_\infty\quad\hbox{weakly in }W^{1,p}(B_2),
 \qquad
 v_k\longrightarrow v_\infty\quad\hbox{strongly in }L^p(B_2).
\end{equation}

By ellipticity and the derivative bound, Arzel\`a--Ascoli and weak-star
compactness give, after a further subsequence,
\begin{equation}\label{sec4:coefficient-convergence}
 A_k\longrightarrow A_\infty
 \quad\hbox{uniformly on compact subsets of }B_2,
 \qquad
 A_k\rightharpoonup^*A_\infty
 \quad\hbox{in }W^{1,\infty}_{\rm loc}(B_2).
\end{equation}
The limit is symmetric, with the same ellipticity bounds and
$\|DA_\infty\|_{L^\infty(B_2)}\le\delta_0$.

From \eqref{sec4:higher-input} and \eqref{sec4:W1p-bound},
\begin{equation}\label{sec4:higher-uniform}
 \|\nabla v_k\|_{L^{p+\gamma}(B_{3/2})}\le C.
\end{equation}
Weak compactness and distributional convergence in
\eqref{sec4:weak-convergence} identify the weak $L^{p+\gamma}(B_{3/2})$
limit of a further subsequence as $\nabla v_\infty$; in particular,
\begin{equation}\label{sec4:limit-higher}
 \nabla v_\infty\in L^{p+\gamma}(B_{3/2}).
\end{equation}

Assume first that $p\ge2$.  By \eqref{sec4:hessian-input-ge2},
ellipticity, and \eqref{sec4:W1p-bound},
\begin{equation}\label{sec4:weighted-Hessian-L1}
 \int_{B_{3/2}}|\nabla v_k|^{p-2}|D^2v_k|\,dx\le C.
\end{equation}
For $p>2$, all Hessian terms below mean noncritical derivatives extended
by zero.  Fix a bounded Lipschitz domain $V\Subset B_{3/2}$ and a component
$i$, and put
\[
 g_k=\nabla v_k,\qquad w_{k,i}=\partial_i v_k,\qquad
 z_{k,i}=|w_{k,i}|^{p-2}w_{k,i}.
\]
For $p=2$, Theorem~\ref{dq:main-theorem} supplies the ordinary weak Hessian,
so \eqref{sec4:weighted-Hessian-L1} and \eqref{sec4:W1p-bound} give
$z_{k,i}=w_{k,i}\in W^{1,1}(V)$ with a uniform bound.

For $p>2$, choose $\chi\in C^\infty([0,\infty))$ with $0\le\chi\le1$,
$\chi=0$ on $[0,1]$ and $\chi=1$ on $[2,\infty)$, and set
\[
 \chi_{k,\tau}=\chi\left(\frac{|g_k|}{\tau}\right),\qquad
 g_{k,\tau}=\chi_{k,\tau}g_k,\qquad
 z_{k,i,\tau}=|(g_{k,\tau})_i|^{p-2}(g_{k,\tau})_i.
\]
The cutoff construction in Corollary~\ref{cor:dq-weighted-sobolev} gives
$g_{k,\tau}\in W^{1,2}(V;\mathbb R^n)$ and
\[
 |D\chi_{k,\tau}|\le\frac C\tau|Dg_k|
 \mathbf1_{\{\tau<|g_k|<2\tau\}}\quad\hbox{in }V.
\]
Since $g_k$ is bounded on $\overline V$, the Sobolev chain and product rules
apply and, using $|g_k|/\tau<2$ on the transition region, yield
\begin{equation}\label{sec4:cutoff-derivative-bound}
 \begin{aligned}
 &\bigl\|Dz_{k,i,\tau}-(p-1)|w_{k,i}|^{p-2}Dw_{k,i}\bigr\|_{L^1(V)}\\
 &\qquad\le C\int_{V\cap\{0<|g_k|<2\tau\}}
 |g_k|^{p-2}|Dg_k|\,dx\longrightarrow0
 \qquad(\tau\downarrow0).
 \end{aligned}
\end{equation}
Here $k$ is fixed, and the convergence follows from the integrability in
\eqref{sec4:weighted-Hessian-L1}.  Also,
$z_{k,i,\tau}\to z_{k,i}$ in $L^1(V)$ by dominated convergence, since
$|z_{k,i,\tau}|+|z_{k,i}|\le2|g_k|^{p-1}$.
Passing to the limit in the weak derivative identity gives
\begin{equation}\label{sec4:z-W11}
 z_{k,i}\in W^{1,1}(V),\qquad
 Dz_{k,i}=(p-1)|w_{k,i}|^{p-2}Dw_{k,i}
 \quad\hbox{on }\{|g_k|>0\},
\end{equation}
with zero extension on $\{|g_k|=0\}$, and
\begin{equation}\label{sec4:z-uniform-W11}
 \|z_{k,i}\|_{W^{1,1}(V)}
 \le C\left(\int_V|g_k|^{p-1}\,dx
 +\int_V|g_k|^{p-2}|Dg_k|\,dx\right)\le C.
\end{equation}
Thus, for every $p\ge2$, Lemma~\ref{lem:sec4-nonlinear-rellich},
identification of the limit by weak $W^{1,p}$ convergence, and a diagonal
subsequence over the finitely many components and a Lipschitz exhaustion
of $B_{3/2}$ give
\begin{equation}\label{sec4:strong-pminus1}
 \nabla v_k\longrightarrow\nabla v_\infty
 \qquad\hbox{strongly in }L^{p-1}_{\rm loc}(B_{3/2}).
\end{equation}
Interpolating with the uniform $L^{p+\gamma}$ bound
\eqref{sec4:higher-uniform} yields
\begin{equation}\label{sec4:strong-Lp-ge2}
 \nabla v_k\longrightarrow\nabla v_\infty
 \qquad\hbox{strongly in }L^p_{\rm loc}(B_{3/2})
 \qquad (p\ge2).
\end{equation}

If $1<p<2$, Theorem~\ref{dq:main-theorem} gives
$v_k\in W^{2,p}_{\rm loc}(B_2)$, so Sobolev--Stampacchia implies
$D^2v_k=0$ a.e. on $\{\nabla v_k=0\}$.  Thus
\eqref{sec4:hessian-input-lt2} and \eqref{sec4:W1p-bound} yield
\[
 \|\nabla v_k\|_{W^{1,1}(V)}\le C
 \qquad(V\Subset B_{3/2}).
\]
Rellich compactness on interior balls and weak $W^{1,p}$ convergence give,
after passing to a diagonal subsequence,
\[
 \nabla v_k\longrightarrow\nabla v_\infty
 \qquad\hbox{strongly in }L^1_{\rm loc}(B_{3/2}).
\]
Interpolating with \eqref{sec4:higher-uniform} gives
\begin{equation}\label{sec4:strong-Lp-lt2}
 \nabla v_k\longrightarrow\nabla v_\infty
 \qquad\hbox{strongly in }L^p_{\rm loc}(B_{3/2})
 \qquad (1<p<2).
\end{equation}
Together with \eqref{sec4:strong-Lp-ge2}, this proves
\begin{equation}\label{sec4:strong-Lp}
 \nabla v_k\longrightarrow\nabla v_\infty
 \qquad\hbox{strongly in }L^p_{\rm loc}(B_{3/2})
\end{equation}
for every $p>1$.

For symmetric matrices $B,\widetilde B$ with the same ellipticity bounds, the mean
value theorem gives
\begin{align}
 \left|(B\xi\cdot\xi)^{p/2}-(\widetilde B\xi\cdot\xi)^{p/2}\right|
 &\le C|B-\widetilde B||\xi|^p, \label{sec4:matrix-energy}\\
 \left|\mathcal A_B(\xi)-\mathcal A_{\widetilde B}(\xi)\right|
 &\le C|B-\widetilde B||\xi|^{p-1}. \label{sec4:matrix-flux}
\end{align}
Here $\mathcal A_B(\xi)=(B\xi\cdot\xi)^{p/2-1}B\xi$.  For fixed $B$,
\begin{equation}\label{sec4:gradient-energy}
 \left|(B\xi\cdot\xi)^{p/2}-(B\zeta\cdot\zeta)^{p/2}\right|
 \le C(|\xi|+|\zeta|)^{p-1}|\xi-\zeta|.
\end{equation}
By H\"older's inequality, \eqref{sec4:coefficient-convergence} and
\eqref{sec4:strong-Lp} therefore give, for every $V\Subset B_{3/2}$,
\begin{equation}\label{sec4:energy-convergence}
 \int_V(A_k\nabla v_k\cdot\nabla v_k)^{p/2}\,dx
 \longrightarrow
 \int_V(A_\infty\nabla v_\infty\cdot\nabla v_\infty)^{p/2}\,dx.
\end{equation}
Taking $V=B_1$ in \eqref{sec4:energy-convergence} and using the doubling
condition and \eqref{sec4:normalized-energy}, we obtain
\begin{equation}\label{sec4:positive-limit-energy}
 \int_{B_1}(A_\infty\nabla v_\infty\cdot\nabla v_\infty)^{p/2}\,dx
 \ge\delta.
\end{equation}

For the radial energy, fix $\eta\in(0,1/8)$ and set
\[
 K_\eta=B_{3/2-\eta}\setminus\overline{B_{1+\eta}},
 \qquad e(x)=\frac{x}{|x|}.
\]
For every $p>1$, the radial density, defined to be zero at $\xi=0$, satisfies
\begin{equation}\label{sec4:matrix-radial}
 \left|
 (B\xi\cdot\xi)^{p/2-1}(e\cdot\xi)^2
 -(\widetilde B\xi\cdot\xi)^{p/2-1}(e\cdot\xi)^2
 \right|
 \le C|B-\widetilde B||\xi|^p.
\end{equation}
Moreover,
\begin{align}
 &\left|
 (B\xi\cdot\xi)^{p/2-1}(e\cdot\xi)^2
 -(B\zeta\cdot\zeta)^{p/2-1}(e\cdot\zeta)^2
 \right| \notag\\
 &\hspace{35mm}\le C(|\xi|+|\zeta|)^{p-1}|\xi-\zeta|.
 \label{sec4:gradient-radial}
\end{align}
\begin{equation}\label{sec4:radial-domination}
 0\le (B\xi\cdot\xi)^{p/2-1}(e\cdot\xi)^2\le C|\xi|^p.
\end{equation}
The matrix estimate follows by the mean value theorem.  For $\xi\ne0$,
direct differentiation and ellipticity give
\[
 \left|D_\xi\bigl[(B\xi\cdot\xi)^{p/2-1}(e\cdot\xi)^2\bigr]\right|
 \le C|\xi|^{p-1}.
\]
Together with \eqref{sec4:radial-domination}, this gives a $C^1$ extension
at zero with derivative zero, proving \eqref{sec4:gradient-radial}
for all $p>1$.  These estimates, H\"older's inequality,
\eqref{sec4:coefficient-convergence}, and \eqref{sec4:strong-Lp} give strong
$L^1(K_\eta)$ convergence of the radial densities for each fixed $\eta$.
Nonnegativity and \eqref{sec4:normalized-radial} therefore force the
limiting density to vanish a.e. on $K_\eta$.  Letting $\eta\downarrow0$
along a sequence gives the same conclusion on the annulus.
Its weight is positive where $\nabla v_\infty\ne0$,
and the radial derivative is zero where $\nabla v_\infty=0$.  Thus
\begin{equation}\label{sec4:radial-zero}
 \frac{x}{|x|}\cdot\nabla v_\infty=0
 \qquad\hbox{for a.e. }x\in B_{3/2}\setminus B_1.
\end{equation}

It remains to pass to the limit in the equation.  For fixed elliptic $B$,
\begin{equation}\label{sec4:gradient-flux}
 |\mathcal A_B(\xi)-\mathcal A_B(\zeta)|
 \le
 \begin{cases}
 C(|\xi|+|\zeta|)^{p-2}|\xi-\zeta|,&p\ge2,\\
 C|\xi-\zeta|^{p-1},&1<p<2.
 \end{cases}
\end{equation}
Together with \eqref{sec4:matrix-flux},
\eqref{sec4:coefficient-convergence}, and \eqref{sec4:strong-Lp}, this gives
\begin{equation}\label{sec4:flux-convergence}
 \mathcal A_{A_k}(x,\nabla v_k)
 \longrightarrow
 \mathcal A_{A_\infty}(x,\nabla v_\infty)
 \quad\hbox{strongly in }L^{p/(p-1)}_{\rm loc}(B_{3/2}).
\end{equation}
For every nonnegative $\varphi\in C_c^1(B_{3/2})$, the equation for $v_k$
and $\widetilde f_k\ge0$ give
\[
 \int_{B_{3/2}}
 \mathcal A_{A_k}(x,\nabla v_k)\cdot\nabla\varphi\,dx\ge0.
\]
Passing to the limit by \eqref{sec4:flux-convergence},
\begin{equation}\label{sec4:limit-supersolution}
 \int_{B_{3/2}}
 \mathcal A_{A_\infty}(x,\nabla v_\infty)\cdot\nabla\varphi\,dx\ge0
 \qquad(\varphi\ge0).
\end{equation}
Thus $v_\infty$ is a weak supersolution in $B_{3/2}$.

By polar-coordinate Sobolev slicing and \eqref{sec4:radial-zero}, for a.e.
$\theta\in S^{n-1}$ the slice $r\mapsto v_\infty(r\theta)$ belongs to
$W^{1,p}(1,3/2)$ and has zero weak derivative. Thus there is a measurable
$g:S^{n-1}\to\mathbb R$ such that
\begin{equation}\label{sec4:angular-representative}
 v_\infty(r\theta)=g(\theta)
 \quad\hbox{for a.e. }(r,\theta)\in(1,3/2)\times S^{n-1}.
\end{equation}
Hence $v_\infty$ has the same essential infimum on every complete
spherical layer compactly contained in the annulus.

Let $\mathcal A=B_{3/2}\setminus\overline{B_1}$. By
\cite[Theorem~3.17]{hkm}, there is a unique weak solution $h$ of
\[
 -\operatorname{div}\mathcal A_{A_\infty}(x,\nabla h)=0
 \quad\hbox{in }\mathcal A,
 \qquad h-v_\infty\in W^{1,p}_0(\mathcal A).
\]
The comparison principle \cite[Lemma~3.18]{hkm}, together with
\eqref{sec4:limit-supersolution}, gives $v_\infty\ge h$ a.e. in $\mathcal A$.
By \cite[Theorem~3.70]{hkm}, $h$ has a locally continuous representative,
which is bounded below
on $\overline{B_{11/8}}\setminus B_{9/8}$, so
\[
 c_0:=\operatorname*{ess\,inf}_{B_{11/8}\setminus B_{9/8}}v_\infty
 >-\infty.
\]
By \eqref{sec4:angular-representative},
\begin{equation}\label{sec4:annular-nonnegative}
 v_\infty-c_0\ge0\quad\hbox{a.e. in }\mathcal A,
 \qquad
 \operatorname*{ess\,inf}_{B_{11/8}\setminus B_{9/8}}
 (v_\infty-c_0)=0.
\end{equation}
Cover $\overline{B_{11/8}}\setminus B_{9/8}$ by finitely many balls
$B_{\rho/2}(x_j)$ with $B_{4\rho}(x_j)\Subset\mathcal A$.
By \eqref{sec4:annular-nonnegative} and finiteness of the cover, at least
one of these balls has $\operatorname*{ess\,inf}_{B_{\rho/2}(x_j)}
(v_\infty-c_0)=0$. The weak Harnack inequality
\cite[Theorem~3.51; Theorem~3.59 and the ensuing discussion for $p\ge n$]{hkm},
applied to $v_\infty-c_0$ in $B_{4\rho}(x_j)$, gives, for some $q>0$,
\[
 \left(\frac1{|B_\rho(x_j)|}
 \int_{B_\rho(x_j)}(v_\infty-c_0)^q\,dx\right)^{1/q}
 \le C\operatorname*{ess\,inf}_{B_{\rho/2}(x_j)}(v_\infty-c_0)=0.
\]
Thus $v_\infty=c_0$ a.e. in $B_\rho(x_j)$. Since $n\ge2$, $\mathcal A$
is connected. The canonical lower semicontinuous representative of
$v_\infty-c_0$ is superharmonic for the limiting operator by
\cite[Theorems~3.63 and~7.16]{hkm}.
The strong minimum principle \cite[Theorem~7.12]{hkm} therefore yields
\begin{equation}\label{sec4:annular-constant}
 v_\infty=c_0
 \qquad\hbox{a.e. in }B_{3/2}\setminus\overline{B_1}.
\end{equation}

Set $w=v_\infty-c_0$.  Since $w=0$ a.e. in the exterior annulus and
$w\in W^{1,p}(B_{3/2})$, the exterior and interior Sobolev traces on
$\partial B_1$ agree and are both zero.  Therefore
\[
 w|_{B_1}\in W^{1,p}_0(B_1).
\]
Using $w^-:=\max\{-w,0\}$ as a test function in the supersolution
inequality in $B_1$
\cite[Lemma~3.11]{hkm} (or comparing with the zero solution) yields
\[
 w\ge0\qquad\hbox{a.e. in }B_1.
\]
Together with \eqref{sec4:annular-constant}, this shows that $w$ is a nonnegative weak
supersolution in the connected ball $B_{3/2}$ and vanishes on the nonempty
open set $B_{3/2}\setminus\overline{B_1}$.  The strong minimum principle
therefore gives
\[
 v_\infty=c_0\qquad\hbox{a.e. in }B_{3/2}.
\]
Consequently,
\[
 \int_{B_1}(A_\infty\nabla v_\infty\cdot\nabla v_\infty)^{p/2}\,dx=0,
\]
contradicting \eqref{sec4:positive-limit-energy}.  This contradiction proves
the lemma.
\end{proof}

\section{Interior estimates}
\label{sec:final-theorem}

In this section we combine the estimates of the preceding sections
with interpolation to obtain interior energy and higher-integrability
bounds, as well as H\"older estimates under the dimension restrictions
stated below.

\noindent We use the notation, norms, weak derivatives and zero-extension conventions of
Sections~2--3.

\begin{theorem}[]
Let $p>1$ and let $u$ be a stable regular solution of
\[
 -\operatorname{div}\bigl((A\nabla u\cdot\nabla u)^{p/2-1}A\nabla u\bigr)
 =f(u)\quad\hbox{in }B_1,
\]
where $f\in C^1(\mathbb R)$, $f\ge0$, and
$A\in W^{1,\infty}_{\rm loc}(B_1;\mathbb R^{n\times n})$ is real symmetric with
$\lambda I\le A\le\Lambda I$. Set $L=\|DA\|_{L^\infty(B_{3/4})}$.
There exist $\gamma_0=\gamma_0(n,p,\lambda,\Lambda)>0$ and
$C=C(n,p,\lambda,\Lambda,L)>0$ such that, in every dimension,
\[
 \left(\int_{B_{1/2}}(A\nabla u\cdot\nabla u)^{p/2}\,dx\right)^{1/p}
 +\|\nabla u\|_{L^{p+\gamma_0}(B_{1/2})}
 \le C\|u\|_{L^1(B_1)}.
\]
If, in addition,
\[
 n<p+\frac{4p}{p-1}\quad\hbox{if }p\ge2,
 \qquad n<5p\quad\hbox{if }1<p<2,
\]
then there exists $\alpha=\alpha(n,p,\lambda,\Lambda)\in(0,1)$ such that,
with $C$ enlarged if necessary and retaining the same dependence,
\[
 \|u\|_{C^\alpha(B_{1/16})}\le C\|u\|_{L^1(B_1)},
\]
and
\[
 \int_{B_r(y)}|\nabla u|^p\,dx
 \le Cr^{n-p+p\alpha}\|u\|_{L^1(B_1)}^p,
 \qquad y\in B_{1/8},\quad 0<r\le\frac1{16}.
\]
All exponents and constants are independent of $f$, $u$, and derivatives
of $A$ of order two or higher.
\end{theorem}

For the weighted argument, assume $1<p<n$, the solution and
coefficient hypotheses of Section~2 with $f\ge0$, and
\begin{equation}
 A(0)=I,\qquad \|DA\|_{L^\infty(B_1)}\le\varepsilon,
 \qquad 0\le\varepsilon\le\varepsilon_0\le\min\{1,\delta_0\}.
 \tag{5.1}\label{s5:normalized-setting}
\end{equation}
Here $\delta_0$ is the Section~\ref{sec:hessian-estimate} threshold for this
ellipticity class and $\varepsilon_0$ is structural; the bounds following
\eqref{normass} apply.

The weighted argument keeps the origin and coordinates fixed, including on
auxiliary covering balls. We normalize at other centers only after obtaining
complete-gradient energy decay. Constants in this normalized setting depend
only on $n,p,\lambda,\Lambda$ and the indicated structural parameters.

Let $0<\rho\le1/8$, $B_{4\rho}\Subset B_1$, $0<a\le n-p$, and set
\[
 c_p(a)=
 \begin{cases}
 pa-(p-1)a^2/4,&p\ge2,\\
 pa-a^2/4,&1<p<2.
 \end{cases}
\]
The nonnegativity of $f$ and coefficient smallness enter through
Corollary~\ref{cor:local-pplus1-hessian}.

By \eqref{local-pplus1-hessian} and ellipticity, for
$B_{2s}(z)\Subset B_1$ and $0<s\le1$,
\[
 \int_{B_s(z)}(A\nabla u\cdot\nabla u)^{p/2-1}|\nabla u|\,|D^2u|\,dx
 \le\frac Cs\int_{B_{2s}(z)}(A\nabla u\cdot\nabla u)^{p/2}\,dx.
\]
Consequently,
\begin{equation}
\begin{aligned}
 &\int_{B_{2\rho}}r^{2-a}(A\nabla u\cdot\nabla u)^{p/2-1}
       |\nabla u|\,|D^2u|\,dx\\
 &\qquad\le C\int_{B_{4\rho}}r^{1-a}
       (A\nabla u\cdot\nabla u)^{p/2}\,dx.
\end{aligned}
\tag{5.2}\label{s5:weighted-hessian-margin}
\end{equation}
Indeed, with $t_j=2^{-j}(2\rho)$, cover each annulus
$\{t_j/2\le r\le t_j\}$ by at most $M(n)$ balls of radius
$s_j=t_j/16$ centered in that annulus. Their doubles lie in
$\{3t_j/8<r<9t_j/8\}\Subset B_{4\rho}$, where $r\asymp t_j$.
Multiply the local estimate by $t_j^{2-a}$ and sum, using
$t_j^{2-a}/s_j=16t_j^{1-a}$ and bounded overlap of the enlarged annuli.
Monotone convergence proves \eqref{s5:weighted-hessian-margin}; its right
side is finite since $a<n$ and $u\in C^1_{\rm loc}$. This local boundedness
does not enter the constant.

Choose $\zeta\in C_c^\infty(B_{2\rho})$, $0\le\zeta\le1$, with
$\zeta=1$ on $B_\rho$ and $|\nabla\zeta|\le C/\rho$, and use
$\eta=r^{-a/2}\zeta$ in the proof of Theorem~\ref{raw}, with the radial
regularization of Remark~\ref{radial-weight-regularization}.
The Hessian error \eqref{BHess}, from $x_ku_{ki}\eta^2$ in
\eqref{derivativeceta}, retains $\zeta^2$ and is bounded by $C\varepsilon$
times the left side of \eqref{s5:weighted-hessian-margin}. The remaining
coefficient errors are supported in $B_{2\rho}$, while all $\nabla\zeta$
terms lie in $B_{2\rho}\setminus B_\rho$.
Comparing \eqref{rhsbelow} with \eqref{upperge2} or \eqref{upperlt2},
using \eqref{s5:weighted-hessian-margin}, and enlarging only the nonnegative
right-hand integrals therefore gives, whenever $c_p(a)>0$,
\begin{equation}
\begin{aligned}
 &c_p(a)\int_{B_\rho}r^{-a}
       (A\nabla u\cdot\nabla u)^{p/2-1}u_r^2\,dx\\
 &\quad +(n-p-a)\int_{B_\rho}r^{-a}
       (A\nabla u\cdot\nabla u)^{p/2}\,dx\\
 &\le C_a\int_{B_{4\rho}\setminus B_\rho}r^{-a}
       (A\nabla u\cdot\nabla u)^{p/2}\,dx
      +C_a\varepsilon\int_{B_{4\rho}}r^{1-a}
       (A\nabla u\cdot\nabla u)^{p/2}\,dx.
\end{aligned}
\tag{5.3}\label{s5:reduced-weighted}
\end{equation}
As in Remark~\ref{radial-weight-regularization}, first remove the critical-set
cutoff at fixed radial truncation, then the radial truncation;
\eqref{s5:weighted-hessian-margin} supplies the integrable Hessian majorant.
Here $C_a$ depends only on $n,p,\lambda,\Lambda,a$ and the fixed cutoff profile.
Step~2 uses the subcritical and critical choices of $a$.

\Needspace{17\baselineskip}
\begin{lemma}[]\label{lem:cube-interpolation}
Let $p\ge2$, $Q=(0,1)^n$, and $\sigma=|\nabla w|^{p-2}$. Assume
\[
 w\in C^1(Q)\cap W^{1,p-1}(Q),\qquad
 \nabla w\in W^{1,1}_{\sigma}(Q),
\]
where the weighted space is the completion for the $L^1$ norm and the seminorm
$\int_Q\sigma|D(\cdot)|$. Derivatives are identified on $\{\nabla w\ne0\}$,
with their weighted densities extended by zero on the critical set when $p>2$.
Then, for $0<\tau<1$,
\[
 \int_Q|\nabla w|^{p-1}\,dx
 \le C\tau\int_Q|\nabla w|^{p-2}|D^2w|\,dx
      +C\tau^{1-p}\int_Q|w|^{p-1}\,dx,
\]
where $C=C(n,p)$. For a translated cube $Q_\ell$ of side length $\ell>0$,
\[
 \int_{Q_\ell}|\nabla w|^{p-1}\,dx
 \le C\tau\ell\int_{Q_\ell}|\nabla w|^{p-2}|D^2w|\,dx
      +C(\tau\ell)^{1-p}\int_{Q_\ell}|w|^{p-1}\,dx.
\]
\end{lemma}

\begin{proof}
Put $g=\nabla w$ and $z_i=|g_i|^{p-2}g_i$. On compact subsets of
$\{g\ne0\}$ the continuous weight is bounded below, so weighted
approximation gives $g\in W^{1,1}_{\rm loc}(\{g\ne0\})$.
For $p>2$, choose a smooth $0\le\chi\le1$ equal to zero on $[0,1]$
and one on $[2,\infty)$, and put $g_\nu=g\chi(|g|/\nu)$.
Gluing with zero gives $g_\nu\in W^{1,1}(Q)$ and
$|Dg_\nu|\le C|Dg|\mathbf1_{\{|g|>\nu\}}$.
Truncating large values justifies the chain rule for
$z_{i,\nu}=|(g_\nu)_i|^{p-2}(g_\nu)_i$. Its derivative is bounded by
$C|g|^{p-2}|Dg|$, and the extra cutoff terms by
$C|g|^{p-2}|Dg|\mathbf1_{\{\nu<|g|<2\nu\}}$.
These bounds are integrable. Dominated convergence in the function and
weak derivative identities yields
\[
 z_i\in W^{1,1}(Q),\qquad
 Dz_i=(p-1)|g_i|^{p-2}Dg_i\quad\text{on }\{g\ne0\},
 \qquad |Dz_i|\le C|g|^{p-2}|Dg|,
\]
with zero extension. At $p=2$ this follows directly from $g\in W^{1,1}(Q)$.
Thus $|z_i|=|w_i|^{p-1}$ has absolutely continuous coordinate slices.

For an interval $I$ of length $h$ and a continuously differentiable $\phi$,
\[
 \inf_I|\phi'|^{p-1}\le Ch^{-p}\int_I|\phi|^{p-1}.
\]
If the infimum is nonzero, $\phi'$ has fixed sign; integrate between the
first and last quarters of $I$ and use H\"older to obtain this bound.
The fundamental theorem for the absolutely continuous function
$|\phi'|^{p-1}$ then gives
\[
 \int_I|\phi'|^{p-1}
 \le Ch\int_I|\phi'|^{p-2}|\phi''|
      +Ch^{1-p}\int_I|\phi|^{p-1},
\]
where the derivative term is understood through the weak chain rule above.
Partition $(0,1)$ into $k$ equal intervals with $\tau^{-1}\le k<2\tau^{-1}$.
Apply this estimate on almost every coordinate section, integrate over
the other variables and sum in $i$, using
$|w_i|^{p-2}|w_{ii}|\le|\nabla w|^{p-2}|D^2w|$.
This proves the unit-cube estimate; $x=x_0+\ell z$ gives the factors
$\tau\ell$ and $(\tau\ell)^{1-p}$ on $Q_\ell$.
\end{proof}

\Needspace{10\baselineskip}
\begin{proof}[Proof of Theorem~5.1]
Step~1 gives the energy and higher-integrability estimates in every dimension
and completes the case $n\le p$. Steps~2--4 complete the proof when $n>p$.

\subsubsection*{Step 1.}

We prove
\begin{equation}
 \left(\int_{B_{1/2}}(A\nabla u\cdot\nabla u)^{p/2}\,dx\right)^{1/p}
 \le C\|u\|_{L^1(B_1)}
 \tag{5.4}
\end{equation}
and
\begin{equation}
 \|\nabla u\|_{L^{p+\gamma_0}(B_{1/2})}
 \le C\|u\|_{L^1(B_1)}.
 \tag{5.5}
\end{equation}
We first prove (5.15) for small $DA$, then handle general
$L=\|DA\|_{L^\infty(B_{3/4})}$ by scalar rescaling and a finite covering,
without affine normalization.

\medskip
\noindent
Let $v$ be a stable regular solution with $f\in C^1(\mathbb R)$,
$f\ge0$, the same ellipticity bounds, and $\|DA\|_{L^\infty(B_1)}\le\delta_0$,
where $\delta_0$ is chosen as in Theorem~\ref{thm:hessian-estimate}.
The uniform $W^{1,\infty}$ bound on $A$ and
Lemma~\ref{lem:direct-anisotropic-higher}, together with ellipticity, give
a common $\gamma_0>0$ and
\begin{equation}
 \|\nabla v\|_{L^{p+\gamma_0}(B_{1/2})}
 \le C\|\nabla v\|_{L^p(B_1)}.
 \tag{5.6}
\end{equation}
Here $\gamma_0$ and $C$ depend only on $n,p,\lambda,\Lambda$; constants
with a small-parameter subscript may also depend on that parameter.

\paragraph{Case I: $p\ge2$.}
Covering $B_{3/4}$ by fixed interior balls and applying the scaled
Hessian estimate gives
\begin{equation}
 \int_{B_{3/4}}|\nabla v|^{p-2}|D^2v|\,dx
 \le C\|\nabla v\|_{L^p(B_1)}^{p-1}.
 \tag{5.7}
\end{equation}
For $p>2$, the weighted density is extended by zero on $\{\nabla v=0\}$,
as in Corollary~\ref{cor:dq-weighted-sobolev}. For $p=2$, the Appendix
gives $v\in W^{2,2}_{\rm loc}$, and Sobolev--Stampacchia gives
$D^2v=0$ a.e. on this set, so the same full-ball integral is valid.

As in \cite[proof of Theorem~1.1, Step~1, Case~1]{cms}, choose finitely many
congruent cubes $Q_j$ of fixed side length
$\ell_0=\ell_0(n)>0$ covering $B_{1/2}$, with
$\overline{Q_j}\subset B_{3/4}$. Corollary~\ref{cor:dq-weighted-sobolev}
allows us to apply Lemma~\ref{lem:cube-interpolation} on each cube.
Summing and absorbing the fixed side-length factors gives
\begin{equation}
 \int_{B_{1/2}}|\nabla v|^{p-1}\,dx
 \le C\tau\int_{B_{3/4}}|\nabla v|^{p-2}|D^2v|\,dx
      +C\tau^{1-p}\int_{B_{3/4}}|v|^{p-1}\,dx,
 \qquad 0<\tau<1.
 \tag{5.8}
\end{equation}
Interpolation between $L^{p+\gamma_0}$ and $L^{p-1}$, followed by (5.6)
and Young's inequality, yields, for $0<\eta<1$,
\begin{equation}
 \|\nabla v\|_{L^p(B_{1/2})}^{p-1}
 \le \eta\|\nabla v\|_{L^p(B_1)}^{p-1}
      +C_\eta\|\nabla v\|_{L^{p-1}(B_{1/2})}^{p-1}.
 \tag{5.9}
\end{equation}
Use this with $\eta/2$, insert (5.7)--(5.8), and then choose $\tau$
small enough that the additional gradient coefficient is at most $\eta/2$.
Thus
\begin{equation}
 \|\nabla v\|_{L^p(B_{1/2})}^{p-1}
 \le \eta\|\nabla v\|_{L^p(B_1)}^{p-1}
      +C_\eta\|v\|_{L^{p-1}(B_1)}^{p-1}.
 \tag{5.10}
\end{equation}

Put $\bar v=|B_1|^{-1}\int_{B_1}v\,dx$ and $w=v-\bar v$.
The shifted nonlinearity $\widehat f(t)=f(t+\bar v)$ is nonnegative,
with $\widehat f(w)=f(v)$ and $\widehat f'(w)=f'(v)$.
Since the gradient and stability test space are unchanged, $w$ is a stable
regular energy solution in the same class. Applying (5.10) to $w$ and
raising both sides to the power $1/(p-1)$ gives
\[
 \|\nabla v\|_{L^p(B_{1/2})}
 \le \eta^{1/(p-1)}\|\nabla v\|_{L^p(B_1)}
      +C_\eta\|w\|_{L^{p-1}(B_1)}.
\]
For $p>2$, interpolation, Poincar\'e's inequality, and
$\|w\|_{L^1(B_1)}\le2\|v\|_{L^1(B_1)}$ give, for every $t>0$,
\[
 \|w\|_{L^{p-1}(B_1)}
 \le \|w\|_{L^1(B_1)}^{1/(p-1)^2}
      \|w\|_{L^p(B_1)}^{1-1/(p-1)^2}
 \le t\|\nabla v\|_{L^p(B_1)}+C_t\|v\|_{L^1(B_1)}.
\]
Choosing $\eta$ first and then $t$ makes the coefficient of
$\|\nabla v\|_{L^p(B_1)}$ arbitrarily small. For $p=2$, the bound on $\|w\|_{L^1}$ gives this directly.

\paragraph{Case II: $1<p<2$.}
Interpolation between $L^{p+\gamma_0}$ and $L^1$, followed by (5.6), gives
\begin{equation}
 \|\nabla v\|_{L^p(B_{1/2})}
 \le \eta\|\nabla v\|_{L^p(B_1)}
      +C_\eta\|\nabla v\|_{L^1(B_{1/2})},
 \qquad 0<\eta<1.
 \tag{5.11}
\end{equation}
Theorem~\ref{dq:main-theorem} gives
$v\in W^{2,p}_{\rm loc}\subset W^{2,1}_{\rm loc}$.
On the fixed ball $U=B_{1/2}$, the interpolation inequality
\[
 \|Dw\|_{L^1(U)}
 \le \tau\|D^2w\|_{L^1(U)}+C\tau^{-1}\|w\|_{L^1(U)},
 \qquad w\in W^{2,1}(U),\quad 0<\tau<\tau_0,
\]
requires no boundary condition; see \cite[Theorem~7.28]{gilbarg-trudinger}.
Indeed, the extension construction in \cite[Theorem~7.25]{gilbarg-trudinger} gives
$F\in W^{2,1}(\mathbb R^n)$ satisfying $\|F\|_1\le C\|w\|_1$ and
\[
\|D^2F\|_1\le C(\|D^2w\|_1+\|Dw\|_1+\|w\|_1).
\]
Mollification gives $\|DF\|_1\le Ch\|D^2F\|_1+Ch^{-1}\|F\|_1$;
choosing $h$ proportional to $\tau$ and absorbing $Ch\|Dw\|_1$ proves
the inequality, with $C,\tau_0$ depending only on $U,n$.
Applied to $v$, it yields
\begin{equation}
 \|\nabla v\|_{L^1(B_{1/2})}
 \le \tau\|D^2v\|_{L^1(B_{1/2})}
      +C\tau^{-1}\|v\|_{L^1(B_1)}.
 \tag{5.12}
\end{equation}
Sobolev--Stampacchia gives $D^2v=0$ a.e. on $\{\nabla v=0\}$, so
Theorem~\ref{thm:hessian-estimate} and ellipticity imply
\[
 \|D^2v\|_{L^1(B_{1/2})}\le C\|\nabla v\|_{L^p(B_1)}.
\]
Insert this and (5.12) into (5.11), first choosing $\eta$ small and then
$\tau<\tau_0$ small enough to make the gradient coefficient arbitrarily small.

Both cases therefore give, for every $0<\eta<1$,
\begin{equation}
 \|\nabla v\|_{L^p(B_{1/2})}
 \le \eta\|\nabla v\|_{L^p(B_1)}
      +C_\eta\|v\|_{L^1(B_1)}.
 \tag{5.13}
\end{equation}
For $B_\rho(y)\subseteq B_1$, the rescaled solution
$v_{\rho,y}(z)=v(y+\rho z)$ has coefficients $A_{\rho,y}(z)=A(y+\rho z)$
and nonlinearity $\rho^pf$. Scalar dilation preserves the regular energy
class and carries the stability test space to the corresponding rescaled
space; moreover,
$\|DA_{\rho,y}\|_\infty\le\rho\delta_0\le\delta_0$.
Apply (5.13) to $v_{\rho,y}$, take the $p$-th power, and relabel the
arbitrarily small coefficient as $\eta$. Scaling back gives
\begin{equation}
 \rho^{p+n(p-1)}\int_{B_{\rho/2}(y)}|\nabla v|^p\,dx
 \le \eta\rho^{p+n(p-1)}\int_{B_\rho(y)}|\nabla v|^p\,dx
      +C_\eta\left(\int_{B_\rho(y)}|v|\,dx\right)^p.
 \tag{5.14}
\end{equation}
For open balls $B\subset\mathbb R^n$, set
$S(B)=\int_{B\cap B_1}|\nabla v|^p\,dx$. The scaling exponent in (5.14) is $p+n(p-1)>0$.
Then $S$ is nonnegative and finite-cover subadditive, and
$S(B_1)<\infty$ since $v\in W^{1,p}(B_1)$.
By (5.14) and $S(B_{\rho/4}(y))\le S(B_{\rho/2}(y))$,
the hypothesis of Simon's absorption lemma \cite{simon-scaling},
in the form of \cite[Lemma~A.4]{cms}, holds with
$C_0=C_\eta\|v\|_{L^1(B_1)}^p$ once
$\eta=\eta(n,p)>0$ is sufficiently small.
(The case $\|v\|_{L^1(B_1)}=0$ is immediate.)
The lemma therefore gives
\begin{equation}
 \|\nabla v\|_{L^p(B_{1/2})}\le C\|v\|_{L^1(B_1)}.
 \tag{5.15}\label{s5:small-coefficient-L1}
\end{equation}

\medskip
\noindent
Choose
\begin{equation}\label{step1-local-scale}
 R_L=c_0\min\left\{1,\frac{\delta_0}{1+L}\right\},
 \qquad 0<R_L\le\frac18,
 \tag{5.16}
\end{equation}
where $c_0=c_0(n,p,\lambda,\Lambda)$ is fixed small enough for the chosen
norm on $DA$. For $y\in B_{1/2}$, we have
$B_{R_L}(y)\subset B_{5/8}\Subset B_{3/4}$, and
\[
 v_y(z)=u(y+R_Lz),\qquad A_y(z)=A(y+R_Lz),\qquad f_y(t)=R_L^pf(t)
\]
defines a stable regular solution with $f_y\ge0$ and
$\|DA_y\|_\infty\le R_LL\le\delta_0$.
Estimate (5.15) therefore gives
\begin{equation}\label{step1-local-Lp}
 \|\nabla u\|_{L^p(B_{R_L/2}(y))}
 \le C R_L^{-1-n+n/p}\|u\|_{L^1(B_{R_L}(y))}.
 \tag{5.17}
\end{equation}
The scaled form of (5.6) on $B_{1/4}\subset B_{1/2}$, followed by (5.15),
also gives
\[
 \|\nabla v_y\|_{L^{p+\gamma_0}(B_{1/4})}
 \le C\|\nabla v_y\|_{L^p(B_{1/2})}
 \le C\|v_y\|_{L^1(B_1)}.
\]
Hence
\[
 \|\nabla u\|_{L^{p+\gamma_0}(B_{R_L/4}(y))}
 \le C R_L^{-1-n+n/(p+\gamma_0)}\|u\|_{L^1(B_{R_L}(y))}.
\]
Cover $B_{1/2}$ by at most $C(n)R_L^{-n}$ balls $B_{R_L/4}(y_j)$,
with centers in $B_{1/2}$ and bounded overlap of $B_{R_L}(y_j)$.
Sum the $p$-th powers in \eqref{step1-local-Lp} and the $(p+\gamma_0)$-th
powers in the last estimate. In each case, the sum of powers of the local
$L^1$ integrals is bounded by the power of their sum, and bounded overlap
controls that sum by $C\|u\|_{L^1(B_1)}$. Ellipticity then gives (5.4),
and the higher-exponent estimate gives (5.5). The constants are independent
of $f$ and $u$; only $C$ depends on $R_L$ and $L$, while $\gamma_0$ remains structural.

\medskip\noindent\textbf{The case $n\le p$.}

If $n\le p$, then $p+\gamma_0>n$; set
\[
 \alpha=1-\frac{n}{p+\gamma_0}\in(0,1).
\]
Poincar\'e's inequality followed by Morrey's embedding
\cite[Theorem~7.26(ii)]{gilbarg-trudinger}, applied to $u-u_{B_{1/2}}$, together with
$|u_{B_{1/2}}|\le |B_{1/2}|^{-1}\|u\|_{L^1(B_1)}$ and (5.5), gives
\[
 \|u\|_{C^\alpha(B_{1/16})}
 \le C\bigl(\|\nabla u\|_{L^{p+\gamma_0}(B_{1/2})}
             +|u_{B_{1/2}}|\bigr)
 \le C\|u\|_{L^1(B_1)}.
\]
For $y\in B_{1/8}$ and $0<r\le1/16$, we have $B_r(y)\subset B_{1/2}$.
H\"older's inequality and (5.5) yield
\[
\begin{aligned}
 \int_{B_r(y)}|\nabla u|^p\,dx
 &\le |B_r|^{\gamma_0/(p+\gamma_0)}
       \|\nabla u\|_{L^{p+\gamma_0}(B_{1/2})}^{p}\\
 &\le Cr^{n-p+p\alpha}\|u\|_{L^1(B_1)}^p,
\end{aligned}
\]
since $n\gamma_0/(p+\gamma_0)=n-p+p\alpha$.
Here $\alpha$ depends only on $n,p,\lambda,\Lambda$, and
$C=C(n,p,\lambda,\Lambda,L)$, as required. Together with (5.4)--(5.5),
this proves all assertions when $n\le p$.

Henceforth assume $n>p$.

\subsubsection*{Step 2. }

Steps~2--3 use the fixed normalized setting (5.1). The dimension assumption
gives $n-p<4p/(p-1)$ for $p\ge2$ and
$n-p<4p$ for $1<p<2$, hence $c_p(a)>0$ for every $0<a\le n-p$.

\medskip\noindent
Fix
\begin{equation}
 \theta=\frac12\min\{1,n-p\},\qquad a=n-p-\theta.
 \tag{5.18}
\end{equation}
Thus $0<\theta<\min\{1,n-p\}$, $0<a<n-p$, and $n-p-a=\theta$.
In \eqref{s5:reduced-weighted}, the coefficient error is bounded by
\[
 \widetilde C\varepsilon\rho
 \int_{B_{4\rho}}r^{p+\theta-n}(A\nabla u\cdot\nabla u)^{p/2}\,dx,
\]
since $r\le4\rho$. Choose a structural scale $\sigma_0>0$ such that
\begin{equation}
 \sigma_0\le\frac1{16},\qquad
 \widetilde C\varepsilon_0\sigma_0\le\frac\theta2.
 \tag{5.19}
\end{equation}
For $0<\rho\le\sigma_0$, split this error into the inner ball and its
complement in $B_{4\rho}$. Absorb the inner part into the complete energy,
leaving coefficient at least $\theta/2$, and drop the nonnegative radial term.
This gives
\begin{equation}
 \int_{B_\rho}r^{p+\theta-n}(A\nabla u\cdot\nabla u)^{p/2}\,dx
 \le C_\theta\int_{B_{4\rho}\setminus B_\rho}
       r^{p+\theta-n}(A\nabla u\cdot\nabla u)^{p/2}\,dx.
 \tag{5.20}\label{s5:subcritical}
\end{equation}
The reduction applies because $B_{4\rho}\Subset B_1$. The constants and
$\sigma_0$ are independent of $u,f$ and the actual $\varepsilon$;
no uniformity as the fixed $\theta$ tends to zero is needed.

\medskip\noindent

Define
\[
 D(\rho)=\rho^{p-n}\int_{B_\rho}(A\nabla u\cdot\nabla u)^{p/2}\,dx,
 \qquad
 B(\rho)=\int_{B_\rho}r^{p-n}
       (A\nabla u\cdot\nabla u)^{p/2-1}u_r^2\,dx.
\]
These quantities are finite by ellipticity, local $C^1$ regularity and
the local integrability of $r^{p-n}$, with the usual zero extension of the
radial density at $\nabla u=0$ when $1<p<2$.
Take $a=n-p$ in \eqref{s5:reduced-weighted}. Now $c_p(n-p)>0$ but
$n-p-a=0$, so the resulting estimate controls only the radial energy:
\begin{equation}
\begin{aligned}
 B(\rho)
 &\le C\int_{B_{4\rho}\setminus B_\rho}
       r^{p-n}(A\nabla u\cdot\nabla u)^{p/2}\,dx\\
 &\quad+C\varepsilon\int_{B_{4\rho}}
       r^{p+1-n}(A\nabla u\cdot\nabla u)^{p/2}\,dx,
 \qquad 0<\rho\le\sigma_0.
\end{aligned}
\tag{5.21}\label{s5:critical-with-error}
\end{equation}
The inner error is transferred to the annulus using $1-\theta>0$
and \eqref{s5:subcritical}:
\[
\begin{aligned}
 \int_{B_\rho}r^{p+1-n}(A\nabla u\cdot\nabla u)^{p/2}\,dx
 &\le\rho^{1-\theta}
       \int_{B_\rho}r^{p+\theta-n}(A\nabla u\cdot\nabla u)^{p/2}\,dx\\
 &\le C_\theta\rho
       \int_{B_{4\rho}\setminus B_\rho}
       r^{p-n}(A\nabla u\cdot\nabla u)^{p/2}\,dx.
\end{aligned}
\]
On the same annulus, $r\le4\rho$ bounds the remaining error.
Since $\varepsilon\rho\le\varepsilon_0\sigma_0$, we conclude that
\begin{equation}
 B(\rho)\le C\int_{B_{4\rho}\setminus B_\rho}
       r^{p-n}(A\nabla u\cdot\nabla u)^{p/2}\,dx,
 \qquad 0<\rho\le\sigma_0.
 \tag{5.22}\label{s5:critical-discrete-input}
\end{equation}

\subsubsection*{Step 3. }

Fix the structural radii and the upper-scale energy
\[
 R_\ast=\min\left\{\frac18,4\sigma_0\right\},\qquad
 \rho_\ast=\frac{R_\ast}{16},\qquad
 E_0=\int_{B_{1/2}}(A\nabla u\cdot\nabla u)^{p/2}\,dx,
\]
and, for $j\ge0$, set
\begin{equation}
 \rho_j=4^{-j}\rho_\ast,\qquad D_j=D(\rho_j),\qquad B_j=B(\rho_j).
 \tag{5.23}
\end{equation}
Since $\rho_j\le\rho_\ast\le\sigma_0/4$, \eqref{s5:critical-discrete-input}
applies at every discrete radius. For $j\ge1$, $4\rho_j=\rho_{j-1}$ and
$p-n<0$ give
\[
 B_j\le C\rho_j^{p-n}
       \int_{B_{\rho_{j-1}}}(A\nabla u\cdot\nabla u)^{p/2}\,dx
       =C4^{n-p}D_{j-1}.
\]
Thus, for a structural $K\ge1$,
\begin{equation}
 B_j\le K D_{j-1},\qquad j\ge1.
 \tag{5.24}\label{s5:discrete-radial-bound}
\end{equation}
The nonnegative radial density also gives $B_j\le B_{j-1}$.

Fix $\delta\in(0,1)$, for example $\delta=1/2$.
If $D_j\ge\delta D_{j-1}$, then
\[
 \int_{B_{\rho_j}}(A\nabla u\cdot\nabla u)^{p/2}\,dx
 \ge\delta4^{p-n}
       \int_{B_{2\rho_j}}(A\nabla u\cdot\nabla u)^{p/2}\,dx.
\]
Apply Lemma~\ref{thm:anisotropic-lemma41} to $u(\rho_j\,\cdot)$ on $B_2$
with $\delta_L=\delta4^{p-n}\in(0,1)$. Its dimension condition holds since
$n>p>1$, and
\[
 \|D[A(\rho_j\,\cdot)]\|_{L^\infty(B_2)}
 \le\rho_j\varepsilon_0\le\delta_0.
\]
The scalar rescaling preserves the ellipticity bounds, regularity,
equation and stability, with nonnegative $C^1$ nonlinearity $\rho_j^pf$;
the stability tests pull back in the corresponding full test spaces.
Scaling the comparison back and inserting the weight on the annulus gives
\[
\begin{aligned}
 D_j
 &\le C_{\delta_L}\rho_j^{p-n}
       \int_{B_{3\rho_j/2}\setminus B_{\rho_j}}
       (A\nabla u\cdot\nabla u)^{p/2-1}u_r^2\,dx\\
 &\le H(B_{j-1}-B_j).
\end{aligned}
\]
Here we used $B_{3\rho_j/2}\subset B_{\rho_{j-1}}$ and absorbed
$(3/2)^{n-p}$ into the structural constant $H\ge1$.
Together with the alternative $D_j<\delta D_{j-1}$, this proves
\begin{equation}
 D_j\le\delta D_{j-1}+H(B_{j-1}-B_j),\qquad j\ge1.
 \tag{5.25}\label{s5:discrete-recurrence}
\end{equation}

With $H,K$ fixed, choose $q>1$ sufficiently close to $1$ that
\begin{equation}
 q\le2,\qquad q<4^{1-\theta},\qquad
 \delta q+HKq(q-1)\le\frac{1+\delta}{2}<1.
 \tag{5.26}
\end{equation}
This is possible since the last left-hand side tends to $\delta$ as
$q\downarrow1$. Combining \eqref{s5:discrete-radial-bound} and
\eqref{s5:discrete-recurrence} gives
\[
\begin{aligned}
 D_j+HqB_j
 &\le[\delta+HK(q-1)]D_{j-1}+HB_{j-1}\\
 &\le q^{-1}(D_{j-1}+HqB_{j-1}).
\end{aligned}
\]
Since $B_{4\rho_\ast}\subset B_{1/2}$, inclusion and
\eqref{s5:critical-discrete-input} at $\rho_\ast$ give the initial bounds
without using a solution-dependent gradient bound:
\[
 D_0\le\rho_\ast^{p-n}E_0,\qquad
 B_0\le C\rho_\ast^{p-n}E_0.
\]
With $\rho_\ast$ structurally fixed, iteration and nonnegativity give
\begin{equation}
 D_j\le q^{-j}(D_0+HqB_0)\le Cq^{-j}E_0,\qquad j\ge0.
 \tag{5.27}
\end{equation}

Set
\begin{equation}
 \beta=\frac{\log q}{\log4},\qquad 0<\beta<1-\theta<1<p.
 \tag{5.28}
\end{equation}
For $0<r\le\rho_\ast$, choose $j\ge0$ with $\rho_{j+1}<r\le\rho_j$.
Monotonicity of the unnormalized complete energy gives
\[
 D(r)\le\left(\frac{\rho_j}{r}\right)^{n-p}D_j
 \le4^{n-p}D_j.
\]
Using $\rho_j<4r$ and $q^{-j}=(\rho_j/\rho_\ast)^\beta$, we obtain
\begin{equation}
 D(r)\le C\left(\frac r{\rho_\ast}\right)^\beta E_0,
 \qquad 0<r\le\rho_\ast.
 \tag{5.29}
\end{equation}
By the definition of $D$, (5.4), and the structural choice of $\rho_\ast$,
\begin{equation}
 \int_{B_r}(A\nabla u\cdot\nabla u)^{p/2}\,dx
 \le Cr^{n-p+\beta}\left(\int_{B_1}|u|\,dx\right)^p,
 \qquad 0<r\le\rho_\ast.
 \tag{5.30}\label{s5:one-center-decay}
\end{equation}
The constant and exponent are uniform in the normalized class, independent
of $f,u$, the actual $\varepsilon$, and higher derivatives of $A$.
Only this complete-energy bound is carried to new affine coordinates.

\subsubsection*{Step 4. }

Return to the original coefficient field. Fix $\varepsilon_0,\rho_\ast,\beta$
for the common transformed ellipticity class
$\lambda'=\lambda/\Lambda$, $\Lambda'=\Lambda/\lambda$, and set
\[
 L=\|DA\|_{L^\infty(B_{3/4})},\qquad
 R_0=\min\left\{\frac14,\frac{\lambda\varepsilon_0}{1+L}\right\},
 \qquad s=\frac{R_0}{\sqrt\Lambda}.
\]
For $y\in B_{1/8}$, let $S_y=A(y)^{1/2}$ and define
\[
 v_y(z)=u(y+sS_yz),\qquad
 \widetilde A_y(z)=S_y^{-1}A(y+sS_yz)S_y^{-1}.
\]
The image of $B_1$ lies in $B_{R_0}(y)\subset B_{3/8}\Subset B_{3/4}$.
Since $S_y$ is constant in $z$, the weak chain rule uses only first
derivatives of $A$ and gives
\[
\begin{gathered}
 \widetilde A_y(0)=I,\qquad
 \lambda'I\le\widetilde A_y\le\Lambda'I,\\
 \|D\widetilde A_y\|_{L^\infty(B_1)}
 \le s\|S_y^{-1}\|_{\rm op}^2\|S_y\|_{\rm op}L
 \le\frac{R_0L}{\lambda}\le\varepsilon_0.
\end{gathered}
\]
Under $x=y+sS_yz$,
\[
 \nabla_zv_y=sS_y\nabla_xu,\qquad
 dx=s^n\det(S_y)\,dz,\qquad
 \widetilde A_y\nabla v_y\cdot\nabla v_y=s^2A\nabla u\cdot\nabla u,
\]
so the transformed equation is
\[
 -\operatorname{div}_z\bigl(
 (\widetilde A_y\nabla v_y\cdot\nabla v_y)^{p/2-1}
 \widetilde A_y\nabla v_y\bigr)
 =\widetilde f_y(v_y),\qquad \widetilde f_y(t)=s^pf(t).
\]
Thus $\widetilde f_y\in C^1$ and $\widetilde f_y\ge0$. Each transformed
stability term for $\psi(z)$ equals the corresponding term for its pullback
$\varphi(x)$ times $s^{p-n}/\det(S_y)$. For $p\ge2$, affine pullback
identifies the completions of the smooth test classes with equivalent
$L^2$ and weighted-gradient norms; for $1<p<2$, it preserves $W^{1,2}_0$,
finite weighted energy, and the critical-set gradient condition.
Zero extension from the interior ellipsoid gives admissible original tests,
and $f'(u)$ is bounded there. Hence $v_y$ is stable on the full prescribed
test space, and (5.30) applies with constants independent of $y$.

Writing $E_t(y)=y+sS_yB_t$, ellipticity gives
\[
 B_{s\sqrt\lambda\,t}(y)\subset E_t(y)
 \subset B_{s\sqrt\Lambda\,t}(y),\qquad
 \lambda^{n/2}\le\det(S_y)\le\Lambda^{n/2}.
\]
For $0<r\le r_\ast:=s\sqrt\lambda\,\rho_\ast$, put
$t=r/(s\sqrt\lambda)\le\rho_\ast$. By this inclusion, the Jacobian
identity, (5.30), and ellipticity,
\[
\begin{aligned}
 \int_{B_r(y)}|\nabla u|^p\,dx
 &\le Cs^{n-p}\det(S_y)
       \int_{B_t}(\widetilde A_y\nabla v_y\cdot\nabla v_y)^{p/2}\,dz\\
 &\le Cs^{n-p}\det(S_y)t^{n-p+\beta}
       \left(\frac1{s^n\det(S_y)}\int_{E_1(y)}|u|\,dx\right)^p\\
 &\le CR_0^{-np-\beta}r^{n-p+\beta}\|u\|_{L^1(B_1)}^p.
\end{aligned}
\]
For $r_\ast<r\le1/16$, use $B_r(y)\subset B_{3/16}\subset B_{1/2}$,
(5.4), and $r^{n-p+\beta}\ge r_\ast^{n-p+\beta}$ to enlarge the constant.
We obtain
\begin{equation}
 \int_{B_r(y)}|\nabla u|^p\,dx
 \le Cr^{n-p+\beta}\|u\|_{L^1(B_1)}^p,
 \qquad y\in B_{1/8},\quad 0<r\le\frac1{16}.
 \tag{5.31}\label{s5:all-center-morrey}
\end{equation}
Here $C=C(n,p,\lambda,\Lambda,L)$, whereas
\begin{equation}
 \alpha=\frac\beta p\in(0,1)
 \tag{5.32}
\end{equation}
depends only on $n,p,\lambda,\Lambda$, and $n-p+\beta=n-p+p\alpha$.

\medskip\noindent\textbf{The Campanato estimate and the full $C^\alpha$ norm.}

For $z\in B_{1/8}$ and $0<r\le1/16$, Poincar\'e's inequality and (5.31) give
\[
 \int_{B_r(z)}|u-u_{B_r(z)}|^p\,dx
 \le Cr^p\int_{B_r(z)}|\nabla u|^p\,dx
 \le Cr^{n+p\alpha}\|u\|_{L^1(B_1)}^p.
\]
These full balls lie in $B_{3/16}$. Telescoping concentric averages gives
\[
|u(z)-u_{B_r(z)}|\le Cr^\alpha\|u\|_{L^1(B_1)}.
\]
For points at distance $d\le1/64$, compare their radius-$d$ averages
in the radius-$2d$ ball about one of them. Subdividing longer segments
into $16$ equal parts within the convex ball $B_{1/8}$ then yields
\begin{equation}
 [u]_{C^\alpha(B_{1/8})}\le C\|u\|_{L^1(B_1)}.
 \tag{5.33}
\end{equation}
Finally, $|u_{B_{1/8}}|\le C\|u\|_{L^1(B_1)}$ and averaging the
H\"older differences gives
$|u(x)-u_{B_{1/8}}|\le C[u]_{C^\alpha(B_{1/8})}$ for $x\in B_{1/16}$.
Consequently,
\[
 \|u\|_{C^\alpha(B_{1/16})}\le C\|u\|_{L^1(B_1)}.
\]
Together with (5.4) and (5.5), this proves all assertions when $n>p$.
\end{proof}

\appendix
\titlecontents{section}[0pt]{\addvspace{2pt}\filright}
              {\contentspush{Appendix~\thecontentslabel\ }}
              {}{\titlerule*[8pt]{.}\contentspage}
\section{Difference-quotient regularity for all \(p>1\)}
\label{dq:difference-quotient-appendix}

This appendix establishes the local regularity results and weighted
Hessian estimates used in the preceding sections.

Let $\Omega\subset\mathbb R^n$ be open, $p>1$, and
$A\in W^{1,\infty}_{\rm loc}(\Omega;\mathbb R^{n\times n})$
be real symmetric. We use its locally Lipschitz representative and assume
that, for some $0<\lambda\le\Lambda<\infty$,
\begin{equation}
\lambda|\xi|^2\le A(x)\xi\cdot\xi\le\Lambda|\xi|^2
\qquad\text{for every }x\in\Omega,\ \xi\in\mathbb R^n.
\label{dq:ellipticity}
\end{equation}
For $1<p<2$, the flux
$(A(x)\xi\cdot\xi)^{\frac{p-2}{2}}A(x)\xi$
is extended continuously by zero at $\xi=0$.
Let $f\in C^1(\R)$ and let $u\in W^{1,p}_{\rm loc}(\Omega)$ satisfy
$f(u)\in L^\infty_{\rm loc}(\Omega)$ and
\begin{equation}
\int_\Omega
(A(x)\nabla u\cdot\nabla u)^{\frac{p-2}{2}}
A(x)\nabla u\cdot\nabla\varphi\,dx
=
\int_\Omega f(u)\varphi\,dx
\label{dq:weak-form}
\end{equation}
for every $\varphi\in C_c^\infty(\Omega)$.
Applying the local regularity results in Section~\ref{intro:solution-class}
on interior balls gives a representative
$u\in C^{1,\vartheta}_{\rm loc}(\Omega)$ for some $\vartheta>0$.
We use this representative below.

The proof adapts the difference-quotient arguments presented in
\cite[Theorems~4.1 and~4.6]{lindqvist} to the variable matrix and reaction term
above; the additional coefficient and reaction estimates are proved below.

\begin{theorem}[]
\label{dq:main-theorem}\label{dq:thm:main}
Let
\[
B_\rho\Subset B_{\rho_1}\Subset B_R\Subset B_S\Subset\Omega
\]
be concentric balls, and put
\begin{equation}
K:=\|DA\|_{L^\infty(B_S)},
\qquad
L_S:=\sup_{|t|\leq\|u\|_{L^\infty(B_S)}}|f'(t)|,
\qquad
M:=\|\nabla u\|_{L^\infty(B_S)}.
\label{dq:LS-def}
\end{equation}

\medskip\noindent\textbf{Case 1: $p>2$.}
The field
\[
V:=\bigl(A(x)\nabla u\cdot\nabla u\bigr)^{\frac{p-2}{4}}\nabla u
\]
belongs to $W^{1,2}(B_\rho;\R^n)$, and
\begin{align}
\int_{B_\rho}|DV|^2\,dx
&\leq
C\left(\frac{1}{(R-\rho)^2}+K^2\right)
\int_{B_S}\bigl(A(x)\nabla u\cdot\nabla u\bigr)^{p/2}\,dx
+CL_S\int_{B_S}|\nabla u|^2\,dx.
\label{dq:V-estimate}
\end{align}
Recovering $\nabla u$ from $V$ locally on the noncritical set gives
\begin{equation}
\int_{B_\rho\cap\{|\nabla u|>0\}}
\bigl(A(x)\nabla u\cdot\nabla u\bigr)^{\frac{p-2}{2}}|D^2u|^2\,dx
<\infty.
\label{dq:weighted-Hessian}
\end{equation}

\medskip\noindent\textbf{Case 2: $1<p<2$.}
One has $u\in W^{2,p}(B_{\rho_1})$ and
\begin{equation}
\int_{B_\rho\cap\{|\nabla u|>0\}}
\bigl(A(x)\nabla u\cdot\nabla u\bigr)^{\frac{p-2}{2}}|D^2u|^2\,dx
<\infty.
\label{dq:eq:mainestimate}
\end{equation}

\medskip\noindent\textbf{Case 3: $p=2$.}
One has $u\in W^{2,2}(B_\rho)$.
In~\eqref{dq:V-estimate}, $C=C(n,p,\lambda,\Lambda)$. No stability or sign assumption on $f$ is used.
\end{theorem}

\begin{lemma}[]
\label{dq:weak-extension}\label{dq:boundedness-lemma}
For every open $U\Subset\Omega$, \eqref{dq:weak-form} holds for
$\varphi\in W_0^{1,p}(U)$, extended by zero. The quantities in
\eqref{dq:LS-def}, as well as $\|u\|_{L^\infty(B_S)}$, are finite.
\end{lemma}

\begin{proof}
The finiteness of $K$ follows from $A\in W^{1,\infty}_{\rm loc}$;
the other bounds follow from $u\in C^1(\Omega)$ and $f\in C^1(\R)$
by compactness. Ellipticity bounds the flux by $C|\nabla u|^{p-1}$;
both the flux and $f(u)$ belong to $L^{p'}(U)$, where $p'=p/(p-1)$.
Density of $C_c^\infty(U)$ in $W_0^{1,p}(U)$ therefore extends the identity.
\end{proof}

Let $B\in\R^{n\times n}$ be real symmetric and satisfy
$\lambda|\zeta|^2\le B\zeta\cdot\zeta\le\Lambda|\zeta|^2$ for every
$\zeta\in\R^n$. At $p=2$, zero powers are understood as $1$.

\begin{lemma}
\label{dq:matrix-inequalities}
For $p\ge2$ and $\xi,\eta\in\R^n$,
\begin{align}
&\Big[(B\xi\cdot\xi)^{\frac{p-2}{2}}B\xi
 -(B\eta\cdot\eta)^{\frac{p-2}{2}}B\eta\Big]\cdot(\xi-\eta)
 \notag\\
&\qquad\ge c\left|
 (B\xi\cdot\xi)^{\frac{p-2}{4}}\xi
 -(B\eta\cdot\eta)^{\frac{p-2}{4}}\eta\right|^2,
\label{dq:matrix-monotonicity-V}\\
&\Big[(B\xi\cdot\xi)^{\frac{p-2}{2}}B\xi
 -(B\eta\cdot\eta)^{\frac{p-2}{2}}B\eta\Big]\cdot(\xi-\eta)
 \notag\\
&\qquad\ge c(|\xi|+|\eta|)^{p-2}|\xi-\eta|^2,
\label{dq:matrix-strong-monotonicity}\\
&\left|(B\xi\cdot\xi)^{\frac{p-2}{2}}B\xi
 -(B\eta\cdot\eta)^{\frac{p-2}{2}}B\eta\right|
 \notag\\
&\qquad\le C\left[(B\xi\cdot\xi)^{\frac{p-2}{4}}
 +(B\eta\cdot\eta)^{\frac{p-2}{4}}\right]
 \notag\\
&\hspace{30mm}\times\left|(B\xi\cdot\xi)^{\frac{p-2}{4}}\xi
 -(B\eta\cdot\eta)^{\frac{p-2}{4}}\eta\right|.
\label{dq:matrix-flux-difference}
\end{align}
Here $c,C>0$ depend only on $p,\lambda,\Lambda$.
\end{lemma}

\begin{proof}
Set $z=B^{1/2}\xi$ and $w=B^{1/2}\eta$. The standard Euclidean
estimates (see \cite[Section~12]{lindqvist}) give
\[
 c(|z|+|w|)^{p-2}\le\int_0^1|w+t(z-w)|^{p-2}\,dt
 \le(|z|+|w|)^{p-2}.
\]
The lower bound follows from the quarter of the segment adjacent to the
larger endpoint, where the norm is at least $(|z|+|w|)/4$.
The eigenvalues of $D(|y|^{p-2}y)$ lie between $|y|^{p-2}$ and
$(p-1)|y|^{p-2}$, while $\|D(|y|^{(p-2)/2}y)\|_{\op}
\le(p/2)|y|^{(p-2)/2}$. Segment integration and Cauchy--Schwarz give
\begin{align*}
 \bigl(|z|^{p-2}z-|w|^{p-2}w\bigr)\cdot(z-w)
 &\ge |z-w|^2\int_0^1|w+t(z-w)|^{p-2}\,dt\\
 &\ge c(|z|+|w|)^{p-2}|z-w|^2,\\
 \left||z|^{(p-2)/2}z-|w|^{(p-2)/2}w\right|^2
 &\le\frac{p^2}{4}|z-w|^2\int_0^1|w+t(z-w)|^{p-2}\,dt\\
 &\le\frac{p^2}{4}
 \bigl(|z|^{p-2}z-|w|^{p-2}w\bigr)\cdot(z-w).
\end{align*}
The same argument with exponent $(p-2)/2$ and the upper derivative bound give
\begin{align*}
 \left||z|^{(p-2)/2}z-|w|^{(p-2)/2}w\right|
 &\ge c(|z|+|w|)^{(p-2)/2}|z-w|,\\
 \left||z|^{p-2}z-|w|^{p-2}w\right|
 &\le(p-1)(|z|+|w|)^{p-2}|z-w|\\
 &\le C(|z|+|w|)^{(p-2)/2}
 \left||z|^{(p-2)/2}z-|w|^{(p-2)/2}w\right|.
\end{align*}
Since $(|z|+|w|)^{(p-2)/2}\le C_p(|z|^{(p-2)/2}+|w|^{(p-2)/2})$,
this proves the Euclidean flux-difference estimate.
The substitution preserves the monotonicity pairing; the Euclidean flux
and $V$ differences are $B^{-1/2}$ and $B^{1/2}$ times their matrix
counterparts. Ellipticity therefore gives all three stated estimates.
\end{proof}

\begin{lemma}
\label{dq:matrix-variation}\label{dq:lem:matrixvariation}
Let $p>1$ and let $B,\widehat B$ satisfy the same symmetry and
ellipticity bounds. For every $\xi\in\R^n$,
\begin{align}
&\left|(B\xi\cdot\xi)^{\frac{p-2}{2}}B\xi
 -(\widehat B\xi\cdot\xi)^{\frac{p-2}{2}}\widehat B\xi\right|
 \le C\|B-\widehat B\|_{\op}|\xi|^{p-1},
\label{dq:matrix-flux-variation}\\
&\left|(B\xi\cdot\xi)^{\frac{p-2}{4}}\xi
 -(\widehat B\xi\cdot\xi)^{\frac{p-2}{4}}\xi\right|
 \le C\|B-\widehat B\|_{\op}|\xi|^{p/2}.
\label{dq:matrix-V-variation}
\end{align}
Here $C=C(p,\lambda,\Lambda)$; both fields are extended continuously
by zero at $\xi=0$.
\end{lemma}

\begin{proof}
For $\xi\ne0$, the segment $B_t=\widehat B+t(B-\widehat B)$ retains
the ellipticity bounds. Differentiating the two fields gives
\begin{align*}
 \frac{d}{dt}\left[(B_t\xi\cdot\xi)^{\frac{p-2}{2}}B_t\xi\right]
 &=\frac{p-2}{2}(B_t\xi\cdot\xi)^{\frac{p-4}{2}}
   \bigl((B-\widehat B)\xi\cdot\xi\bigr)B_t\xi\\*
 &\quad +(B_t\xi\cdot\xi)^{\frac{p-2}{2}}(B-\widehat B)\xi,\\*
 \frac{d}{dt}\left[(B_t\xi\cdot\xi)^{\frac{p-2}{4}}\xi\right]
 &=\frac{p-2}{4}(B_t\xi\cdot\xi)^{\frac{p-6}{4}}
   \bigl((B-\widehat B)\xi\cdot\xi\bigr)\xi.
\end{align*}
Ellipticity bounds these derivatives respectively by
$C\|B-\widehat B\|_{\op}|\xi|^{p-1}$ and
$C\|B-\widehat B\|_{\op}|\xi|^{p/2}$ for every $p>1$.
Integration proves both estimates; continuity handles $\xi=0$.
\end{proof}

\begin{lemma}
\label{dq:lem:monotonicity}
For $1<p<2$ and $B$ as above,
\begin{align}
&\Big[(B\xi\cdot\xi)^{\frac{p-2}{2}}B\xi
 -(B\eta\cdot\eta)^{\frac{p-2}{2}}B\eta\Big]\cdot(\xi-\eta)
 \notag\\
&\qquad\ge c_p(B\xi\cdot\xi+B\eta\cdot\eta)^{\frac{p-2}{2}}
 |B^{1/2}(\xi-\eta)|^2,
\label{dq:eq:singmon}
\end{align}
where $c_p=(p-1)2^{(p-2)/2}$ and the right-hand side is zero when
$\xi=\eta=0$.
\end{lemma}

\begin{proof}
Set $z=B^{1/2}\xi$, $w=B^{1/2}\eta$ and, for $\varepsilon>0$,
$F_\varepsilon(y)=(\varepsilon^2+|y|^2)^{(p-2)/2}y$.
Since $p-2<0$ and $(y\cdot a)^2\le(\varepsilon^2+|y|^2)|a|^2$,
\[
\begin{aligned}
 DF_\varepsilon(y)a\cdot a
 &=(\varepsilon^2+|y|^2)^{\frac{p-2}{2}}|a|^2
   +(p-2)(\varepsilon^2+|y|^2)^{\frac{p-4}{2}}(y\cdot a)^2\\
 &\ge(p-1)(\varepsilon^2+|y|^2)^{\frac{p-2}{2}}|a|^2.
\end{aligned}
\]
For $z,w$ not both zero, integrate along $w+t(z-w)$. The bound
$|w+t(z-w)|\le|z|+|w|$ and the negative exponent give
\begin{align*}
 \bigl(F_\varepsilon(z)-F_\varepsilon(w)\bigr)\cdot(z-w)
 &\ge(p-1)|z-w|^2\int_0^1
   (\varepsilon^2+|w+t(z-w)|^2)^{\frac{p-2}{2}}\,dt\\
 &\ge(p-1)\bigl(\varepsilon^2+(|z|+|w|)^2\bigr)^{\frac{p-2}{2}}
   |z-w|^2.
\end{align*}
Letting $\varepsilon\downarrow0$ and using
$(|z|+|w|)^2\le2(|z|^2+|w|^2)$ gives \eqref{dq:eq:singmon}
with the stated $c_p$; the double-zero case follows by definition.
\end{proof}

For $1<p<2$ and $a,b\ge0$, Young's inequality with conjugate exponents
$p/(2-p)$ and $p/[2(p-1)]$ gives
\begin{equation}
 a^{2-p}b^2\le\frac{2-p}{p}a^p
 +\frac{2(p-1)}p b^{\frac p{p-1}}
 \le C_p\left(a^p+b^{\frac p{p-1}}\right).
\label{dq:eq:youngspecial}
\end{equation}

\begin{proof}[Proof of Theorem~\ref{dq:main-theorem}]
Fix a coordinate direction $e_k$ and write
\begin{equation}
\begin{gathered}
u_h(x)=u(x+he_k),\qquad A_h(x)=A(x+he_k),\qquad\\
\delta_hu=\frac{u_h-u}{h},\qquad
\delta_h\nabla u=\frac{\nabla u_h-\nabla u}{h}.
\end{gathered}
\label{dq:uh-dq}
\end{equation}
For $0<|h|<S-R$, the segment formula and Jensen's inequality give
\begin{equation}
\begin{gathered}
\delta_hu(x)=\int_0^1\partial_ku(x+the_k)\,dt,\\
\|\delta_hu\|_{L^p(B_R)}\le\|\nabla u\|_{L^p(B_S)},\qquad
\|\delta_hu\|_{L^2(B_R)}\le\|\nabla u\|_{L^2(B_S)}.
\end{gathered}
\label{dq:dq-Lp}
\end{equation}
In particular,
\begin{equation}
|\delta_hu|\le M,\qquad
\left|\frac{f(u_h)-f(u)}h\right|\le L_S|\delta_hu|.
\label{dq:dq-pointwise-bounded}%
\makeatletter\ltx@label{dq:f-local-Lipschitz}\makeatother
\end{equation}
Choose $\zeta\in C_c^\infty(B_R)$ with $0\le\zeta\le1$; the region
where $\zeta=1$ will be specified in each case. Translate the weak
equation by $he_k$ and test both equations with
$\zeta^2(u_h-u)\in W_0^{1,p}(B_R)$, using
Lemma~\ref{dq:weak-extension}. Subtracting and dividing by $h^2$ yields
\begin{align}
&\int_{B_R}\zeta^2
\frac{(A_h\nabla u_h\cdot\nabla u_h)^{\frac{p-2}{2}}A_h\nabla u_h
 -(A(x)\nabla u\cdot\nabla u)^{\frac{p-2}{2}}A(x)\nabla u}{h}
\cdot\delta_h\nabla u\,dx\notag\\
&=-2\int_{B_R}\zeta\,\delta_hu
\frac{(A_h\nabla u_h\cdot\nabla u_h)^{\frac{p-2}{2}}A_h\nabla u_h
 -(A(x)\nabla u\cdot\nabla u)^{\frac{p-2}{2}}A(x)\nabla u}{h}
\cdot\nabla\zeta\,dx\notag\\
&\quad+\int_{B_R}\zeta^2\frac{f(u_h)-f(u)}h\,\delta_hu\,dx.
\label{dq:eq:exactidentity}
\end{align}
The coefficient change, separated by adding and subtracting
$(A_h\nabla u\cdot\nabla u)^{(p-2)/2}A_h\nabla u$, satisfies
\begin{equation}
\left|\frac{(A_h\nabla u\cdot\nabla u)^{\frac{p-2}{2}}A_h\nabla u
 -(A(x)\nabla u\cdot\nabla u)^{\frac{p-2}{2}}A(x)\nabla u}{h}\right|
\le CK|\nabla u|^{p-1}
\label{dq:coefficient-vector-bound}%
\makeatletter\ltx@label{dq:eq:coefexplicit}\makeatother
\end{equation}
by \eqref{dq:matrix-flux-variation} and the Lipschitz bound for $A$.
All constants below are independent of $h$ and depend only on
$n,p,\lambda,\Lambda$, unless a dependence on $\varepsilon$ is indicated.

\medskip\noindent\textbf{Case 1: $p>2$.}
Set $\widehat\rho=(\rho+R)/2$, choose $\zeta=1$ on
$B_{\widehat\rho}$ with $|\nabla\zeta|\le C/(R-\rho)$, and take
$0<|h|<h_0:=\tfrac12\min\{S-R,\widehat\rho-\rho\}$.
In \eqref{dq:eq:exactidentity}, split each flux difference as
\begin{align*}
&\frac{(A_h\nabla u_h\cdot\nabla u_h)^{\frac{p-2}{2}}A_h\nabla u_h
 -(A(x)\nabla u\cdot\nabla u)^{\frac{p-2}{2}}A(x)\nabla u}{h}\\*
&\quad=\frac{(A_h\nabla u_h\cdot\nabla u_h)^{\frac{p-2}{2}}A_h\nabla u_h
 -(A_h\nabla u\cdot\nabla u)^{\frac{p-2}{2}}A_h\nabla u}{h}\\*
&\qquad+\frac{(A_h\nabla u\cdot\nabla u)^{\frac{p-2}{2}}A_h\nabla u
 -(A(x)\nabla u\cdot\nabla u)^{\frac{p-2}{2}}A(x)\nabla u}{h}.
\end{align*}
The two summands separate the gradient and coefficient changes.
By \eqref{dq:matrix-monotonicity-V} and
\eqref{dq:matrix-strong-monotonicity} with $B=A_h$, averaging the
two lower bounds gives
\begin{align*}
&\int_{B_R}\zeta^2
\frac{(A_h\nabla u_h\cdot\nabla u_h)^{\frac{p-2}{2}}A_h\nabla u_h
 -(A_h\nabla u\cdot\nabla u)^{\frac{p-2}{2}}A_h\nabla u}{h}
\cdot\delta_h\nabla u\,dx\\*
&\quad\ge c\int_{B_R}\zeta^2\left|
\frac{(A_h\nabla u_h\cdot\nabla u_h)^{\frac{p-2}{4}}\nabla u_h
 -(A_h\nabla u\cdot\nabla u)^{\frac{p-2}{4}}\nabla u}{h}\right|^2\,dx\\*
&\qquad+c\int_{B_R}\zeta^2
(|\nabla u_h|+|\nabla u|)^{p-2}|\delta_h\nabla u|^2\,dx.
\end{align*}
Here $c>0$ depends only on $p,\lambda,\Lambda$.
For the cutoff term arising from the first summand, use
\eqref{dq:matrix-flux-difference}, ellipticity, and H\"older with
exponents $p,2p/(p-2),2$. Its absolute value is bounded by
\begin{align*}
&C\|\nabla\zeta\|_\infty\|\delta_hu\|_{L^p(B_R)}
\left\|(|\nabla u_h|+|\nabla u|)^{\frac{p-2}{2}}
\right\|_{L^{2p/(p-2)}(B_R)}\\*
&\qquad\times\left\|\zeta
\frac{(A_h\nabla u_h\cdot\nabla u_h)^{\frac{p-2}{4}}\nabla u_h
 -(A_h\nabla u\cdot\nabla u)^{\frac{p-2}{4}}\nabla u}{h}
\right\|_{L^2(B_R)}.
\end{align*}
By \eqref{dq:dq-Lp}, translation, and ellipticity, the product of
the $L^p$ and $L^{2p/(p-2)}$ norms is at most
$C[\int_{B_S}(A(x)\nabla u\cdot\nabla u)^{p/2}\,dx]^{1/2}$.
Young's inequality therefore bounds this cutoff term by
\begin{align*}
&\varepsilon\int_{B_R}\zeta^2\left|
\frac{(A_h\nabla u_h\cdot\nabla u_h)^{\frac{p-2}{4}}\nabla u_h
 -(A_h\nabla u\cdot\nabla u)^{\frac{p-2}{4}}\nabla u}{h}\right|^2\,dx\\*
&\qquad+C_\varepsilon\|\nabla\zeta\|_\infty^2
\int_{B_S}(A(x)\nabla u\cdot\nabla u)^{p/2}\,dx.
\end{align*}
For the coefficient-gradient term, \eqref{dq:coefficient-vector-bound}
and weighted Young give
\begin{align*}
CK\int_{B_R}\zeta^2|\nabla u|^{p-1}|\delta_h\nabla u|\,dx
&\le\varepsilon\int_{B_R}\zeta^2
(|\nabla u_h|+|\nabla u|)^{p-2}|\delta_h\nabla u|^2\,dx\\
&\quad+C_\varepsilon K^2
\int_{B_S}(A(x)\nabla u\cdot\nabla u)^{p/2}\,dx,
\end{align*}
where the remainder follows from
\[
|\nabla u|^{2p-2}(|\nabla u_h|+|\nabla u|)^{2-p}
\le |\nabla u|^p,
\]
with the product defined as zero when both gradients vanish.
For the coefficient-cutoff term, H\"older, \eqref{dq:dq-Lp}, and
$2ab\le a^2+b^2$ yield
\begin{align*}
&CK\int_{B_R}\zeta|\delta_hu|\,|\nabla u|^{p-1}|\nabla\zeta|\,dx\\*
&\quad\le C\bigl(K^2+\|\nabla\zeta\|_\infty^2\bigr)
\int_{B_S}(A(x)\nabla u\cdot\nabla u)^{p/2}\,dx.
\end{align*}
Finally, \eqref{dq:f-local-Lipschitz} and the $L^2$ bound in
\eqref{dq:dq-Lp} give
\begin{align*}
\left|\int_{B_R}\zeta^2\frac{f(u_h)-f(u)}h\,\delta_hu\,dx\right|
&\le L_S\int_{B_R}\zeta^2|\delta_hu|^2\,dx\\
&\le L_S\int_{B_S}|\nabla u|^2\,dx.
\end{align*}
Absorbing the two errors with a sufficiently small structural
$\varepsilon$ and restricting to $B_{\widehat\rho}$ gives
\begin{align}
&\int_{B_{\widehat\rho}}\left|
\frac{(A_h\nabla u_h\cdot\nabla u_h)^{\frac{p-2}{4}}\nabla u_h
 -(A_h\nabla u\cdot\nabla u)^{\frac{p-2}{4}}\nabla u}{h}\right|^2\,dx
\notag\\*
&\quad+\int_{B_{\widehat\rho}}
(|\nabla u_h|+|\nabla u|)^{p-2}|\delta_h\nabla u|^2\,dx
\notag\\*
&\le C\left((R-\rho)^{-2}+K^2\right)
\int_{B_S}(A(x)\nabla u\cdot\nabla u)^{p/2}\,dx
+CL_S\int_{B_S}|\nabla u|^2\,dx.
\label{dq:uniform-fixed-coefficient}
\end{align}
To obtain the actual spatial difference quotient of $V$, add
\[
\frac{(A_h\nabla u\cdot\nabla u)^{\frac{p-2}{4}}\nabla u
 -(A(x)\nabla u\cdot\nabla u)^{\frac{p-2}{4}}\nabla u}{h}
\]
to the quotient in the first integral of
\eqref{dq:uniform-fixed-coefficient}. By \eqref{dq:matrix-V-variation},
this correction is at most $CK|\nabla u|^{p/2}$ in absolute value.
Also $|V|^2\le\lambda^{-1}(A(x)\nabla u\cdot\nabla u)^{p/2}$.
Thus the local difference-quotient criterion
\cite[Lemma~7.24]{gilbarg-trudinger}, applied in every direction, gives
$V\in W^{1,2}(B_\rho;\R^n)$ with \eqref{dq:V-estimate}.

For Hessian recovery, let
$\Omega_\delta=B_\rho\cap\{A(x)\nabla u\cdot\nabla u>\delta\}$.
On this open set,
\[
 V\cdot A(x)V=(A(x)\nabla u\cdot\nabla u)^{p/2},
\]
so
\begin{equation}
\nabla u=(V\cdot A(x)V)^{-\frac{p-2}{2p}}V.
\label{dq:invert-V}
\end{equation}
On each $U\Subset\Omega_\delta$, the right-hand side is a Lipschitz
function of $(A,V)$ near its bounded range, since
$V\cdot A(x)V\ge\delta^{p/2}$. The Sobolev chain and product rules
\cite[Section~7.4]{gilbarg-trudinger} therefore give
$u\in W^{2,2}_{\rm loc}(\Omega_\delta)$.
Put $\beta=(p-2)/4$. For $\xi\ne0$, we differentiate
with respect to $\xi$ in the direction $z$ and take the scalar product
with $A(x)z$ to obtain
\begin{align*}
&D_\xi[(A(x)\xi\cdot\xi)^\beta\xi][z]\cdot A(x)z\\
&\quad=(A(x)\xi\cdot\xi)^\beta A(x)z\cdot z
 +2\beta(A(x)\xi\cdot\xi)^{\beta-1}(A(x)\xi\cdot z)^2\\
&\quad\ge\lambda(A(x)\xi\cdot\xi)^\beta|z|^2.
\end{align*}
Since $|A(x)z|\le\Lambda|z|$, it follows that
\begin{equation}
|D_\xi[(A(x)\xi\cdot\xi)^\beta\xi][z]|
\ge\frac\lambda\Lambda(A(x)\xi\cdot\xi)^\beta|z|.
\label{dq:gradient-derivative-coercive}
\end{equation}
The coefficient part of $\partial_kV$ is
\[
\beta(A(x)\nabla u\cdot\nabla u)^{\beta-1}
((\partial_kA(x))\nabla u\cdot\nabla u)\nabla u,
\]
whose norm is at most $C|DA|(A(x)\nabla u\cdot\nabla u)^{p/4}$
by ellipticity.
Applying \eqref{dq:gradient-derivative-coercive} to
$z=\partial_k\nabla u$ in the chain rule therefore gives
\[
 (A(x)\nabla u\cdot\nabla u)^{\frac{p-2}{4}}
 |\partial_k\nabla u|
 \le C|\partial_kV|
 +C|DA|(A(x)\nabla u\cdot\nabla u)^{p/4}.
\]
Squaring and summing over $k$ gives
\begin{equation}
(A(x)\nabla u\cdot\nabla u)^{\frac{p-2}{2}}|D^2u|^2
\le C|DV|^2+CK^2(A(x)\nabla u\cdot\nabla u)^{p/2}
\quad\text{a.e. on }\Omega_\delta.
\label{dq:anisotropic-Hessian-bound}
\end{equation}
The local weak Hessians agree on overlaps. Integrating and letting
$\delta\downarrow0$ proves \eqref{dq:weighted-Hessian} by monotone
convergence and \eqref{dq:V-estimate}.

\medskip\noindent\textbf{The endpoint $p=2$.}
Here $V=\nabla u$. With the same cutoff and translations, the fixed-matrix
term in \eqref{dq:eq:exactidentity} is at least
$\lambda\int_{B_R}\zeta^2|\delta_h\nabla u|^2\,dx$.
Apply Cauchy--Schwarz and Young to the gradient-cutoff and
coefficient-gradient terms, absorbing each with coefficient
$\varepsilon<\lambda/2$. Their remainders are bounded by
$C_\varepsilon(\|\nabla\zeta\|_\infty^2+K^2)
\int_{B_S}|\nabla u|^2\,dx$.
The coefficient-cutoff and reaction estimates above remain valid, so
\[
\int_{B_\rho}|\delta_h\nabla u|^2\,dx
\le C\big((R-\rho)^{-2}+K^2+L_S\big)
\int_{B_S}|\nabla u|^2\,dx.
\]
The difference-quotient criterion gives $u\in W^{2,2}(B_\rho)$.
All zero powers equal $1$, including on the critical set.

\medskip\noindent\textbf{Case 2: $1<p<2$.}
Choose $\zeta\in C_c^\infty(B_R)$, $0\le\zeta\le1$, with
$\zeta=1$ on $B_{\rho_1}$ and
\begin{equation}
 |\nabla\zeta|\le\frac{C}{R-\rho_1},\qquad
 |D^2\zeta|\le\frac{C}{(R-\rho_1)^2}.
 \label{dq:eq:cutoff}
\end{equation}
Throughout this case, $0<|h|<h_0<\min\{S-R,\rho_1-\rho\}$.
The averaged-flux identity used in \cite[Lemma~4.5 and the proof of
Theorem~4.6]{lindqvist} reads, for $g\in L^{p'}(B_S)$ and
$\psi\in W_0^{1,p}(B_R)$,
\begin{equation}
 \int_{B_R}\psi\,\frac{g(x+he_k)-g(x)}h\,dx
 =-\int_{B_R}\partial_k\psi\int_0^1g(x+the_k)\,dt\,dx.
 \label{dq:eq:case2-IBP}
\end{equation}
The average has the difference quotient as its $L^{p'}$ weak
$k$-derivative, so density gives the identity. Define the vector average
\begin{align}
 \overline G_h(x):={}&\int_0^1
 \big(A(x+the_k)\nabla u(x+the_k)\cdot\nabla u(x+the_k)\big)^{\frac{p-2}{2}}\notag\\[-1mm]
 &\hspace{16mm}\times A(x+the_k)\nabla u(x+the_k)\,dt.
 \label{dq:average-flux-definition}
\end{align}
Applying \eqref{dq:eq:case2-IBP} componentwise with
$\psi_i=\zeta\,\delta_hu\,\partial_i\zeta$ rewrites the cutoff term in
\eqref{dq:eq:exactidentity} as
\begin{equation}
 2\int_{B_R}\overline G_h\cdot
 \big[\zeta\,\partial_k(\delta_hu)\,\nabla\zeta
 +\delta_hu\,\partial_k\zeta\,\nabla\zeta
 +\zeta\,\delta_hu\,\nabla(\partial_k\zeta)\big]\,dx.
 \label{dq:eq:ibpexplicit}
\end{equation}
Here $\partial_k(\delta_hu)=\delta_h(\partial_ku)$ uses only first
derivatives of $u$. Ellipticity and Jensen's inequality give
\begin{equation}
\begin{aligned}
 |\overline G_h(x)|^{p'}&\le C\int_0^1|\nabla u(x+the_k)|^p\,dt,\\
 \|\overline G_h\|_{L^{p'}(B_R)}
 &\le C\|\nabla u\|_{L^p(B_S)}^{p-1}.
\end{aligned}
 \label{dq:eq:averageflux}
\end{equation}

\medskip\noindent\emph{Coercivity and absorption.}
Split the left-hand flux difference in \eqref{dq:eq:exactidentity}
at \[
(A_h\nabla u\cdot\nabla u)^{(p-2)/2}A_h\nabla u.
\]
By Lemma~\ref{dq:lem:monotonicity}, its gradient part satisfies
\begin{align}
 &\frac{(A_h\nabla u_h\cdot\nabla u_h)^{\frac{p-2}{2}}A_h\nabla u_h
 -(A_h\nabla u\cdot\nabla u)^{\frac{p-2}{2}}A_h\nabla u}{h}
 \cdot\delta_h\nabla u\notag\\
 &\qquad\ge c
 \big(\mu^2+A_h\nabla u_h\cdot\nabla u_h+A_h\nabla u\cdot\nabla u\big)^{\frac{p-2}{2}}
 |A_h^{1/2}\delta_h\nabla u|^2
 \label{dq:eq:mainlowerexplicit}
\end{align}
for every $\mu>0$. The parameter $\mu$ only decreases the negative-power
weight; the flux and equation remain unchanged. For the coefficient part,
\eqref{dq:coefficient-vector-bound} and ellipticity bound its absolute
value by the left-hand side below. Weighted Young's inequality, using
the density in \eqref{dq:eq:mainlowerexplicit}, gives
\begin{align*}
 &CK\int_{B_R}\zeta^2|\nabla u|^{p-1}
 |A_h^{1/2}\delta_h\nabla u|\,dx\\
 &\quad\le\varepsilon\int_{B_R}\zeta^2
 \big(\mu^2+A_h\nabla u_h\cdot\nabla u_h+A_h\nabla u\cdot\nabla u\big)^{\frac{p-2}{2}}
 |A_h^{1/2}\delta_h\nabla u|^2\,dx\\
 &\qquad+C_\varepsilon K^2\int_{B_R}\zeta^2|\nabla u|^{2p-2}
 \big(\mu^2+A_h\nabla u_h\cdot\nabla u_h+A_h\nabla u\cdot\nabla u\big)^{\frac{2-p}{2}}\,dx.
\end{align*}
Since the base is at least $\lambda|\nabla u|^2$ and
$(p-1)+(2-p)/2=p/2$, the remainder is at most $C_\varepsilon K^2$
times the integral of its $p/2$-power against $\zeta^2$. Moreover,
\begin{equation}
 \big(\mu^2+A_h\nabla u_h\cdot\nabla u_h+A_h\nabla u\cdot\nabla u\big)^{p/2}
 \le C\big(\mu^p+|\nabla u_h|^p+|\nabla u|^p\big).
 \label{dq:eq:weightp}
\end{equation}
Translation therefore bounds the coefficient remainder by
$C_\varepsilon K^2\big(\mu^p|B_R|
+\int_{B_S}(A(x)\nabla u\cdot\nabla u)^{p/2}\,dx\big)$.

For the first term of \eqref{dq:eq:ibpexplicit}, use
$|\partial_k(\delta_hu)|\le\lambda^{-1/2}|A_h^{1/2}\delta_h\nabla u|$
and weighted Young's inequality. Its remainder is controlled by
\eqref{dq:eq:youngspecial}, with exponents $p/(2-p)$ and $p/[2(p-1)]$:
\begin{align}
 &\big(\mu^2+A_h\nabla u_h\cdot\nabla u_h+A_h\nabla u\cdot\nabla u\big)^{\frac{2-p}{2}}
 |\overline G_h|^2\notag\\*
 &\qquad\le C
 \big(\mu^2+A_h\nabla u_h\cdot\nabla u_h+A_h\nabla u\cdot\nabla u\big)^{p/2}
 +C|\overline G_h|^{p'}.
 \label{dq:eq:firstIBPpower}
\end{align}
Together with \eqref{dq:eq:cutoff}, \eqref{dq:eq:averageflux}, and
\eqref{dq:eq:weightp}, this yields
\begin{align}
 &\left|2\int_{B_R}\zeta\,\partial_k(\delta_hu)\,
 \overline G_h\cdot\nabla\zeta\,dx\right|\notag\\*
 &\quad\le\varepsilon\int_{B_R}\zeta^2
 \big(\mu^2+A_h\nabla u_h\cdot\nabla u_h+A_h\nabla u\cdot\nabla u\big)^{\frac{p-2}{2}}
 |A_h^{1/2}\delta_h\nabla u|^2\,dx\notag\\*
 &\qquad+\frac{C_\varepsilon}{(R-\rho_1)^2}
 \left(\mu^p|B_R|+\int_{B_S}(A(x)\nabla u\cdot\nabla u)^{p/2}\,dx\right).
 \label{dq:eq:firstIBPestimate}
\end{align}
H\"older's inequality, \eqref{dq:dq-Lp}, and \eqref{dq:eq:averageflux}
bound the sum of the absolute values of the remaining two terms in
\eqref{dq:eq:ibpexplicit} by
\[
 \frac{C}{(R-\rho_1)^2}
 \|\delta_hu\|_{L^p(B_R)}\|\overline G_h\|_{L^{p'}(B_R)}
 \le\frac{C}{(R-\rho_1)^2}
 \int_{B_S}(A(x)\nabla u\cdot\nabla u)^{p/2}\,dx.
\]
Finally, \eqref{dq:f-local-Lipschitz}
and $|\delta_hu|^2\le M^{2-p}|\delta_hu|^p$ give
\begin{equation}
 \left|\int_{B_R}\zeta^2\frac{f(u_h)-f(u)}h\,\delta_hu\,dx\right|
 \le CL_SM^{2-p}\int_{B_S}(A(x)\nabla u\cdot\nabla u)^{p/2}\,dx.
 \label{dq:eq:reaction}
\end{equation}
Taking $2\varepsilon<c$ absorbs both weighted quadratic errors, with
$\zeta^2$ retained, and \eqref{dq:eq:exactidentity} becomes
\begin{align}
 &\int_{B_R}\zeta^2
 \big(\mu^2+A_h\nabla u_h\cdot\nabla u_h+A_h\nabla u\cdot\nabla u\big)^{\frac{p-2}{2}}
 |A_h^{1/2}\delta_h\nabla u|^2\,dx\notag\\
 &\quad\le C\left(\frac1{(R-\rho_1)^2}+K^2\right)
 \left(\mu^p|B_R|+\int_{B_S}(A(x)\nabla u\cdot\nabla u)^{p/2}\,dx\right)
 \notag\\
 &\qquad+CL_SM^{2-p}\int_{B_S}(A(x)\nabla u\cdot\nabla u)^{p/2}\,dx,
 \label{dq:eq:muestimate}
\end{align}
with constants independent of $h$ and $\mu$.

\medskip\noindent\emph{Removing the weight regularization and recovering the Hessian.}
For fixed $h$, let $\mu\downarrow0$. The densities increase and vanish
identically where $\nabla u_h=\nabla u=0$; their limit there is defined
as zero. Monotone convergence and $\zeta=1$ on $B_{\rho_1}$ give
\begin{equation}
\begin{aligned}
 &\int_{B_{\rho_1}}
 \big(A_h\nabla u_h\cdot\nabla u_h+A_h\nabla u\cdot\nabla u\big)^{\frac{p-2}{2}}
 |A_h^{1/2}\delta_h\nabla u|^2\,dx\\
 &\quad\le C\left(\frac1{(R-\rho_1)^2}+K^2+L_SM^{2-p}\right)
 \int_{B_S}(A(x)\nabla u\cdot\nabla u)^{p/2}\,dx.
\end{aligned}
\label{dq:eq:weighteddq}
\end{equation}
At a zero of the weight's base, ellipticity also gives
$\delta_h\nabla u=0$. On the positive-base set, factor
\begin{align*}
 |A_h^{1/2}\delta_h\nabla u|^p
 &=\left[
 \big(A_h\nabla u_h\cdot\nabla u_h+A_h\nabla u\cdot\nabla u\big)^{\frac{p-2}{2}}
 |A_h^{1/2}\delta_h\nabla u|^2\right]^{p/2}\\*
 &\quad\times
 \big(A_h\nabla u_h\cdot\nabla u_h+A_h\nabla u\cdot\nabla u\big)^{\frac{p(2-p)}4}.
\end{align*}
The same zero convention extends this identity to the double-zero set.
Ellipticity and H\"older's inequality with exponents $2/p$ and $2/(2-p)$ yield
\begin{align}
 \int_{B_{\rho_1}}|\delta_h\nabla u|^p\,dx
 &\le C\left[\int_{B_{\rho_1}}
 \big(A_h\nabla u_h\cdot\nabla u_h+A_h\nabla u\cdot\nabla u\big)^{\frac{p-2}{2}}
 |A_h^{1/2}\delta_h\nabla u|^2\,dx\right]^{p/2}\notag\\*
 &\quad\times\left[\int_{B_{\rho_1}}
 \big(A_h\nabla u_h\cdot\nabla u_h+A_h\nabla u\cdot\nabla u\big)^{p/2}\,dx
 \right]^{\frac{2-p}{2}}\notag\\*
 &\le C\left(\frac1{(R-\rho_1)^2}+K^2+L_SM^{2-p}\right)^{p/2}
 \int_{B_S}(A(x)\nabla u\cdot\nabla u)^{p/2}\,dx.
 \label{dq:eq:W2p}
\end{align}
The last inequality uses \eqref{dq:eq:weighteddq}, ellipticity and
translation. The difference-quotient criterion gives
$u\in W^{2,p}(B_{\rho_1})$.
The difference-quotient estimate and smooth approximation
\cite[Lemma~7.23 and Theorem~7.9]{gilbarg-trudinger}
yield, for each $k$,
\[
 \delta_h^{(k)}\nabla u\longrightarrow D_k\nabla u
 \quad\text{strongly in }L^p(B_\rho;\mathbb R^n),
\]
and almost everywhere along a subsequence.
Continuity of $A$, $\nabla u$, and the matrix square root then gives,
on $\{|\nabla u|>0\}$, the limiting density
\[
 2^{\frac{p-2}{2}}(A(x)\nabla u\cdot\nabla u)^{\frac{p-2}{2}}
 |A(x)^{1/2}D_k\nabla u|^2.
\]
Fatou's lemma on this fixed set, applied for each $k$ and summed, yields
\begin{equation}
\begin{aligned}
 &\int_{B_\rho\cap\{|\nabla u|>0\}}
 (A(x)\nabla u\cdot\nabla u)^{\frac{p-2}{2}}
 |A(x)^{1/2}D^2u|^2\,dx\\
 &\quad\le C\left(\frac1{(R-\rho_1)^2}+K^2+L_SM^{2-p}\right)
 \int_{B_S}(A(x)\nabla u\cdot\nabla u)^{p/2}\,dx.
\end{aligned}
\label{dq:eq:onesided}
\end{equation}
Since $|A(x)^{1/2}D^2u|^2\ge\lambda|D^2u|^2$, this proves
\eqref{dq:eq:mainestimate}.
\end{proof}

\begin{cor}[]
\label{cor:dq-weighted-sobolev}
Assume the hypotheses of Theorem~\ref{dq:main-theorem} and let $p\ge2$.  If
$U\Subset\Omega$ is a bounded Lipschitz domain and
\[
 \sigma=|\nabla u|^{p-2},
\]
then every component of $\nabla u$ belongs to the weighted space
$W^{1,2}_\sigma(U)$, defined as the completion of those functions
in $C^1(U)$ for which the following norm is finite:
\[
 \|h\|_{L^2(U)}+
 \left(\int_U\sigma|\nabla h|^2\,dx\right)^{1/2}.
\]
In particular,
\[
 \nabla u\in W^{1,1}_\sigma(U;\mathbb R^n),
\]
so the weighted Sobolev hypothesis in Lemma~\ref{lem:cube-interpolation} is
satisfied on every cube compactly contained in $\Omega$.
\end{cor}

\begin{proof}
For $p=2$, apply Theorem~\ref{dq:main-theorem}. For $p>2$, use the
noncritical Hessian $D^2u$, extended by zero on $\{\nabla u=0\}$,
in the weighted expressions below. Ellipticity and
\eqref{dq:weighted-Hessian} give
$\int_U|\nabla u|^{p-2}|D^2u|^2\,dx<\infty$.
Choose $\chi\in C^\infty([0,\infty))$ with $0\le\chi\le1$,
$\chi=0$ on $[0,1]$ and $\chi=1$ on $[2,\infty)$, and set
\[
 \chi_\tau(x)=\chi\left(\frac{|\nabla u(x)|}{\tau}\right).
\]
The recovery argument gives $\nabla u\in W^{1,2}_{\rm loc}$ on the
noncritical region. For fixed $\tau>0$, $\chi_\tau\nabla u$ vanishes on
$U\cap\{|\nabla u|\le\tau\}$. The compact set
$\overline U\cap\{|\nabla u|\ge\tau\}$ lies in the noncritical region.
A finite covering of this set gives local chain-rule derivatives,
which agree with zero on overlaps with $U\cap\{|\nabla u|<\tau\}$.
These derivatives therefore patch to give the weak derivative on $U$. For
$E_\tau=\{\tau<|\nabla u|<2\tau\}$, this gives
\[
 D(\chi_\tau\nabla u)=\chi_\tau D^2u+\nabla u\otimes D\chi_\tau,
 \qquad |D\chi_\tau|\le\frac C\tau|D^2u|\mathbf1_{E_\tau}.
\]
Since $|\nabla u|\le2\tau$ on $E_\tau$, one has
$|D(\chi_\tau\nabla u)|\le C|D^2u|\mathbf1_{\{|\nabla u|>\tau\}}$ and
\[
 \int_U|D(\chi_\tau\nabla u)|^2\,dx
 \le C\tau^{2-p}\int_U|\nabla u|^{p-2}|D^2u|^2\,dx<\infty.
\]
Together with the boundedness of $\chi_\tau\nabla u$, this proves
$\chi_\tau\nabla u\in W^{1,2}(U;\mathbb R^n)$.
The same derivative formula also gives
\[
 \int_U|\nabla u|^{p-2}|D(\chi_\tau\nabla u)-D^2u|^2\,dx
 \le C\int_{U\cap\{0<|\nabla u|<2\tau\}}
 |\nabla u|^{p-2}|D^2u|^2\,dx\longrightarrow0,
\]
and $\|\chi_\tau\nabla u-\nabla u\|_{L^2(U)}\le2\tau|U|^{1/2}\to0$.
Since $\sigma=|\nabla u|^{p-2}$ is bounded on $U$, smooth $W^{1,2}(U)$
approximation of each $\chi_\tau\nabla u$
\cite[Theorem~7.9]{gilbarg-trudinger} also converges in the weighted norm.
A diagonal choice, combined with the convergence above, gives
$\nabla u\in W^{1,2}_\sigma(U;\mathbb R^n)$.
Finally, $\sigma\in L^1(U)$ and Cauchy--Schwarz gives
\[
\int_U\sigma|Dh|\le(\int_U\sigma|Dh|^2)^{1/2}(\int_U\sigma)^{1/2}.
\]
Applying this to the smooth approximation errors and using
$L^2(U)\hookrightarrow L^1(U)$ gives
$\nabla u\in W^{1,1}_\sigma(U;\mathbb R^n)$.
\end{proof}

\Needspace{12\baselineskip}
\section*{Acknowledgments}

The author would like to thank Professor Fa Peng for many helpful discussions.
The author is supported by the National Key R\&D Program of China
(No.~2025YFA1018400) and the National Natural Science Foundation of China
(No.~12571209).

\begingroup
\small
\raggedright
\urlstyle{same}

\endgroup

\bigskip
\noindent
\begin{minipage}{\textwidth}
\small
Weican Luo\\
School of Mathematical Sciences, Beihang University, Changping District,
Shahe Higher Education Park, South Third Street No.~9, Beijing 102206, PR China\\
\textit{E-mail:} \href{mailto:weicanluo@buaa.edu.cn}{\nolinkurl{weicanluo@buaa.edu.cn}}
\end{minipage}


\begin{thebibliography}{23}
\addcontentsline{toc}{section}{References}

\bibitem{antonini-ciraolo-farina}
C.~A.~Antonini, G.~Ciraolo, and A.~Farina,
Interior regularity results for inhomogeneous anisotropic quasilinear equations,
\emph{Math. Ann.} \textbf{387} (2023), 1745--1776.
\href{https://doi.org/10.1007/s00208-022-02500-x}{\nolinkurl{doi:10.1007/s00208-022-02500-x}}.

\bibitem{araujo-zhang}
D.~J.~Ara\'ujo and L.~Zhang,
Optimal $C^{1,\alpha}$ estimates for a class of elliptic quasilinear equations,
\emph{Commun. Contemp. Math.} \textbf{22} (2020), Paper No.~1950014.
\href{https://doi.org/10.1142/S0219199719500147}{\nolinkurl{doi:10.1142/S0219199719500147}}.

\bibitem{brezis-vazquez}
H.~Brezis and J.~L.~V\'azquez,
Blow-up solutions of some nonlinear elliptic problems,
\emph{Rev. Mat. Univ. Complut. Madrid} \textbf{10} (1997), 443--469.
\href{https://doi.org/10.5209/rev_REMA.1997.v10.n2.17459}{\nolinkurl{doi:10.5209/rev_REMA.1997.v10.n2.17459}}.

\bibitem{cabre-low-dim}
X. Cabr\'e,
Regularity of minimizers of semilinear elliptic problems up to dimension 4,
\emph{Comm. Pure Appl. Math.} \textbf{63} (2010), 1362--1380.
\href{https://doi.org/10.1002/cpa.20327}{\nolinkurl{doi:10.1002/cpa.20327}}.

\bibitem{cabre-capella-sanchon}
X. Cabr\'e, A. Capella, and M. Sanch\'on,
Regularity of radial minimizers of reaction equations involving the $p$-Laplacian,
\emph{Calc. Var. Partial Differential Equations} \textbf{34} (2009), 475--494.
\href{https://doi.org/10.1007/s00526-008-0192-3}{\nolinkurl{doi:10.1007/s00526-008-0192-3}}.

\bibitem{cfrs}
X. Cabr\'e, A. Figalli, X. Ros-Oton, and J. Serra,
Stable solutions to semilinear elliptic equations are smooth up to dimension 9,
\emph{Acta Math.} \textbf{224} (2020), 187--252.
\href{https://doi.org/10.4310/ACTA.2020.v224.n2.a1}{\nolinkurl{doi:10.4310/ACTA.2020.v224.n2.a1}}.

\bibitem{cms}
X. Cabr\'e, P. Miraglio, and M. Sanch\'on,
Optimal regularity of stable solutions to nonlinear equations involving the $p$-Laplacian,
\emph{Adv. Calc. Var.} \textbf{15} (2022), 749--785.
\href{https://doi.org/10.1515/acv-2020-0055}{\nolinkurl{doi:10.1515/acv-2020-0055}}.
Corrected version: \href{https://arxiv.org/abs/2006.01445v3}{\nolinkurl{arXiv:2006.01445v3}} (2022).

\bibitem{cabre-sanchon}
X. Cabr\'e and M. Sanch\'on,
Semi-stable and extremal solutions of reaction equations involving the $p$-Laplacian,
\emph{Commun. Pure Appl. Anal.} \textbf{6} (2007), 43--67.
\href{https://doi.org/10.3934/cpaa.2007.6.43}{\nolinkurl{doi:10.3934/cpaa.2007.6.43}}.

\bibitem{castorina-riey-sciunzi}
D.~Castorina, G.~Riey, and B.~Sciunzi,
Hopf Lemma and regularity results for quasilinear anisotropic elliptic equations,
\emph{Calc. Var. Partial Differential Equations} \textbf{58} (2019), Paper No.~95.
\href{https://doi.org/10.1007/s00526-019-1528-x}{\nolinkurl{doi:10.1007/s00526-019-1528-x}}.

\bibitem{castorina-sanchon}
D. Castorina and M. Sanch\'on,
Regularity of stable solutions of $p$-Laplace equations through geometric Sobolev type inequalities,
\emph{J. Eur. Math. Soc.} \textbf{17} (2015), 2949--2975.
\href{https://doi.org/10.4171/JEMS/576}{\nolinkurl{doi:10.4171/JEMS/576}}.

\bibitem{crandall-rabinowitz}
M.~G. Crandall and P.~H. Rabinowitz,
Some continuation and variational methods for positive solutions of nonlinear elliptic eigenvalue problems,
\emph{Arch. Rational Mech. Anal.} \textbf{58} (1975), 207--218.
\href{https://doi.org/10.1007/BF00280741}{\nolinkurl{doi:10.1007/BF00280741}}.

\bibitem{damascelli-sciunzi}
L.~Damascelli and B.~Sciunzi,
Regularity, monotonicity and symmetry of positive solutions of $m$-Laplace equations,
\emph{J. Differential Equations} \textbf{206} (2004), 483--515.
\href{https://doi.org/10.1016/j.jde.2004.05.012}{\nolinkurl{doi:10.1016/j.jde.2004.05.012}}.

\bibitem{dong-zhu}
H.~Dong and H.~Zhu,
Gradient estimates for singular $p$-Laplace type equations with measure data,
\emph{J. Eur. Math. Soc.} \textbf{26} (2024), 3939--3985.
\href{https://doi.org/10.4171/JEMS/1400}{\nolinkurl{doi:10.4171/JEMS/1400}}.

\bibitem{erneta}
I.~U. Erneta,
Stable solutions to semilinear elliptic equations for operators with variable coefficients,
\emph{Commun. Pure Appl. Anal.} \textbf{22} (2023), 530--571.
\href{https://doi.org/10.3934/cpaa.2022162}{\nolinkurl{doi:10.3934/cpaa.2022162}}.

\bibitem{farina-sciunzi-valdinoci}
A.~Farina, B.~Sciunzi, and E.~Valdinoci,
Bernstein and De Giorgi type problems: new results via a geometric approach,
\emph{Ann. Scuola Norm. Sup. Pisa Cl. Sci.} (5) \textbf{7} (2008), 741--791.
\href{https://doi.org/10.2422/2036-2145.2008.4.06}{\nolinkurl{doi:10.2422/2036-2145.2008.4.06}}.

\bibitem{gilbarg-trudinger}
D.~Gilbarg and N.~S.~Trudinger,
\emph{Elliptic Partial Differential Equations of Second Order},
reprint of the 1998 edition, Classics in Mathematics, Springer, Berlin, 2001.
\href{https://doi.org/10.1007/978-3-642-61798-0}{\nolinkurl{doi:10.1007/978-3-642-61798-0}}.

\bibitem{hkm}
J.~Heinonen, T.~Kilpel\"ainen, and O.~Martio,
\emph{Nonlinear Potential Theory of Degenerate Elliptic Equations},
Dover Publications, Mineola, NY, 2006.
\url{https://users.jyu.fi/~terok/NLPT/NLPT.html}.

\bibitem{lindqvist}
P.~Lindqvist,
\emph{Notes on the $p$-Laplace Equation}, second edition,
Report~161, University of Jyv\"askyl\"a, Department of Mathematics and Statistics, 2017.
\url{https://jyx.jyu.fi/handle/123456789/54658}.

\bibitem{miao-peng-zhou}
Q.~Miao, F.~Peng, and Y.~Zhou,
A global second-order Sobolev regularity for $p$-Laplacian type equations with variable coefficients in bounded domains,
\emph{Calc. Var. Partial Differential Equations} \textbf{62} (2023), Paper No.~191.
\href{https://doi.org/10.1007/s00526-023-02538-y}{\nolinkurl{doi:10.1007/s00526-023-02538-y}}.

\bibitem{nedev}
G. Nedev,
Regularity of the extremal solution of semilinear elliptic equations,
\emph{C. R. Acad. Sci. Paris S\'er. I Math.} \textbf{330} (2000), 997--1002.
\href{https://doi.org/10.1016/S0764-4442(00)00289-5}{\nolinkurl{doi:10.1016/S0764-4442(00)00289-5}}.

\bibitem{simon-scaling}
L.~Simon,
Schauder estimates by scaling,
\emph{Calc. Var. Partial Differential Equations} \textbf{5} (1997), 391--407.
\href{https://doi.org/10.1007/s005260050072}{\nolinkurl{doi:10.1007/s005260050072}}.

\bibitem{tolksdorf}
P.~Tolksdorf,
Regularity for a more general class of quasilinear elliptic equations,
\emph{J. Differential Equations} \textbf{51} (1984), 126--150.
\href{https://doi.org/10.1016/0022-0396(84)90105-0}{\nolinkurl{doi:10.1016/0022-0396(84)90105-0}}.

\bibitem{uhlenbeck}
K.~Uhlenbeck,
Regularity for a class of non-linear elliptic systems,
\emph{Acta Math.} \textbf{138} (1977), 219--240.
\href{https://doi.org/10.1007/BF02392316}{\nolinkurl{doi:10.1007/BF02392316}}.

\end{thebibliography}
\end{document}